\documentclass[11pt,leqno]{article}
\usepackage{amsmath,amssymb,amscd,latexsym}

\usepackage[all]{xy}

\newcommand{\C}{{\mathbb{C}}}
\newcommand{\F}{{\mathbb{F}}}
\newcommand{\oF}{\overline{\F}}
\newcommand{\Ge}{\mathbb{G}}
\newcommand{\Pa}{{\mathbb{P}}}
\newcommand{\Q}{{\mathbb{Q}}}
\newcommand{\oQ}{\overline{\Q}}
\newcommand{\Z}{{\mathbb{Z}}}
\newcommand{\oZ}{\overline{\Z}}
\newcommand{\bC}{\mathbf{C}}

\newcommand{\bV}{\mathbf{V}}

\newcommand{\tC}{\tilde{C}}
\newcommand{\tD}{\tilde{D}}

\newcommand{\Alb}{\mathrm{Alb}}
\newcommand{\adj}{\mathrm{adj}}
\newcommand{\Ann}{\mathrm{Ann}}
\newcommand{\cont}{\mathrm{cont}}
\newcommand{\et}{\mathrm{\acute{e}t}}

\newcommand{\orho}{\overline{\rho}}
\newcommand{\os}{\overline{s}}
\newcommand{\good}{\mathrm{good}}
\newcommand{\id}{\mathrm{id}}
\renewcommand{\mod}{\;\mathrm{mod}\;}
\newcommand{\Mod}{\mathbf{Mod}}
\newcommand{\Mor}{\mathrm{Mor}}
\newcommand{\op}{\mathrm{op}}
\newcommand{\rank}{\mathrm{rank}}
\newcommand{\red}{\mathrm{red}}
\newcommand{\rk}{\mathrm{rk}\,}
\newcommand{\sing}{\mathrm{sing}}
\newcommand{\spec}{\mathrm{spec}\,}
\newcommand{\supp}{\mathrm{supp}\,}
\newcommand{\Aut}{\mathrm{Aut}}

\newcommand{\Ext}{\mathrm{Ext}}
\newcommand{\Fr}{\mathrm{Fr}}

\newcommand{\Gal}{\mathrm{Gal}}
\newcommand{\GL}{\mathrm{GL}\,}
\newcommand{\Hom}{\mathrm{Hom}}
\newcommand{\uHom}{\underline{\Hom}}
\newcommand{\Imm}{\mathrm{Im}\,}
\newcommand{\Iso}{\mathrm{Iso}\,}
\newcommand{\Ker}{\mathrm{Ker}\,}
\newcommand{\Lie}{\mathrm{Lie}\,}
\newcommand{\Ob}{\mathrm{Ob}\,}
\newcommand{\Pic}{\mathrm{Pic}}
\newcommand{\PGL}{\mathrm{PGL}}
\newcommand{\Rep}{\mathbf{Rep}\,}
\newcommand{\uVec}{\mathbf{Vec}}
\newcommand{\Ah}{{\mathcal A}}
\newcommand{\Bh}{\mathcal{B}}
\newcommand{\Ch}{{\mathcal C}}
\newcommand{\Eh}{{\mathcal E}}
\newcommand{\Fh}{{\mathcal F}}

\newcommand{\Ih}{{\mathcal I}}
\newcommand{\Jh}{\mathcal{J}}

\newcommand{\Lh}{{\mathcal L}}

\newcommand{\Nh}{{\mathcal N}}
\newcommand{\Oh}{{\mathcal O}}

\newcommand{\Sh}{{\mathcal S}}

\newcommand{\Yh}{\mathcal{Y}}
\newcommand{\Zh}{\mathcal{Z}}

\newcommand{\ea}{\mathfrak{a}}
\newcommand{\eb}{\mathfrak{b}}

\newcommand{\emm}{{\mathfrak{m}}}
\newcommand{\eo}{\mathfrak{o}}
\newcommand{\ep}{\mathfrak{p}}

\newcommand{\eB}{\mathfrak{B}}
\newcommand{\eE}{\mathcal{E}}

\newcommand{\eS}{\mathfrak{S}}
\newcommand{\eT}{\mathfrak{T}}
\newcommand{\eX}{{\mathfrak X}}
\newcommand{\oeX}{\overline{\eX}}
\newcommand{\eZ}{\mathfrak{Z}}

\newcommand{\oK}{\overline{K}}
\newcommand{\oX}{\overline{X}}
\newcommand{\oY}{\overline{Y}}

\newcommand{\tYh}{\tilde{\Yh}}
\newcommand{\oYh}{\overline{\Yh}}
\newcommand{\tgamma}{\tilde{\gamma}}
\newcommand{\tpi}{\tilde{\pi}}
\newcommand{\tvarphi}{\tilde{\varphi}}

\newcommand{\trho}{\tilde{\rho}}

\newcommand{\tX}{\tilde{X}}

\newcommand{\ty}{\tilde{y}}
\newcommand{\tY}{\tilde{Y}}
\newcommand{\tZ}{\tilde{Z}}
\newcommand{\tz}{\tilde{z}}
\newcommand{\hsigma}{\,\!^{\sigma}\!}
\newcommand{\htau}{\,\!^{\tau}}
\newcommand{\ohne}{\setminus}
\newcommand{\silo}{\stackrel{\sim}{\longrightarrow}}
\newcommand{\tei}{\, | \,}
\newcommand{\hullet}{\raisebox{0.03cm}{$\scriptstyle \bullet$}}

\newcommand{\verk}{\raisebox{0.03cm}{\mbox{\scriptsize $\,\circ\,$}}}
\newtheorem{theorem}{Theorem}
\newtheorem{lemma}[theorem]{Lemma}
\newtheorem{prop}[theorem]{Proposition}
\newtheorem{defn}[theorem]{Definition}
\newtheorem{cor}[theorem]{Corollary}

\newenvironment{rem}{\noindent {\bf Remark}}{}
\newenvironment{rems}{\noindent {\bf Remarks}}{}
\newenvironment{definition}{\noindent {\bf Definition}\it}{}
\newenvironment{theoremon}{\noindent {\bf Theorem}\it}{}
\newenvironment{proofof}{\noindent {\bf Proof of}}{\mbox{}\hspace*{\fill}$\Box$}
\newenvironment{proof}{\noindent {\bf Proof}}{\mbox{}\hspace*{\fill}$\Box$}
\begin{document}
\title{Vector bundles on $p$-adic curves and parallel transport}
\author{Christopher Deninger \and Annette Werner}
\date{}
\maketitle
\begin{abstract}
We define functorial isomorphisms of parallel transport along \'etale paths for a class of vector bundles on a $p$-adic curve. All bundles of degree zero whose reduction is strongly semistable belong to this class. In particular, they give rise to representations of the algebraic fundamental group of the curve. This may be viewed as a partial analogue of the classical Narasimhan--Seshadri theory of vector bundles on compact Riemann surfaces.

\centerline{MSC (2000) 14H60, 14H30, 11G20}
\end{abstract}

\section*{Introduction}
On a compact Riemann surface every finite dimensional complex representation of the fundamental group gives rise to a flat vector bundle and hence to a holomorphic vector bundle. By a theorem of Weil, one obtains precisely the holomorphic bundles whose indecomposable components have degree zero \cite{W}. It was proved by Narasimhan and Seshadri \cite{Na-Se} that unitary representations give rise to polystable bundles of degree zero. Moreover, every stable bundle of degree zero comes from an irreducible unitary representation. 

The present paper establishes a partial $p$-adic analogue of this theory, generalized to representations of the fundamental groupoid. The following is our main result. Recall that a vector bundle on a smooth projective curve over a field of characteristic $p$ is called strongly semistable if the pullbacks of $E$ by all non-negative powers of the absolute Frobenius morphism are semistable. Let $X$ be a smooth projective curve over $\oQ_p$ and let $\eo$ be the ring of integers in $\C_p$. A model $\eX$ of $X$ is a finitely presented flat and proper scheme over $\oZ_p$ with generic fibre $X$. The special fibre $\eX_k$ is then a union of projective curves over $k = \oF_p$. We say that a vector bundle $E$ on $X_{\C_p} = X \otimes \C_p$ has strongly semistable reduction of degree zero if the following is true: $E$ can be extended to a vector bundle $\eE$ on $\eX_{\eo} = \eX \otimes \eo$ for some model $\eX$ of $X$ such that the pullback of the special fibre $\eE_k$ of $\eE$ to the normalization of each irreducible component of $\eX_k$ is strongly semistable of degree zero. We say that $E$ has potentially strongly semistable reduction of degree zero if there is a finite \'etale morphism $\alpha : Y \to X$ of smooth projective curves such that $\alpha^* E$ has strongly semistable reduction of degree zero.

\begin{theoremon}
  Let $E$ be a vector bundle on $X_{\C_p}$ with potentially strongly semistable reduction of degree zero. Then there are functorial isomorphisms of ``parallel transport'' along \'etale paths between the fibres of $E_{\C_p}$ on $X_{\C_p}$. In particular one obtains a representation $\rho_{E,x}$ of $\pi_1 (X,x)$ on $E_x$ for every point $x$ in $X (\C_p)$. The parallel transport is compatible with tensor products, duals, internal homs, pullbacks and Galois conjugation.
\end{theoremon}

The theorem applies in particular to line bundles of degree zero on $X_{\C_p}$. In this case the $p$-part of the corresponding character of $\pi_1 (X,x)$ was already constructed by Tate using Cartier duality for the $p$-divisible group of the abelian scheme $\Pic^0_{\eX / \oZ_p}$ cf. \cite{Ta} \S\,4 and \cite{De-We2}. His method does not extend to bundles of higher rank.

Let us now discuss the contents of the paper in more detail. Afterwards we can sketch the proof of the theorem.

In the first section we investigate the category $\Sh_{\eX,D}$ consisting of finitely presented proper $\oZ_p$-morphisms $\pi : \Yh \to \eX$ whose generic fibre is a finite covering of $Y$ which is \'etale outside of a divisor $D$ on $X$. The important point is that for given $\pi$ in $\Sh_{\eX,D}$ there is an object $\pi' : \Yh' \to \eX$ in $\Sh_{\eX,D}$ lying over $\pi$ with better properties, e.g. cohomologically flat of dimension zero or even semistable. We also construct certain coverings $\pi$ using the theory of the Picard functor which are used several times.

In the second section we define and investigate categories $\eB_{X_{\C_p} , D}$ and $\eB^{\sharp}_{X_{\C_p} ,D}$ involving a divisor $D$ on $X$ and also an analogous category $\eB_{\eX_{\eo} , D}$ for a fixed model $\eX$ of $X$. These are defined as follows. The category $\eB_{\eX_{\eo} , D}$ consists of all vector bundles $\eE$ on $\eX_{\eo}$ such that for all $n \ge 1$ there is a covering $\pi$ in $\Sh_{\eX,D}$ with $\pi^* \eE$ trivial modulo $p^n$. In theorem \ref{cd_t16_neu} it is proved that for $\eE$ to lie in $\eB_{\eX_{\eo , D}}$ it suffices that $\pi^*_k \eE_k$ is trivial where $\pi_k$ is the special fibre of some $\pi$.

Next, $\eB_{X_{\C_p} ,D}$ consists of all bundles which are isomorphic to the generic fibre of a bundle $\eE$ in $\eB_{\eX_{\eo}, D}$ for some model $\eX$ of $X$. These categories are additive and stable under extensions. Finally, we define $\eB^{\sharp}_{X_{\C_p} , D}$ as the category of vector bundles on $X_{\C_p}$ whose pullback along $\alpha$ lies in $\eB_{Y_{\C_p} , \alpha^* D}$ for some finite morphism $\alpha : Y \to X$ between smooth projective curves which is \'etale over $X \ohne D$. We obtain an additive category which is closed under extensions and contains all line bundles of degree zero. All vector bundles in $\eB^{\sharp}$ are semistable of degree zero.

The third section is devoted to the definition and study of certain isomorphisms of parallel transport along \'etale paths in $U = X \ohne D$ for the bundles in the category $\eB^{\sharp}_{X_{\C_p}, D}$. In more technical terms, we construct an exact $\otimes$-functor $\rho$ from $\eB^{\sharp}_{X_{\C_p} , D}$ to the category of continuous representations of the \'etale fundamental groupoid $\Pi_1 (U)$ on $\C_p$-vector spaces. The basic idea is this: Consider a bundle $\eE$ in $\eB_{\eX_{\eo} , D}$ and for a given $n \ge 1$ let $\pi : \Yh \to \eX$ be an object of $\Sh_{\eX,D}$ such that $\pi^*_n \eE_n$ is a trivial bundle on $\Yh_n$. Here the index $n$ denotes reduction modulo $p^n$. Consider points $x$ and $x'$ in $X (\C_p) = \eX (\eo)$ and choose a point $y$ in $Y =\Yh_{\C_p}$ above $x$. For an \'etale path $\gamma$ from $x$ to $x'$ i.e. an isomorphism of fibre functors, let $\gamma y$ be the corresponding point above $x'$. For a ``good'' cover $\pi$ we have isomorphisms
\[
\eE_{x_n} \xleftarrow{\overset{y^*_n}{\sim}} \Gamma (\Yh_n , \pi^*_n \eE_n) \xrightarrow{\overset{(\gamma y)^*_n}{\sim}} \eE_{x'_n} \; .
\]
We define the parallel transport $\rho_{\eE} (\gamma) : \eE_x \silo \eE_{x'}$ as the projective limit of the maps $\rho_{\eE,n} (\gamma) = (\gamma y)^*_n \verk (y^*_n)^{-1}$. This parallel transport is then extended to $\eB_{X_{\C_p} ,D}$ and $\eB^{\sharp}_{X_{\C_p} , D}$. We also prove that the functor mapping a bundle $E$ in $\eB^{\sharp}_{X_{\C_p},D}$ to its fibre in a point $x \in U (\C_p)$ is faithful.

Using a Seifert--van Kampen theorem for \'etale groupoids we show that for a bundle $E$ which is in $\eB^{\sharp}$ for two disjoint divisors, one actually obtains a parallel transport along all \'etale paths in $X$.

The proof of the theorem above starts with a characterization of those vector bundles on a purely one-dimensional proper scheme over a finite field $\F_q$ whose pullback to the normalization of each irreducible component is strongly semistable of degree zero: These are exactly the bundles whose pullback by a finite surjective morphism to a purely one-dimensional proper $\F_q$-scheme becomes trivial. For vector bundles on smooth projective curves over finite fields this characterization is due to Lange and Stuhler \cite{LS}. Hence we have to lift finite covers in characteristic $p$ to characteristic zero. The main point here is to construct a morphism of models whose reduction factors over a given power of Frobenius. In fact our method allows us to construct two coverings $\pi$ in $\Sh_{\eX, D}$ and $\tpi$ in $\Sh_{\eX,\tD}$ for two disjoint divisors $D$ and $\tD$ such that $\pi^*_k \eE_k$ and $\tpi^*_k \eE_k$ are both trivial. By the above theory, one gets the parallel transport on all of $X_{\C_p}$. In the case of good reduction M. Raynaud has shown us a direct proof of this fact c.f. theorem \ref{cd_t20_neu}.

For Mumford curves, Faltings (\cite{Fa}) associates a vector bundle on $X$ to every $K$-rational representation of the Schottky group and proves that every semistable vector bundle of degree zero arises in this way. It was shown by Herz \cite{He} that his construction is compatible with ours. 

Recently Faltings has announced a $p$-adic version of non-abelian Hodge theory \cite{Fa2}. He proves an equivalence of categories between vector bundles on $X_{\C_p}$ endowed with a $p$-adic Higgs field and a certain category of ``generalized representations'' which contains the representations of $\pi_1 (X,x)$ as a full subcategory. His methods are different from ours. In particular Faltings uses his theory of almost \'etale extensions. The main open problem in Faltings approach is to characterize the Higgs bundles corresponding to actual representations of $\pi_1 (X,x)$. He shows that with zero Higgs field, line bundles of degree zero and their successive extensions come from $\pi_1 (X,x)$-representations and suggests that perhaps all semistable vector bundles of degree zero are obtained in this way. The main theorem of our paper shows that this is true if in addition the bundle has potentially strongly semistable reduction.

The present preprint improves and  replaces the second part of \cite{De-We}. The first part of \cite{De-We} will be published as \cite{De-We2}.

Finally we would like to draw the reader's attention to possibly related works of Berkovich \cite{Be} \S\,9 on $p$-adic integration, of Ogus and Vologodsky on non-abelian Hodge theory in characteristic $p$ and of Vologodsky \cite{Vo} on Hodge structures on fundamental groups.

{\bf Acknowledgements:} We are indebted to Qing Liu for his detailed answers to our questions about models of curves and for sending us the preprint \cite{Liu2}. We are also grateful to Holger Brenner for helpful discussions about vector bundles. We would like to thank Gerd Faltings for criticising our confinement to smooth models in \cite{De-We} and for making available to us his notes \cite{Fa2}. Moreover, we thank Siegfried Bosch, H\'el\`ene Esnault, Sylvain Maugeais, Niko Naumann, Mohamed Sa\"{\i}di, Peter Schneider, Matthias Strauch and Ulrich Stuhler for interesting discussions and suggestions. Our greatest debt is to Michel Raynaud for his careful reading of the manuscript and for suggesting several improvements of our theorems and their proofs. This concerns in particular theorems \ref{aw_t5}
and \ref{cd_t20_neu} and proposition \ref{cd_t14neu}.

\section{Categories of ``coverings''}

In this section we introduce simplified and generalized versions of the categories of coverings that were used in \cite{De-We} to define the $p$-adic representations attached to certain vector bundles.

In the following, a variety over a field $k$ is a geometrically irreducible and geometrically reduced separated scheme of finite type over $k$. A curve is a one-dimensional variety. Let $R$ be a  valuation ring with quotient field $Q$ of characteristic zero. For a smooth projective curve $X$ over $Q$ consider a model $\eX$ of $X$ over $R$ i.e. a finitely presented, flat and proper scheme over $\spec R$ together with an isomorphism $X = \eX \otimes_R Q$. For a divisor $D$ on $X$ we write $X \ohne D$ for $X \ohne \supp D$.

Consider the following category $\Sh_{\eX,D}$. Objects are finitely presented proper $R$-morphisms $\pi : \Yh \to \eX$ whose generic fibre $\pi_{Q} : \Yh_{Q} \to X$ is finite and such that 
\[
\pi_Q : \pi^{-1}_Q (X \ohne D) \to X \ohne D \quad \mbox{is \'etale} \; .
\]
We set $\Sh_{\eX} = \Sh_{\eX , \emptyset}$. In this case the generic fibre $\pi_Q$ is a finite \'etale covering. A morphism from $\pi_1 : \Yh_1 \to \eX$ to $\pi_2 : \Yh_2 \to \eX$ in $\Sh_{\eX,D}$ is given by a morphism $\varphi : \Yh_1 \to \Yh_2$ such that $\pi_1 = \pi_2 \verk \varphi$. Note that $\varphi$ is finitely presented and proper and that $\varphi_Q$ is finite, and \'etale over $X \ohne D$.

If such a morphism exists, we say that $\pi_1$ dominates $\pi_2$. If in addition $\varphi_Q$ induces an isomorphism of the local rings in two generic points we say that $\pi_1$ {\it strictly} dominates $\pi_2$. In the case where $\Yh_{1Q}$ and $\Yh_{2Q}$ are both smooth projective curves this means that $\varphi_Q$ is an isomorphism. 

It is clear that finite products and finite fibre products exist in $\Sh_{\eX,D}$. Moreover, for every morphism $f : \eX \to \eX'$ of models over $R$ and every divisor $D'$ on $X'$, the fibre product induces a functor $f^{-1} : \Sh_{\eX',D'} \to \Sh_{\eX, f^* D'}$. 

We frequently use the fact that any non-constant morphism of a reduced and irreducible scheme $\eZ$ to a discrete valuation ring is flat, cf. \cite{Liu}, Corollary 4.3.10. Besides, note that if $\eZ$ is flat and of finite presentation over $R$ with irreducible and reduced generic fibre, then $\eZ$ is also irreducible and reduced by \cite{Liu}, Proposition 4.3.8.

We define the full subcategory
\[
\Sh^{\good}_{\eX,D} \subset \Sh_{\eX,D}
\]
to consist of those objects in $\Sh_{\eX,D}$ whose structural morphism $\lambda : \Yh \to \spec R$ is flat and satisfies $\lambda_* \Oh_{\Yh} = \Oh_{\spec R}$ universally and whose generic fibre $\lambda_Q : \Yh_Q \to \spec Q$ is smooth. In particular $\Yh_Q$ is geometrically connected and hence a smooth projective curve, which implies that $\Yh$ is irreducible and reduced.

Let $\Sh^{ss}_{\eX,D}$ denote the full subcategory of $\Sh_{\eX,D}$ consisting of all $\pi : \Yh \to \eX$ such that $\lambda : \Yh \to \spec R$ is a semistable curve whose generic fibre $\Yh_Q$ is a smooth projective curve over $Q$. Recall that $\lambda : \Yh \to \spec R$ is a semistable curve iff $\lambda$ is flat and for all $s \in \spec R$ the geometric fibre $\Yh_{\os}$ is reduced with only ordinary double points as singularities, see \cite{DM} or \cite{Liu}, section 10.3. Note that since $\Yh_Q$ is irreducible and reduced, the scheme $\Yh$ is irreducible and reduced as well. If $R$ is a discrete valuation ring, then $\Yh$ is normal since $\Yh_Q$ is normal, see \cite{Liu}, Proposition 10.3.15.

\begin{theorem}
  \label{cd_t1} Assume that the base ring $R$ is a discrete valuation ring.\\
1) The category $\Sh^{ss}_{\eX,D}$ is a full subcategory of $\Sh^{\good}_{\eX,D}$.\\
2) The objects $\Yh \to \eX$ of $\Sh^{ss}_{\eX,D}$ have the property that $\Pic^0_{\Yh/R}$ exists as a semiabelian scheme which is isomorphic to the identity component of the N\'eron model of the abelian variety $\Pic^0_{\Yh_Q / Q}$.\\
3) For any discrete valuation ring $R'$ dominating $R$ set $\eX' = \eX \otimes_R R'$ and let $D'$ be the inverse image of $D$ in $X'$. The natural base extension functor $\Sh_{\eX,D} \to \Sh_{\eX',D'}$ maps $\Sh^{\good}_{\eX,D}$ into $\Sh^{\good}_{\eX',D'}$ and $\Sh^{ss}_{\eX,D}$ into $\Sh^{ss}_{\eX',D'}$. (More generally this  is true for valuation rings $R$ and $R'$.)\\
4) For any finite number of objects $\pi_i : \Yh_i \to \eX$ in $\Sh_{\eX,D}$ there exists a finite extension $Q' / Q$ such that the objects $\pi_i \otimes_R R'$ of $\Sh_{\eX',D'}$ are all dominated by a single object of $\Sh^{\good}_{\eX',D'}$ and even of $\Sh^{ss}_{\eX',D'}$. Here $R'$ is a discrete valuation ring in $Q'$ dominating $R$.\\
5) For any object $\pi : \Yh \to \eX$ of $\Sh_{\eX,D}$ there exists an extension of discrete valuation rings $R' / R$ as in 4) such that $\pi \otimes_R R'$ is {\bf strictly} dominated by an object of $\Sh^{\good}_{\eX',D'}$ and even of $\Sh^{ss}_{\eX',D'}$. 
\end{theorem}

\begin{proof}
  1) Let $\pi : \Yh \to \eX$ be an object of $\Sh^{ss}_{\eX,D}$. By assumption the geometric fibres of $\Yh$ over $\spec R$ are reduced. Together 
with the flatness of $\lambda : \Yh \to \spec R$ it follows from \cite{EGAIII} 7.8.6 that $\lambda$ is cohomologically flat in dimension zero. This means that the formation of $\lambda_* (\Oh)$ commutes with arbitrary base changes. Since $\lambda$ is proper the sheaf $\lambda_* (\Oh)$ on $\spec R$ is coherent and hence given by the finitely generated $R$-module $\Gamma (\spec R , \lambda_* (\Oh)) = \Gamma (\Yh , \Oh)$. Since $\Yh$ is integral, this module is torsion free, hence free, so that $\lambda_* (\Oh) = \Oh^r_{\spec R}$ for some $r \ge 1$. Since $\Yh_Q$ is a smooth curve, it follows that $r = 1$. Taken together we find that the equation $\lambda_* (\Oh) = \Oh_{\spec R}$ holds universally.\\
2) Since $\Yh$ has semistable reduction over $\spec R$ it follows from \cite{BLR} 9.4, Theorem 1 that $\Pic^0_{\Yh / R}$ is a smooth separated $R$-scheme which is semi-abelian. By \cite{BLR} 9.7, Corollary 2 the connected component of the N\'eron model of $\Pic^0_{\Yh_Q / Q}$ is canonically isomorphic to $\Pic^0_{\Yh / R}$.\\
3) Note here that semistability is by definition preserved under base change.\\
4) Since finite products exist in $\Sh_{\eX,D}$ assertion 4) follows from assertion 5).\\
5) {\bf I} Let us first prove the claim for the category $\Sh^{\good}_{\eX,D}$. This proof will be taken up in a $G$-equivariant context in theorem \ref{cd_t4} below. Let $Q'$ be a finite extension field of $Q$ such that $\Yh$ has a $Q'$-rational point over $X \ohne D$ and such that the irreducible components of $\Yh_{Q'}$ are geometrically irreducible. Let $R'$ be a discrete valuation ring in $Q'$ dominating $R$. Set $\Yh_{R'} = \Yh \otimes_R R'$. Choose an irreducible component of $\Yh_{Q'}$ containing a $Q'$-rational point over $X \ohne D$ and let $\Yh^*$ be its closure in $\Yh_{R'}$ with the reduced scheme structure. Then $\Yh^*$ is integral and we can pass to its normalization $\tYh$ which is finite over $\Yh^*$ by \cite{EGAIV} (7.8.6). $\tYh$ is a proper, flat $R'$-scheme. Since $\tYh \otimes_{R'} Q'$ is the normalization of $\Yh^*_{Q'}$ it has a $Q'$-rational point. By Lipman's resolution of singularities, there is an irreducible regular $R'$-scheme $\Yh^{\vee}$ together with a proper $R'$-morphism $\Yh^{\vee} \to \tYh$ which is an isomorphism on the generic fibre. $\Yh^{\vee}$ is obtained by repeatedly blowing up the singular locus followed by normalization. This process becomes stationary after finitely many steps (see \cite{Li} and also \cite{Liu} 8.3.44). Hence we obtain a regular, irreducible scheme $\Yh^{\vee}$, which is proper and flat over $R'$, together with a proper morphism $\Yh^{\vee} \to \eX'$ strictly dominating $\pi \otimes_R R'$. The $Q'$-rational point in the generic fibre of $\Yh^{\vee}$ induces a section of $\Yh^{\vee} \to \spec R'$ by properness. Now we apply a theorem of Raynaud to deduce that $\Yh^{\vee}$ is cohomologically flat in dimension $0$, see \cite{Ray}, Th\'eor\`eme (8.2.1) (ii) $\Rightarrow$ (iv) or \cite{Liu} 9.1.24 and 9.1.32. Thus $\Yh^{\vee} \to \eX'$ lies in $\Sh^{\good}_{\eX' , D'}$.

{\bf II} Alternatively, at least if the residue field of $R$ is perfect the claim for $\Sh^{\good}_{\eX,D}$ could be proved by using instead of Raynaud's theorem a theorem of Epp. Replacing $\Yh$ by $\tYh$ and $R$ by $R'$ (of {\bf I}) we may assume that $\Yh$ is normal and that $\Yh_Q$ is a smooth projective curve over $Q$. Using \cite{Ep} Theorem 2.0, it can be shown that there is a finite extension $Q'$ of $Q$ and a discrete valuation ring $R'$ in $Q'$ dominating $R$ such that the normalization $\tYh$ of $\Yh \otimes_R R'$ has geometrically reduced fibres. As in the proof of part 1) it follows that the object $\tpi : \tYh \to \Yh  \otimes_R R' \to \eX \otimes_R R' = \eX'$ strictly dominating $\pi' = \pi \otimes_R R' : \Yh \otimes_R R' \to \eX'$ is in $\Sh^{\good}_{\eX',D'}$.

{\bf III} We now prove that after base extension every object $\pi$ of $\Sh_{\eX,D}$ is strictly dominated by an object of $\Sh^{ss}_{\eX,D}$. In view of part 1) this gives a third proof for the assertion on $\Sh^{\good}_{\eX,D}$. We construct $\Yh^{\vee} \to \eX'$ as in {\bf I}. Since $\Yh^{\vee}$ is irreducible, regular and proper and flat over $R'$, a result of Lichtenbaum \cite{Lich} implies that $\Yh^{\vee}$ is projective over $R'$. According to \cite{Liu2}, Theorem 0.2, there is a finite extension $Q^{\dagger}$ of $Q'$ and a discrete valuation ring $R^{\dagger}$ in $Q^{\dagger}$ dominating $R'$ and a semistable model $\Yh^{\dagger}$ of $\Yh^{\vee} \otimes_{R'} Q^{\dagger}$ {\it together with a morphism} $\Yh^{\dagger} \xrightarrow{\varphi} \Yh^{\vee} \otimes_{R'} R^{\dagger}$ over $\spec R^{\dagger}$. The composition
\[
\Yh^{\dagger} \xrightarrow{\varphi} \Yh^{\vee} \otimes_{R'} R^{\dagger} \to \eX^{\dagger} = \eX \otimes_R R^{\dagger}
\]
defines an object of $\Sh^{ss}_{\eX^{\dagger},D^{\dagger}}$ which strictly dominates $\pi^{\dagger} = \pi \otimes_R R^{\dagger}$. 
\end{proof}

The next result is used later to prove that certain categories of vector bundles are stable under extensions and contain all line bundles of degree zero.

As before let $R$ be a discrete valuation ring with quotient field $Q$ of characteristic zero. Consider a smooth projective curve of nonzero genus $X$ over $Q$ with a $Q$-rational point $x$ and a semistable model $\eX$ of $X$ over $\spec R$. Fix some $N \ge 1$ and define an \'etale covering $\alpha : Y \to X$ by the cartesian diagram
\[
\xymatrix{
Y \ar[r]^-i \ar[d]_{\alpha} & \Alb_{X / Q} \ar[d]^N \\
X \ar[r]^-{i_x} & \Alb_{X / Q}
}
\]
Here $i_x$ is the canonical immersion into the Albanese variety corresponding to the rational point $x$. Note that $Y$ is geometrically connected and hence a smooth projective curve.

\begin{prop}
  \label{cd_t2}
In the above situation, there exist \\
$\hullet$ a finite extension $Q' / Q$ and a discrete valuation ring $R'$ in $Q'$ dominating $R$\\
$\hullet$ a semistable model $\Yh'$ of $Y' = Y \otimes_Q Q'$ over $\spec R'$ \\
$\hullet$ a morphism
\[
\pi' : \Yh' \longrightarrow \eX' = \eX \otimes_R R'
\]
such that the following assertions hold:\\
a) The generic fibre $\pi'_{Q'}$ of $\pi'$ is $\alpha' = \alpha \otimes_Q Q'$.\\
b) There is a commutative diagram
\[
\xymatrix{
\Pic^0_{\eX' / R'} \ar[rr]^{\pi^{'*}} \ar[dr]_N & & \Pic^0_{\Yh' / R'} \\
 & \Pic^0_{\eX' / R'} \ar[ur]_g & 
}
\]
for some morphism $g$ with $g (0) = 0$, where $0$ denotes the zero section over $\spec R'$.
\end{prop}

\begin{rem}
After proving the proposition, we saw that in \cite{Fa2} Faltings uses a similar construction to make Higgs bundles on $p$-adic curves ``small''.  
\end{rem}

\begin{proof}
  Let $\Yh_1$ be the normalization of $\eX$ in the function field $Q (Y)$ of $Y$. Then $\Yh_1$ is a model of $Y$ which is equipped with a morphism $\pi_1 : \Yh_1 \to \eX$. According to \cite{EGAIV} 7.8.3 (vi) the morphism $\pi_1$ is finite. We will view $\pi_1$ as an object of $\Sh_{\eX}$. For an extension $R' / R$ as in theorem \ref{cd_t1} part 5) there exists an object $\pi' : \Yh' \to \eX'$ of $\Sh^{ss}_{\eX'}$ strictly dominating $\pi_1 \otimes_R R'$. Changing the identification of $\Yh' \otimes_{R'} Q'$ with $Y' = Y \otimes_Q Q'$ if necessary, we may assume that the generic fibre of $\pi'$ is $\alpha' = \alpha \otimes_Q Q'$. 

The origin in $\Alb_{X / Q} = \widehat{\Pic^0_{X / Q}}$ and the point $x$ of $X$ define a $Q$-rational point $y$ of $Y$ with $i (y) = 0$. Let
\[
i_y : Y \longrightarrow \Alb_{Y / Q}
\]
be the corresponding immersion. By the universal property of the Albanese variety, there is a unique morphism $f : \Alb_{Y/Q} \to \Alb_{X/Q}$ which is necessarily a homomorphism such that $f \verk i_y = i$.

Applying the functor $\Pic^0_{\_ \otimes Q' / Q'}$ to the commutative diagram
\[
\xymatrix{
Y \ar[dd]_{\alpha} \ar[r]^-{i_y} & \Alb_{Y/Q} \ar[dr]^f & \\
 & & \Alb_{X / Q} \ar[dl]^N \\
X \ar[r]^-{i_x} & \Alb_{X / Q} & 
}
\]
we obtain the following commutative diagram, where $f' = f \otimes_Q Q'$:
\[
\xymatrix{
\Pic^0_{Y' / Q'} & \\
 & \Pic^0_{X' / Q'} \ar[ul]_{\hat{f}'} \\
\Pic^0_{X' / Q'} \ar[uu]^{\alpha^{'*}} \ar[ur]_{\hat{N} = N}
}
\]
Let $\Nh$ be the N\'eron model of $\Pic^0_{Y' / Q'}$ over $\spec R'$ and let $\Nh^0$ be its identity component. By theorem \ref{cd_t1} part 2) we know that $\Pic^0_{\eX' / R'}$ and $\Pic^0_{\Yh' / R'}$ exist as smooth and separated schemes and that $\Pic^0_{\Yh' / R'}$ is isomorphic to $\Nh^0$. By the universal property of the N\'eron model, the natural map
\begin{equation}
  \label{eq:cd1}
  \Mor_{R'} (\Pic^0_{\eX' / R'} , \Nh) \silo \Mor_{Q'} (\Pic^0_{X' / Q'} , \Pic^0_{Y' / Q'})
\end{equation}
is bijective. Hence $\hat{f}'$ has a unique extension to a morphism $g : \Pic^0_{\eX' / R'} \to \Nh$. By construction, the composition $g \verk N$ has generic fibre $\hat{f}' \verk N = \alpha^{'*}$. Since $\alpha'$ is the generic fibre of $\pi' : \Yh' \to \eX'$, the induced homomorphism
\[
\pi^{'*} : \Pic^0_{\eX' / R'} \longrightarrow \Pic^0_{\Yh' / R'} 
\]
has generic fibre $\alpha^{'*}$ as well. Using the N\'eron property (\ref{eq:cd1}) it follows that $g \verk N$ is equal to the composition
\[
\Pic^0_{\eX' / R'} \xrightarrow{\pi^{'*}} \Pic^0_{\Yh' / R'} = \Nh^0 \hookrightarrow \Nh \; .
\]
In particular we get that $g (0) = g (N (0)) = 0$ where $0$ denotes the zero sections of $\Pic^0_{\eX' / R'}$ respectively $\Pic^0_{\Yh' / R'}$. Since the special fibre of $\Pic^0_{\eX' / R'}$ is connected, it follows that $g$ is a morphism
\[
g : \Pic^0_{\eX' / R'} \longrightarrow \Nh^0 = \Pic^0_{\Yh' / R'} 
\]
with $g \verk N = \pi^{'*}$ as desired.
\end{proof}

We fix an algebraic closure $\oQ_p$ of $\Q_p$ and consider finite extensions $\Q_p \subset K \subset \oQ_p$. The rings of integers will be denoted by $\eo_K$ and $\eo_{\oK} = \oZ_p$. 

The following corollary of theorem \ref{cd_t1} will be used constantly.

\begin{cor}
  \label{cd_t3}
Let $X$ be a smooth projective curve over $\oQ_p$ and $D$ a divisor on $X$. Let $\eX$ be a model of $X$ over $\spec \oZ_p$.\\
1) Given any finite number of objects $\pi_i : \Yh_i \to \eX$ in $\Sh_{\eX,D}$ (resp. given one object $\pi_1 : \Yh_1 \to \eX$ in $\Sh_{\eX,D}$) there is a finite extension $K$ of $\Q_p$ and a curve $X_K / K$ with model $\eX_{\eo_K} / \eo_K$ and a divisor $D_K$ of $X_K$ such that the following hold: We have $X = X_K \otimes_K \oK$ and $D = D_K \otimes \oK$ and $\eX = \eX_{\eo_K} \otimes_{\eo_K} \oZ_p$ and there is an object $\pi_{\eo_K} : \Yh_{\eo_K} \to \eX_{\eo_K}$ of $\Sh^{\good}_{\eX_{\eo_K} , D_K}$ and even of $\Sh^{ss}_{\eX_{\eo_K} , D_K}$ such that $\pi = \pi_{\eo_K} \otimes_{\eo_K} \oZ_p$ dominates all $\pi_i$ in $\Sh_{\eX,D}$ (resp. dominates $\pi_1$ {\rm strictly}).\\
2)The category $\Sh^{ss}_{\eX,D}$ is a full subcategory of $\Sh^{\good}_{\eX,D}$.\\
3) Any finite number of objects $\pi_i : \Yh_i \to \eX$ in $\Sh_{\eX,D}$ are dominated by a common object $\pi : \Yh \to \eX$ of $\Sh^{\good}_{\eX,D}$ and even of $\Sh^{ss}_{\eX,D}$. Every single object $\pi_1 : \Yh_1 \to \eX$ in $\Sh_{\eX,D}$ is {\rm strictly} dominated by an object of $\Sh^{\good}_{\eX,D}$ and even of $\Sh^{ss}_{\eX,D}$.
\end{cor}

\begin{proof}
  Part 1) follows from theorem \ref{cd_t1}, 4), 5) using noetherian descent as in \cite{EGAIV} \S\,8, in particular (8.8.3) and (8.10.5), together with \cite{EGAIV} (17.7.8) to descend to the category $\Sh_{\eX_1,D_1}$ for some $\eX_1 / \eo_{K_1}$ with divisor $D_1$ where $K_1 \supset \Q_p$ is a finite extension. \\
2) Similarly as above, every object $\pi : \Yh \to \eX$ of $\Sh^{ss}_{\eX,D}$ descends to an object $\pi_{\eo_K} : \Yh_{\eo_K} \to \eX_{\eo_K}$ of $\Sh_{\eX_{\eo_K} , D_K}$ where $K \supset \Q_p$ is finite such that $\Yh_{\eo_K} / \eo_K$ is flat. Since the geometric fibres of $\Yh_{\eo_K} / \eo_K$ and $\Yh / \oZ_p$ can be identified, it follows that $\pi_{\eo_K} : \Yh_{\eo_K} \to \eX_{\eo_K}$ is in $\Sh^{ss}_{\eX_{\eo_K}, D_K}$ and hence in $\Sh^{\good}_{\eX_{\eo_K}, D_K}$ by theorem \ref{cd_t1}, 1). Therefore $\pi = \pi_{\eo_K} \otimes_{\eo_K} \oZ_p$ lies in $\Sh^{\good}_{\eX,D}$ by theorem \ref{cd_t1}, 3). \\
Part 3) follows by combining 1) and theorem \ref{cd_t1}, 3).
\end{proof}

Later we will construct a canonical parallel transport for certain vector bundles. The proof that it is well defined requires the following theorem. Let $\eT_{\eX,D}$ be the following category. Objects are finitely presented proper $G$-equivariant morphisms $\pi : \Yh \to \eX$ over $\spec \oZ_p$ where $G$ is a finite (abstract) group which acts $\oZ_p$-linearly from the left on $\Yh$ and trivially on $\eX$. Moreover the generic fibre $\pi_{\oQ_p}$ is finite and its restriction $\Yh_{\oQ_p} \ohne \pi^* D \to X \ohne D$ is an \'etale $G$-torsor.

A morphism from the $G_1$-equivariant morphism $\pi_1 : \Yh_1 \to \eX$ to the $G_2$-equivariant morphism $\pi_2 : \Yh_2 \to \eX$ in $\eT_{\eX,D}$ is given by a morphism $\varphi : \Yh_1 \to \Yh_2$ with $\pi_1 = \pi_2 \verk \varphi$ together with a homomorphism $\gamma : G_1 \to G_2$ of groups such that $\varphi$ is $G_1$-equivariant if $G_1$ acts on $\Yh_2$ via $\gamma$.

This definition generalizes the category $\eT_{\eX} = \eT_{\eX, \phi}$ used in \cite{De-We} \S\,5. There is an obvious forgetful functor $\eT_{\eX,D} \to \Sh_{\eX,D}$. The full subcategory $\eT^{\good}_{\eX,D}$ of $\eT_{\eX,D}$ consists of those objects which are mapped to objects of $\Sh^{\good}_{\eX,D}$. 

\begin{theorem}
  \label{cd_t4}
For any object $\pi : \Yh \to \eX$ in $\Sh_{\eX,D}$ there is a finite group $G$ and a $G$-equivariant morphism $\pi' : \Yh' \to \eX$ defining an object of $\eT^{\good}_{\eX,D}$ which admits a morphism $\varphi : \Yh' \to \Yh$ with $\pi \verk \varphi = \pi'$. In other words, every object of $\Sh_{\eX,D}$ is dominated by the image of an object in $\eT^{\good}_{\eX,D}$.
\end{theorem}

\begin{proof}
Let us first show that every object $\pi : \Yh \to \eX$ of $\eT_{\eX,D}$ is dominated by an object of $\eT^{\good}_{\eX,D}$. By noetherian descent we can assume that there is a finite extension $K$ of $\Q_p$ in $\oQ_p$ with ring of integers $R$ such that $\pi$ descends to the object $\pi_R : \Yh_R \to \eX_R$ in $\eT_{\eX_R , D_K}$. Denote by $G$ the group acting on $\Yh_R$ over $\eX_R$ such that $Y_K  \ohne \pi^*_K D_K \to X_K \ohne D_K$ is an \'etale $G$-torsor. Now we follow the construction in the proof of theorem \ref{cd_t1}, 5) {\bf I} and consider a geometrically irreducible component of $Y_{K'}$ containing a $K'$-rational point over $X_{K'} \ohne D_{K'}$, where $K'$ is a finite extension of $K$ in $\oQ_p$. Denote by $H \subseteq G$ the stabilizer of this component. Then $H$ acts in a natural way on $\Yh^*$, and also on $\tYh$ and $\Yh^{\vee}$. Therefore $\Yh^{\vee} \to \eX \otimes_R R'$, where $R'$ is the ring of integers in $K'$, is an object of $\eT^{\good}_{\eX_{R'} , D_{K'}}$ dominating $\pi_{R'}$. By base-change to $\oZ_p$, our claim follows. Hence it suffices to show that there is an object $\pi'$ of $\eT_{\eX,D}$ which dominates $\pi$. By corollary \ref{cd_t3}, 1), we may assume that we have $\pi = \pi_{\eo_K} \otimes \oZ_p$ with $\pi_{\eo_K} : \Yh_{\eo_K} \to \eX_{\eo_K}$ in $\Sh^{\good}_{\eX_{\eo_K} , D_K}$. Let $Y'_K$ be the smooth projective curve whose function field is the Galois closure of $K (Y_K)$ over $K (X_K)$. The Galois group $G$ acts on $Y'_K$ over $X_K$. The morphism $Y'_K \to X_K$ is finite and over $X_K \ohne D_K$ it defines a Galois covering with group $G$. Consider the normalization $\tYh_{\eo_K}$ of $\Yh_{\eo_K}$ in $K (Y'_K)$. By \cite{EGAIV}, (7.8.3) (vi) the morphism $\tYh_{\eo_K} \to \Yh_{\eo_K}$ is finite. Hence $\tYh_{\eo_K} \to \Yh_{\eo_K} \to \eX_{\eo_K}$ defines an object of $\Sh_{\eX_{\eo_K} , D_K}$ with generic fibre $\tY_K = Y'_K \to X_K$. 

By the proof of \cite{Liu2}, Lemma 2.4, there exists a model $\Yh'_{\eo_K}$ of $Y'_K$ over $\eo_K$ endowed with an action of $G$ extending the action on $Y'_K$ together with a morphism $\varphi_{\eo_K} : \Yh'_{\eo_K} \to \Yh_{\eo_K}$ which is an isomorphism on the generic fibre.

Let $\pi'_{\eo_K} : \Yh'_{\eo_K} \to \eX_{\eo_K}$ be the composition $\pi'_{\eo_K} = \pi_{\eo_K} \verk \varphi_{\eo_K} : \Yh'_{\eo_K} \to \Yh_{\eo_K} \to \eX_{\eo_K}$. Since $\Yh'_{\eo_K}$ is reduced, $G$-equivariance of the generic fibre of $\pi'_{\eo_K}$ implies $G$-equivariance of $\pi'_{\eo_K}$ cf. \cite{EGAII}, 7.2.21.

 Now put
\[
\Yh' := \Yh'_{\eo_K} \otimes_{\eo_K} \oZ_p \quad \mbox{and} \quad \pi' = \pi'_{\eo_K} \otimes_{\eo_K} \oZ_p : \Yh' \longrightarrow \eX \; .
\]
Then $\pi' : \Yh' \to \eX$ is an object of $\eT_{\eX,D}$ such that $\varphi = \varphi_{\eo_K} \otimes_{\eo_K} \oZ_p : \Yh' \to \Yh$ satisfies $\pi \verk \varphi = \pi'$. 
\end{proof}

\section{Two categories of vector bundles on $p$-adic curves}

Let $\uVec_S$ be the category of vector bundles on a scheme $S$. For a bundle $E$ we often write $E$ for its locally free sheaf of sections $\Oh (E)$. Let $\eo$ be the ring of integers in $\C_p = \Hat{\oQ_p}$ and set $\eo_n = \eo / p^n \eo = \oZ_p / p^n \oZ_p$. For every $\eo$-scheme $\Yh$ we set $\Yh_n = \Yh \otimes_{\eo} \eo_n$. Let $X$ be as before a smooth projective curve over $\oQ_p$ and set $X_{\C_p} = X \otimes_{\oQ_p} \C_p$.

First of all, we show that vector bundles on $X_{\C_p}$ can be extended to vector bundles on suitable models. The elegant argument in the proof was commicated to us by M.~Raynaud.

\begin{theorem}
  \label{aw_t5}
For every vector bundle $E$ on $X_{\C_p}$ and every model $\eX$ of $X$ there exists a model $\eX'$ of $X$ dominating $\eX$ such that $E$ extends to a vector bundle on $\eX'_{\eo}$. If $\eX$ is smooth, then $E$ can be extended to a vector bundle on $\eX_{\eo}$ itself.
\end{theorem}

\begin{proof}
  We can extend $E$ to a quasi-coherent sheaf $\Fh$ of finite presentation on $\eX_{\eo}$, see \cite{Ha}, Appendix, Corollary 2 to Proposition 2. Let $\Jh \subset \Oh_{\eX_{\eo}}$ be the $r$-th Fitting ideal of $\Fh$, where $r$ is the rank of $E$. Since $\Fh$ is of finite presentation, $\Jh$ is quasi-coherent of finite type. Besides, $\Jh \cdot \Oh_{X_{\C_p}}$ is equal to the Fitting ideal of $E$, hence to $\Oh_{X_{\C_p}}$. Therefore there exists some $n \ge 1$ such that $p^n \Oh_{\eX_{\eo}} \subset \Jh$. By appoximating the local generators of $\Jh$ with elements in $\Oh_{\eX}$ modulo $p^n$, we see that $\Jh$ descends to an ideal $\Jh_0 \subset \Oh_{\eX}$. Let $\varphi : \eX' \to \eX$ be the blowing-up of $\Jh_0$. Since $\Jh_0$ is of finite type, $\varphi$ is of finite presentation, so that $\varphi$ is a map in $\Sh_{\eX}$ inducing an isomorphism on the generic fibre. The base change map $\varphi_{\eo} : \eX'_{\eo} \to \eX_{\eo}$ is the blowing-up of $\Jh$. Hence $\varphi^{-1}_{\eo} (\Jh) \Oh_{\eX'_{\eo}}$ is invertible. Since $\varphi^{-1}_{\eo} (\Jh) \Oh_{\eX'_{\eo}}$ is the $r$-th Fitting ideal $F_r (\varphi^*_{\eo} \Fh)$ of $\varphi^*_{\eo} \Fh$, we can apply \cite{Ray-Gru}, (5.4.3) to deduce that $\varphi^*_{\eo} \Fh / \Ann_{\varphi^*_{\eo} \Fh} (F_r (\varphi^*_{\eo} \Fh))$ is locally free of rank $r$ on $\eX'_{\eo}$. Hence it gives rise to a vector bundle $\Eh$ on $\eX'_{\eo}$ with generic fibre $E$.

If $\eX$ is smooth over $\oZ_p$, then $\Pic^0_{\eX / \oZ_p} (\eo) = \Pic^0_{X / \oQ_p} (\C_p)$, so that every line bundle of degree $0$ extends to a line bundle on $\eX_{\eo}$. Besides, $\eX$ carries a line bundle $\Nh$ whose generic fibre has rank one. Hence every line bundle on $X_{\C_p}$ can be extended to $\eX_{\eo}$. The general case follows by induction on the rank of $E$. Namely, there is an exact sequence of vector bundles $0 \to E_1 \to E \to E_2 \to 0$ on $X_{\C_p}$ where $\rk E_i < \rk E$ for $i = 1,2$. By hypothesis, $E_1$ and $E_2$ can be extended to $\eE_1$ and $\eE_2$ on $\eX_0$. By flat base change we have an isomorphism
\[
\Ext^1_{\eX_{\eo}} (\eE_2 , \eE_1) \otimes_{\eo} \C_p \silo \Ext^1_{X_{\C_p}} (E_2 , E_1) \; .
\]
This implies that $E$ is isomorphic to the generic fibre of a vector bundle $\eE$ on $\eX_{\eo}$. Note here that extensions of locally free sheaves are locally free because the cohomology of affine schemes vanishes.
\end{proof}

\begin{defn}
  \label{cd_t5}
{\bf a} For a model $\eX$ of $X$ over $\oZ_p$ and a divisor $D$ in $X$ the category $\eB_{\eX_{\eo},D}$ is defined to be the full subcategory of $\uVec_{\eX_{\eo}}$ consisting of vector bundles $\eE$ on $\eX_{\eo} = \eX \otimes_{\oZ_p} \eo$ with the following property: For every $n \ge 1$ there is an object $\pi : \Yh \to \eX$ of $\Sh_{\eX,D}$ such that $\pi^*_n \eE_n$ is a trivial bundle on $\Yh_n$. Here $\pi_n, \Yh_n$ and $\eE_n$ are the reductions $\mod p^n$ of $\pi, \Yh$ and $\eE$.\\
{\bf b} The full subcategory $\eB_{X_{\C_p},D}$ of $\uVec_{X_{\C_p}}$ consists of all vector bundles on $X_{\C_p}$ which are isomorphic to a bundle of the form $j^* \eE$ with $\eE$ in $\eB_{\eX_{\eo},D}$ for some model $\eX$ of $X$. Here $j$ is the open immersion of $X_{\C_p}$ into $\eX_{\eo}$.\\
{\bf c} The full subcategory $\eB^{\sharp}_{X_{\C_p},D}$ of $\uVec_{X_{\C_p}}$ consists of all vector bundles $E$ on $X_{\C_p}$ such that $\alpha^*_{\C_p} E$ is in $\eB_{Y_{\C_p},\alpha^*D}$ for some finite covering $\alpha : Y \to X$ of $X$ by a smooth projective curve $Y$ over $\oQ_p$ such that $\alpha$ is \'etale over $X \ohne D$.
\end{defn}

\begin{rems}
{\bf a} For $D = \emptyset$ we simply write $\eB_{\eX_{\eo}}$ for $\eB_{\eX_{\eo},D}$ etc.\\
{\bf b} In \cite{De-We} \S\,6 a category $\eB_{\eX_{\eo}}$ was defined as above, but using coverings in $\eT_{\eX}$ instead of $\Sh_{\eX}$. It follows from theorem \ref{cd_t4} that both definitions give the same category. Consequently, also the category $\eB_{X_{\C_p}}$ is the same as the one defined in \cite{De-We} Definition 19. 
\end{rems}

\begin{lemma}
  \label{cd_t7a}
The category $\eB_{X_{\C_p},D}$ consists of all vector bundles isomorphic to $j^* \eE$ with $\eE$ in $\eB_{\eX_{\eo},D}$ and $\eX$ a {\rm semistable} model of $X$ over $\oZ_p$.
\end{lemma}

\begin{proof}
  Given any model $\eX$ of $X$, there is a semistable model $\Yh$ of $X$ strictly dominating $\eX$. This follows from corollary \ref{cd_t3}, 3) applied to $\pi_1 = \id_{\eX}$. Since the pullback of bundles on $\eX$ to $\Yh$ maps $\eB_{\eX_{\eo},D}$ to $\eB_{\Yh_{\eo},D}$ by proposition \ref{cd_t6} below, the assertion follows.
\end{proof}

\begin{lemma}
  \label{cd_t7neu}
Let $f : X \to X'$ be a morphism of smooth, projective curves over $\oQ_p$. For every model $\eX'$ of $X'$ there exists a model $\eX$ of $X$ and a $\oZ_p$-linear morphism $\tilde{f} : \eX \to \eX'$ such that the diagram
\[
\xymatrix{
\eX \ar[d]_{\tilde{f}} & X \ar[l] \ar[d]^f \\
\eX' & X' \ar[l]
}
\]
is commutative.
\end{lemma}

\begin{proof}
  Since $f$ is proper, it is either surjective or maps $X$ to a closed point of $X'$. In the second case, because of properness any model $\eX$ of $X$ will do. Hence we can assume that $f$ is surjective, hence finite. There is a finite extension $K$ of $\Q_p$ in $\oQ_p$ such that $f$ descends to a morphism $f_K : X_K \to X'_K$ of smooth, proper curves over $K$ and such that $\eX'$ descends to a model $\eX'_{\eo_K}$ of $X'_K$.

Define $\eX_{\eo_K}$ as the normalization of the reduced and irreducible scheme $\eX'_{\eo_K}$ in the function field $K (X_K)$ of $X_K$, and let $\tilde{f}_{\eo_K} : \eX_{\eo_K} \to \eX'_{\eo_K}$ be the corresponding finite morphism. Since $f_K : X_K \to X'_K$ is the normalization of $X'_K$ in $K (X_K)$, the generic fibre of $\eX_{\eo_K}$ can be identified with $X_K$ so that the desired diagram commues. Base-change with $\oZ_p$ completes the proof.
\end{proof}

\begin{prop}
  \label{cd_t6}
The categories $\eB_{\eX_{\eo},D}$ resp. $\eB_{X_{\C_p},D}$ and $\eB^{\sharp}_{X_{\C_p},D}$ are full additive subcategories of $\uVec_{\eX_{\eo}}$ resp. $\uVec_{X_{\C_p}}$ which are closed under tensor products, duals, internal homs and exterior powers. For every morphism $f : \eX \to \eX'$ over $\oZ_p$ resp. $f : X \to X'$ over $\oQ_p$ and every divisor $D'$ on $X'$, the pullback functor $f^*$ of vector bundles restricts to an additive exact functor $f^* : \eB_{\eX'_{\eo},D'} \to \eB_{\eX_{\eo},f^* D'}$ resp. $f^* : \eB_{X'_{\C_p},D'} \to \eB_{X_{\C_p},f^*D'}$ and $f^* : \eB^{\sharp}_{X'_{\C_p},D'} \to \eB^{\sharp}_{X_{\C_p} , f^* D'}$. These functors commute with tensor products, duals, internal homs and exterior powers.
\end{prop}

The {\bf proof} is straightforward for $\eB_{\eX_{\eo},D}$ and $\eB_{X_{\C_p},D}$ given corollary \ref{cd_t3}, 3), lemma \ref{cd_t7neu} and the functoriality of the categories $\Sh$. For $\eB^{\sharp}$, note first that given finite morphisms $Y_i \xrightarrow{\alpha_i} X$ for $1 \le i \le n$ \'etale over $X \ohne D$ by smooth projective curves $Y_i$, there is a finite morphism $Y \xrightarrow{\beta} X$ \'etale over $X \ohne D$ by another such curve $Y$ such that $\beta$ factors over each $\alpha_i$: Take the normalization of any irreducible component of $Y_1 \times_X \ldots \times_X Y_n$. Thus the assertions about $\oplus, \otimes$ etc. for $\eB^{\sharp}$ follow from those for $\eB$. Next, given $E$ in $\eB^{\sharp}_{X'_{\C_p},D'}$ and $f : X \to X'$, choose a finite morphism $\alpha' : Y' \to X'$, \'etale over $D'$ such that $\alpha^{'*} E$ lies in $\eB_{Y'_{\C_p}, \alpha'^* D'}$. Let $Y$ be the normalization of an irreducible component of $f^{-1} (Y')$ and consider the commutative diagram:
\[
\xymatrix{
Y \ar[r] \ar[d]_{\alpha} & f^{-1} (Y') = X \times_{X'} Y' \ar[r] & Y' \ar[d]^{\alpha'} \\
X \ar[rr]^f && X'
}
\]
Let $g : Y \to Y'$ be the upper horizontal map. By functoriality of $\eB$ we know that $g^* \alpha^{'*} E$ lies in $\eB_{Y_{\C_p},g^* \alpha'^* D'}$. Hence $f^* E$ is in $\eB^{\sharp}_{X_{\C_p}, f^* D'}$. \endproof

\begin{prop}
  \label{cd_t7}
{\bf a} Let $\alpha : Y \to X$ be a finite morphism, \'etale over $X \ohne D$ of smooth and proper curves over $\oQ_p$. Then a vector bundle $E$ on $X_{\C_p}$ lies in $\eB^{\sharp}_{X_{\C_p},D}$ if and only if $\alpha^* E$ lies in $\eB^{\sharp}_{Y_{\C_p}, \alpha^* D}$. \\
{\bf b} Assume in addition that $\alpha : Y \to X$ is \'etale. For a vector bundle $F$ on $Y_{\C_p}$ let $\alpha_* F$ be the vector bundle on $X_{\C_p}$ corresponding to the locally free sheaf $\alpha_* \Oh (F)$. If $F$ is in $\eB^{\sharp}_{Y_{\C_p}}$ then $\alpha_* F$ is in $\eB^{\sharp}_{X_{\C_p}}$.
\end{prop}

\begin{proof}
  {\bf a} This follows from the functoriality of $\eB^{\sharp}$ in proposition \ref{cd_t6}. \\
{\bf b} Consider $F$ in $\eB^{\sharp}_{Y_{\C_p}}$ and choose a Galois covering $\gamma : Y' \to X$ which factors over $Y$ i.e. $\gamma$ is a composition $\gamma : Y' \xrightarrow{\beta} Y \xrightarrow{\alpha} X$. Let $G$ be the Galois group of $Y'$ over $X$ and let $H$ be the one of $Y'$ over $Y$. For every $\sigma$ in $G$ the adjunction map $F \to \beta_* \beta^* F$ induces a map
\[
\alpha_* F \longrightarrow \alpha_* \beta_* \beta^* F = \gamma_* \beta^* F = \gamma_* \sigma_* \beta^* F \; .
\]
Note here that $\gamma \verk \sigma = \gamma$. This gives a map
\[
\gamma^* \alpha_* F \longrightarrow \sigma_* \beta^* F \; .
\]
For $\tau$ in $H$ we have $\tau_* \beta^* = (\tau^{-1})^* \beta^* = (\beta \verk \tau^{-1})^* = \beta^*$. Hence we obtain a well defined map
\begin{equation}
  \label{eq:cd2}
  \gamma^* \alpha_* F \longrightarrow \bigoplus_{\sigma \in G \mod H} \sigma_* \beta^* F \; .
\end{equation}
Argueing locally, one sees that (\ref{eq:cd2}) is an isomorphism. Now, $\sigma_* \beta^* F = (\sigma^{-1})^* \beta^* F$ belongs to $\eB^{\sharp}_{Y'_{\C_p}}$ by functoriality of $\eB^{\sharp}$. Hence $\gamma^* \alpha_* F$ belongs to this category as well. It follows that $\alpha_* F$ lies in $\eB^{\sharp}_{X_{\C_p}}$ as was to be shown.
\end{proof}

We now prove that our categories are stable under extensions of vector bundles. 

\begin{theorem}
  \label{cd_t8}
The categories $\eB_{\eX_{\eo},D} , \eB_{X_{\C_p},D}$ and $\eB^{\sharp}_{X_{\C_p},D}$ are stable under extensions, e.g. if
\[
0 \longrightarrow E' \longrightarrow E \longrightarrow E'' \longrightarrow 0
\]
is an exact sequence of vector bundles on $X_{\C_p}$ such that $E'$ and $E''$ are objects of the category $\eB_{X_{\C_p},D}$, then $E$ is also contained in $\eB_{X_{\C_p},D}$.
\end{theorem}

\begin{proof}
We give the proof for $\eB_{X_{\C_p},D}$. The case of $\eB_{\eX_{\eo},D}$ is similar. The assertion for $\eB^{\sharp}$ follows formally from the one for $\eB$. Thus, let $E'$ und $E''$ be in $\eB_{X_{\C_p},D}$. 
  By definition, there exist models $\eX'$ and $\eX''$ of $X$ over $\oZ_p$ and vector bundles $\eE'_1$ in $\eB_{\eX'_{\eo},D}$ and $\eE''_1$ in $\eB_{\eX''_{\eo},D}$ such that
\[
E' \simeq j^*_{\eX'_{\eo}} \eE'_1 \quad \mbox{and} \quad E'' \simeq j_{\eX''_{\eo}} \eE''_1 \; ,
\]
where $j_{\eX'_{\eo}}$ and $j_{\eX''_{\eo}}$ are the open immersions of the generic fibre $X_{\C_p}$ into $\eX'_{\eo} = \eX' \otimes_{\oZ_p} \eo$ respectively $\eX''_{\eo} = \eX'' \otimes_{\oZ_p} \eo$. 

Applying Proposition \ref{cd_t17} below, there exists a model $\eX$ of $X$ over $\oZ_p$ together with morphisms over $\oZ_p$
\[
\eX' \xleftarrow{p_1} \eX \xrightarrow{p_2} \eX''
\]
restricting to the identity on the generic fibres. By functoriality $\eE' = p^*_1 \eE'_1$ and $\eE'' = p^*_2 \eE''_1$ lie in $\eB_{\eX_{\eo},D}$. 

Reducing to cohomology and using flat base change one sees that $j^*_{\eX_{\eo}}$ induces an isomorphism $\Ext^1_{\eX_{\eo}} (\eE'' , \eE') \otimes_{\eo} \C_p \xrightarrow{\sim} \Ext^1_{X_{\C_p}} (E'' , E')$. Hence there is some $k \ge 0$ such that the extension class we get by multiplying $p^k$ with the class in $\Ext^1_{X_{\C_p}} (E'' , E')$ induced by $F$ comes from $\Ext^1_{\eX_{\eo}} (\eE'' , \eE')$. 

Hence pullback by $p^k$-multiplication on $E''$ induces an extension
\[
\xymatrix{
0 \ar[r] & E' \ar[r] \ar@{=}[d] & E_1 \ar[r] \ar[d]_{\wr} & E'' \ar[r] \ar[d]^{p^k}_{\wr} & 0 \\
0 \ar[r] & E' \ar[r] & E \ar[r] & E'' \ar[r] & 0
}
\]
on $X_{\C_p}$ for which there is an exact sequence
\[
0 \longrightarrow \eE' \longrightarrow \eE \longrightarrow \eE'' \longrightarrow 0
\]
of vector bundles on $\eX_{\eo}$ such that $j^*_{\eX_{\eo}} \eE \simeq E_1 \simeq E$. Note here that any extension of a locally free sheaf by another one is again locally free. The reason is that locally every such extension splits because the coherent cohomology of affine schemes vanishes.

Let us fix some $n \ge 1$. Since $\eE'$ and $\eE''$ lie in $\eB_{\eX_{\eo},D}$, we find objects $\pi' : \Yh' \to \eX$ and $\pi'' : \Yh'' \to \eX$ of $\Sh_{\eX,D}$ such that $\pi^{'*}_n \eE'_n$ is trivial on $\Yh'_n = \Yh' \otimes_{\oZ_p} \eo_n$ and $\pi^{''*}_n \eE''_n$ is trivial on $\Yh''_n = \Yh'' \otimes_{\oZ_p} \eo_n$.

By Corollary \ref{cd_t3}, 1), there is a finite extension $K$ of $\Q_p$ with the following properties:\\
$\hullet$ $X,D$ and $\eX$ descend to a curve $X_K / K$ a divisor $D_K$ on $X_K$ and a model $\eX_R / R$ respectively, where $R = \eo_K$. \\
$\hullet$ there is an object $\pi_{R} : \Yh_{R} \to \eX_R$ of $\Sh^{ss}_{\eX_R , D_K}$ such that 
\[
\pi = \pi_{R} \otimes_{R} \oZ_p : \Yh = \Yh_{R} \otimes_{R} \oZ_p \to \eX
\]
dominates both $\pi'$ and $\pi''$. \\
$\hullet$ the generic fibre $Y_K$ of $\Yh_{R}$ has a rational point.

Now $\pi^*_n \eE'_n$ and $\pi^*_n \eE''_n$ are trivial bundles on $\Yh_n$. If $r'$ respectively $r''$ denote their ranks the extension:
\begin{equation}
  \label{eq:cd3}
  0 \longrightarrow \pi^*_n \eE'_n \longrightarrow \pi^*_n \eE_n \longrightarrow \pi^*_n \eE''_n \longrightarrow 0
\end{equation}
gives rise to a class in $\Ext^1_{\Yh_n} (\Oh^{r''} , \Oh^{r'}) \simeq H^1 (\Yh_n , \Oh)^{r'\,r''}$. 

{\bf Claim} There exists an object $\sigma : \eZ \to \eX$ in $\Sh_{\eX,D}$ and a morphism $\rho : \eZ \longrightarrow  \Yh$ in $\Sh_{\eX,D}$, such that the induced map $\rho^*_n : H^1 (\Yh_n , \Oh) \to H^1 (\eZ_n , \Oh)$ is trivial. 

Assume that the claim holds. Then $\rho^*_n$ applied to the extension (\ref{eq:cd3}) is trivial, which implies that $\sigma^*_n \eE_n = \rho^*_n \pi^*_n \eE_n$ is a trivial vector bundle on $\eZ_n$. Since this argument can be done for every $n \ge 1$ it follows that $\eE$ lies in $\eB_{\eX_{\eo},D}$, which implies that $E$ is contained in the category $\eB_{X_{\C_p},D}$. The theorem follows.

Hence it remains to prove the claim. If the genus of $Y_K$ is zero, then $Y_K \cong \Pa^1_K$ since $Y_K$ was assumed to have a rational point. Hence $\chi (Y_K , \Oh) = 1$ and therefore $\chi (Y_{\kappa} , \Oh) = 1$ where $Y_{\kappa}$ is the special fibre of $\Yh_R$. Since $\lambda_* \Oh_{\Yh_R} = \Oh_R$ holds universally we have $H^0 (\Yh_{\kappa} , \Oh) = \kappa$ and therefore $H^1 (\Yh_{\kappa} , \Oh) = 0$. Now \cite{Mu} Corollary 3 on p. 53 implies that $H^1 (\Yh_n , \Oh) = 0$. In proving the claim we can therefore assume from now on that the genus of $Y_K$ is nonzero. Let us first show that it suffices to find a morphism $\rho : \eZ \to \Yh$ in $\Sh_{\eX,D}$ such that
\[
\rho^* : H^1 (\Yh , \Oh) \longrightarrow H^1 (\eZ , \Oh)
\]
satisfies $\rho^* (H^1 (\Yh , \Oh)) \subseteq p^n H^1 (\eZ , \Oh)$. 

Namely, consider the commutative diagram
\[
\xymatrix{
H^1 (\Yh , \Oh) \otimes_{\oZ_p} \eo_n \ar[rr]^{\rho^* \otimes \eo_n} \ar[d] & & H^1 (\eZ , \Oh) \otimes_{\oZ_p} \eo_n \ar[d] \\
H^1 (\Yh_n , \Oh) \ar[rr]_-{\rho^*_n} & & H^1 (\eZ_n , \Oh)
}
\]
By assumption, the upper horizontal map is zero. Hence $\rho^*_n = 0$, if the left vertical map $H^1 (\Yh, \Oh) \otimes_{\oZ_p} \eo_n \longrightarrow H^1 (\Yh_n , \Oh)$ is surjective. Since $\oZ_p$ is flat over $R$ and therefore $\eo_n = \oZ_p / p^n \oZ_p$ is flat over $R / p^n R$, it suffices by flat base change to prove surjectivity of
\[
H^1 (\Yh_{R} , \Oh) \otimes_{R} R / p^n R \to H^1 (\Yh_{R} \otimes_{R} R / p^n R , \Oh) \; .
\]

Let $k$ be the residue field of $R$. By Nakayama's lemma, it suffices to prove surjectivity after tensoring with $k$. Consider the commutative triangle
\[
\xymatrix{
H^1 (\Yh_{R} , \Oh) \otimes_{R} k \ar[rr] \ar[dr] & & H^1 (\Yh_{R} \otimes_{R} R / p^n R , \Oh) \otimes_{R / p^n R} k \ar[dl] \\
 & H^1 (\Yh_k , \Oh) &
}
\]
Both vertical maps are isomorphisms by \cite{Mu}, Corollary 3 on p.~53 since $\Yh_k$ is one-dimensional and hence has vanishing second cohomology. Hence the horizontal map is a fortiori surjective.

By proposition \ref{cd_t2} applied to the smooth projective curve $Y_K$ over $K$ and its semistable model $\Yh_{R}$ over $R$ with $N = p^n$ there exist the following:\\
$\bullet$ a finite extension $K'$ of $K$ in $\oQ_p$ with ring of integers $R' = \eo_{K'}$\\
$\bullet$ an object
\[
\rho_{R'} : \eZ_{R'} \longrightarrow \Yh_{R'} = \Yh_{R} \otimes_{R} R'
\]
of $\Sh^{ss}_{\Yh_{R'}}$ such that there is a commutative diagram
\[
\xymatrix{
\Pic^0_{\Yh_{R'} / R'} \ar[rr]^{\rho^*_{R'}} \ar[dr]_{p^n} && \Pic^0_{\eZ_{R'} / R'} \\
 & \Pic^0_{\Yh_{R'} / R'} \ar[ur]_g
}
\]
for some morphism $g$ with $g (0) = 0$. 

Note that the Lie algebra of a group functor coincides with the Lie algebra of its identity component, if the latter exists (see e.g. \cite{SGA3I}, expose $\mathrm{VI_B}$, remarque 3.2). Hence we can apply \cite{BLR}, 8.4, Theorem 1, to the proper, flat $R'$-schemes $\Yh_{R'}$ and $\eZ_{R'}$ which as in the proof of theorem \ref{cd_t1}, 1) are both cohomologically flat in dimension $0$ over $\spec R'$. Hence we obtain a commutative diagram with horizontal isomorphisms
\[
\xymatrix{
\Lie \Pic^0_{\Yh_{R'} / R'} \ar[r]^{\sim} \ar[d]_{\Lie p^n} \ar@/_5pc/[dd]_{\Lie \rho^*_{R'}} & H^1 (\Yh_{R'} , \Oh) \ar[d] \ar@/^5pc/[dd]^{\rho^*_{R'}} \\
\Lie \Pic^0_{\Yh_{R'} / R'} \ar[r]^{\sim} \ar[d]_{\Lie g} & H^1 (\Yh_{R'}, \Oh) \ar[d] \\
\Lie \Pic^0_{\eZ_{R'} / R'} \ar[r]^{\sim} & H^1 (\eZ_{R'} , \Oh)
}
\]
Since $\Lie p^n$ is $p^n$-multiplication, we deduce that 
\[
\rho^*_{R'} (H^1 (\Yh_{R'} , \Oh)) \subseteq p^n H^1 (\eZ_{R'} , \Oh) \; ,
\]
and by flat base change that $\rho^* (H^1 (\Yh, \Oh)) \subset p^n H^1 (\eZ , \Oh)$ which completes the proof.
\end{proof}

Note that in the following theorem and its proof we have changed our usual notation somewhat. 

\begin{theorem}
  \label{cd_t9}
{\bf a} For any smooth projective curve $\oX$ over $\oQ_p$ the category $\eB^{\sharp}_{X_{\C_p}}$, where $X_{\C_p} = \oX \otimes_{\oQ_p} \C_p$, contains all line bundles $L$ of degree zero on $X_{\C_p}$.\\
{\bf b} If $\oX$ has a smooth model over $\oZ_p$, then $\eB_{X_{\C_p}}$ contains all line bundles of degree zero on $X_{\C_p}$. 
\end{theorem}

\begin{proof}
  We may assume that $\oX$ has positive genus. By the semistable reduction theorem there is a finite extension $K$ of $\Q_p$ and a smooth projective curve $X$ over $K$ with $X (K) \neq \emptyset$ together with a semistable model $\eX$ over $\eo_K$ such that $\oX = X \otimes_K \oQ_p$. In particular $\eX$ is cohomologically flat of dimension zero over $\eo_K$. According to \cite{BLR} 9.4 Theorem 1, $\Pic^0_{\eX / \eo_K}$ is a semiabelian scheme over $\eo_K$. Hence $\Pic^0_{\eX / \eo_K} (\eo)$ is an open subgroup of $\Pic^0_{X / K} (\C_p) = \Pic^0 (X_{\C_p})$ the group of isomorphism classes of line bundles on $X_{\C_p}$ of degree zero.

{\bf Claim} If the class of $L$ in $\Pic^0 (X_{\C_p})$ lies in $\Pic^0_{\eX / \eo_K} (\eo)$ then $L$ is in $\eB_{X_{\C_p}}$.

\begin{proofof}
  {\bf the claim} By assumption $L$ is the generic fibre of a line bundle $\Lh$ on $\eX_{\eo}$ giving rise to a class in $\Pic^0_{\eX / \eo_K} (\eo)$. Note that according to \cite{BLR} 8.1 Proposition 4, we have $\Pic (\eX_{\eo}) = \Pic_{\eX / \eo_K} (\eo)$. Now, $\eo_n = \oZ_p / p^n \oZ_p = \varinjlim_{F / K} \eo_F / p^n \eo_F$ where $F$ runs over the finite extensions of $K$ in $\oQ_p$. The rings $\eo_F / p^n \eo_F$ are finite, hence $\Pic^0_{\eX / \eo_K} (\eo_F / p^n \eo_F)$ is a finite group. It follows that 
\[
\Pic^0_{\eX / \eo_K} (\eo_n) = \varinjlim_{F / K} \Pic^0_{\eX / \eo_K} (\eo_F / p^n \eo_F)
\]
is a torsion group. Let $\Lh_n = \Lh \otimes_{\eo} \eo_n$ be the reduction $\mod p^n$ of $\Lh$ to a line bundle on $\eX_n = \eX_{\eo} \otimes_{\eo} \eo_n = \eX \otimes_{\eo_K} \eo_n$. It defines a class in $\Pic^0_{\eX_n / \eo_n} (\eo_n) = \Pic^0_{\eX / \eo_K} (\eo_n)$ which must have finite order. Hence there is some $N \ge 1$ such that $\Lh^{\otimes N}_n \simeq \Oh$. By proposition \ref{cd_t2} applied to $Q = K , R = \eo_K$ and $X , \eX$, there is a finite extension $K \subset K' \subset \oQ_p$ with ring of integers $R' = \eo_{K'}$ and an object $\pi_{R'} : \Yh_{R'} \to \eX_{R'} = \eX \otimes_{\eo_K} R'$ of $\Sh^{ss}_{\eX_{R'}}$ together with a commutative diagram, where $g (0) = 0$:
\[
\xymatrix{
\Pic^0_{\eX_{R'} / R'} \ar[rr]^{\pi^*_{R'}} \ar[dr]_N & & \Pic^0_{\Yh_{R'} / R'} \\
 & \Pic^0_{\eX_{R'} / R'} \ar[ur]_g
}
\]
Moreover we can assume that $Y_{K'} = \Yh_{R'} \otimes_{R'} K'$ has a $K'$-rational point. For the object $\pi = \pi_{R'} \otimes_{R'} \oZ_p : \oYh = \Yh_{R'} \otimes_{R'} \oZ_p \to \oeX = \eX \otimes_{\eo_K} \oZ_p$ of $\Sh_{\oeX}$ we therefore get the commutative diagram
\[
\xymatrix{
\Pic^0_{\eX_n / \eo_n} (\eo_n) \ar[rr]^{\pi^*_n} \ar[dr]_N & & \Pic^0_{\Yh_n / \eo_n} (\eo_n) \\
 & \Pic^0_{\eX_n / \eo_n} (\eo_n) \ar[ur]_{G_n}
}
\]
where $G_n (0) = 0$. Hence we find
\[
\pi^*_n [\Lh_n] = G_n (N [\Lh_n]) = G_n ([\Lh^{\otimes N}_n]) = G_n (0) = 0 \; .
\]
It follows that $\pi^*_n \Lh_n$ is a trivial bundle on $\Yh_n$. Since this construction can be done for every $n \ge 1$ the bundle $\Lh$ belongs to $\eB_{\eX_{\eo}}$ and therefore $L$ is an object of $\eB_{X_{\C_p}}$.
\end{proofof}

We can now proceed with the proof of the theorem. Part {\bf b} follows from the claim for $\eX$ smooth. In order to prove {\bf a}, let $L$ be any line bundle of degree zero on $X_{\C_p}$. By a result of Coleman (Theorem 4.1. in \cite{Co}), the cokernel of the inclusion map
\[
\Pic^0_{\eX / \eo_K} (\eo) \hookrightarrow \Pic^0_{X / K} (\C_p)
\]
is torsion. Hence there exists an integer $N \ge 1$ such that $L^{\otimes N}$ is the generic fibre of some line bundle $\Lh_1$ on $\eX_{\eo}$ giving rise to a class in $\Pic^0_{\eX / \eo_K} (\eo)$. With notations as before, we have for this $N$ a commutative diagram
\[
\xymatrix{
\Pic^0_{\oeX / \oZ_p} \ar[rr]^{\pi^*} \ar[dr]_N & & \Pic^0_{\oYh / \oZ_p} \\
 & \Pic^0_{\oeX / \oZ_p} \ar[ur]_G & 
}
\]
where $G (0) = 0$. Since $\pi_{R'}$ is in $\Sh_{\eX_{R'}}$, the generic fibre $\alpha$ of $\pi : \oYh \to \oeX$ is a finite \'etale covering $\alpha : \oY = \oYh \otimes_{\oZ_p} \oQ_p \to \oX$ of $\oX$ by the smooth projective curve $\oY$. It suffices to show that $\alpha^*_{\C_p} L$ belongs to $\eB_{Y_{\C_p}}$. Under the inclusion
\[
\Pic^0_{\oYh / \oZ_p} (\eo) \hookrightarrow \Pic^0 (Y_{\C_p}) \; ,
\]
the element $G ([\Lh_1])$ is mapped to $\alpha^*_{\C_p} ([L])$. By the claim applied to $\oY$ and the pair $Y_{K'} , \Yh_{R'}$ instead of $\oX$ and $X , \eX$ it follows that $\alpha^*_{\C_p} ([L])$ lies in $\eB_{Y_{\C_p}}$ as was to be shown.
\end{proof}

\begin{rem}
  By the preceding results the category $\eB_{X_{\C_p}}$ contains all unipotent vector bundles on $X_{\C_p}$, i.e. all bundles obtained by successive extensions of the trivial line bundle.

More generally, the category $\eB^{\sharp}_{X_{\C_p}}$ contains all successive extensions of line bundles of degree zero.
\end{rem}

The following insight is due to Faltings without proof in his setting of $p$-adic Higgs bundles \cite{Fa2}. We give a proof below.

\begin{theorem}
  \label{cd_t13neu}
Let $D$ be a divisor on a smooth projective curve $X$ over $\oQ_p$. Then every bundle in $\eB^{\sharp}_{X_{\C_p},D}$ is semistable of degree zero.
\end{theorem}

\begin{proof}
  By the definition of semistability it suffices to show the assertion for every bundle $E'$ in $\eB_{X_{\C_p},D}$. 

We may assume that $E' = \eE' \otimes_{\eo} \C_p$ for a bundle $\eE'$ in $\eB_{\eX_{\eo},D}$ for a model $\eX$ of $X$. By corollary \ref{cd_t3},\,3) there exists an object $\pi : \Yh \to \eX$ of $\Sh^{ss}_{\eX,D}$ such that $\pi^*_1 \eE'_1$ is a trivial bundle on $\Yh_1 = \Yh \otimes \eo / p$, where $\eE'_1 = \eE' \otimes \eo / p$. Since the generic fibre $\pi_{\oQ_p}$ of $\pi$ is finite it suffices to show that $E = \pi^*_{\Q_p} E'$ is semistable of degree zero on $\Yh_{\C_p}$. Setting $\eE = \pi^* \eE'$ we have $E = \eE \otimes_{\eo} \C_p$.

Besides, $\eE_1 = \eE \otimes \eo / p$ is a trivial bundle on $\Yh_1$. We have to show that $E$ has degree zero and that every subbundle $L \subset F$ has degree $\deg L \le 0$. 

Let $K$ be a finite extension of $\Q_p$ such that $\Yh$ descends to a model $\Yh_{\eo_K}$ over $\eo_K$ of its generic fibre $Y$, i.e. $\Yh = \Yh_{\eo_K} \otimes_{\eo_K} \oZ_p$. Since $\Yh_{\eo_K} / \eo_K$ has the same geometric fibres as $\Yh / \oZ_p$ it is also semistable. The scheme $\Yh_{\eo}$ is the projective limit of the semistable $A$-schemes $\Yh_A = \Yh_{\eo_K} \otimes_{\eo_K} A$, where $A$ runs over the finitely generated normal $\eo_K$-subalgebras of $\eo$. Moreover $\Yh_1$ is the projective limit of the schemes $\Yh_{A_1} = \Yh_A \otimes_A A_1$, where $A_1 = A / p A$.

Consider the family $(\Yh_{\eo} , \eE , L \subset F , \eE_1 \overset{\alpha}{\silo} \Oh^r_{\Yh_1})$ where $\alpha$ is some isomorphism of locally free $\Oh_{\Yh_1}$-sheaves. By \cite{EGAIV} (8.5.5), (8.9.1), (8.5.2), (11.2.6) there exists a normal finitely generated $\eo_K$-algebra $A$ in $\eo$ with quotient field $Q (A)$ such that the family descends to a family 
\[
(\Yh_A , \eE_A , L_{Q (A)} \subset \eE_{Q (A)}, \eE_{A_1} \overset{\alpha_{A_1}}{\silo} \Oh^r_{\Yh_{A_1}}) \; , \quad \mbox{where}
\]
\begin{itemize}
\item $\Yh_A$ is a proper semistable curve over $A$
\item $\eE_A$ is a vector bundle on $\Yh_A$ and $\eE_{Q (A)} = \eE_A \otimes_A Q (A)$
\item $L_{Q (A)}$ is a vector bundle on $\Yh_{Q (A)} = \Yh_A \otimes_A Q (A)$ which is a subbundle of $\eE_{Q (A)}$
\item $\alpha_{A_1}$ is an isomorphism of locally free $\Oh_{\Yh_{A_1}}$-modules where $\eE_{A_1} = \eE_A \otimes_A A_1$.
\end{itemize}

We need a prime ideal $\ep$ of $A$ of height one containing the maximal ideal $(\pi_K)$ of $\eo_K$. Since $A \subset \eo$, the special fibre $(\spec A) \otimes \eo_K / \pi_K$ is non-empty. Any prime ideal $\ep$ in $A$ corresponding to the generic point of an irreducible component of $(\spec A) \otimes \eo_K / \pi_K$ will do, cf. \cite{Liu} theorem 4.3.12. Note that $\ep \supset p A$. Since $A$ is normal, $A_{\ep} \subset \C_p$ is a discrete valuation ring containing $\eo_K$. Note that in general $A_{\ep} \not\subset \eo$. 

Let $R$ be the strict henselization of $A_{\ep}$ in the algebraic closure of $Q (A)$ in $\C_p$. Then $R$ is a discrete valuation ring in $\C_p$ with quotient field $Q \subseteq \C_p$ whose residue field $\kappa \supset \eo_K / \pi_K$ is separably closed. Let $(\Yh_R , \eE_R , L_Q \subset \eE_Q)$ be the base change of $(\Yh_A , \eE_A , L_{Q (A)} \subset \eE_{Q (A)})$ via $A \subset R$ resp. $Q (A) \subset Q$. The restriction $\eE_{\kappa}$ of $\eE_R$ to the special fibre $\Yh_{\kappa} = \Yh_R \otimes_R \kappa$ is trivial because $\eE_{A_1}$ is trivial and $A \subset R$ induces a map $A_1 \to R/p \to \kappa$ since $p \in \ep_R$. By Riemann--Roch, $\deg (\eE_Q) = \chi (\eE_Q) - r \chi (\Oh_{\Yh_Q})$ where $r$ is the rank of $\eE$. By \cite{EGAIII}, 7.9.4, the Euler characteristic of vector bundles on $\Yh_R$ is locally constant in the fibres, which implies $\deg \eE_Q = \chi (\eE_{\kappa}) - r \chi (\Oh_{\Yh_{\kappa}}) = 0$. Since $E = \eE_Q \otimes_Q \C_p$, it follows that $\deg E = 0$. 
Similarly, $\deg L = \deg L_Q$. It remains therefore to show that $\deg L_Q \le 0$. Using the next result the theorem follows.
\end{proof}

The proof of the following proposition is due to Raynaud. It replaces a more involved argument in an earlier version of this paper.

\begin{prop}
  \label{cd_t14neu}
Let $R$ be a discrete valuation ring with quotient field $Q$ and separably closed residue field $\kappa$. Let $Z$ be a smooth projective curve over $Q$ with a model $\Zh$ over $R$. Consider a vector bundle $\eE$ on $\Zh$ whose special fibre $\eE_{\kappa}$ is a trivial bundle on $\Zh_{\kappa}$. Then its generic fibre $E = \eE_Q$ is semistable of degree zero.
\end{prop}

\begin{proof}
  By assumption $\det \eE_{\kappa}$ is a trivial line bundle. Hence we have
  \begin{eqnarray*}
    \deg E = \deg \det E & = & \chi (Z , \det \eE_Q) - \chi (Z , \Oh) \\
& = & \chi (\Zh_{\kappa} , \det \eE_{\kappa}) - \chi (\Zh_{\kappa} , \Oh) = 0
  \end{eqnarray*}
since the Euler characteristics are constant in the fibres. 

It suffices to show that for every exact sequence $0 \to E_1 \to E \to E_2 \to 0$ of vector bundles on $Z$ we have $\deg (E_2) \ge 0$. Consider the canonical extension $\Fh_1 \subset \eE$ of $E_1$ in $E$ cf. \cite{EGA1} (9.4.1). For every open subset $U$ of $\Zh$ we have
\[
\Gamma (U , \Fh_1) = \{ s \in \Gamma (U, \eE) \tei s \, |_{U \cap Z} \in \Gamma (U \cap Z , E_1) \} \; .
\]
The sheaf $\Fh_1$ is a coherent, torsion free $\Oh_{\Zh}$-module. Let $\Fh_2 = \eE / \Fh_1$ be the quotient, so that $0 \to \Fh_1 \to \eE \to \Fh_2 \to 0$ is an exact sequence of coherent sheaves on $\Zh$ with generic fibre $0 \to E_1 \to E \to E_2 \to 0$. If $r$ is the rank of $\Fh_2$, we blow up the $r$-th Fitting ideal $\Ih$ of $\Fh_2$ and get a proper morphism
\[
\varphi : \Zh' \longrightarrow \Zh
\]
which is an isomorphism on the generic fibres. 

If we denote by $\Ih'$ the ideal $\varphi^{-1} (\Ih) \cdot \Oh_{\Zh'}$ (which coincides with the $r$-th Fitting ideal of $\varphi^* \Fh_2$), then $\eE_2 = \varphi^* \Fh_2 / \Ann_{\varphi^* \Fh_2} (\Ih')$ is a locally free sheaf on $\Zh'$ by \cite{Ray-Gru}, (5.4.3). Let $\Fh$ be the coherent sheaf on $\Zh'$ such that the sequence 
\[
0 \longrightarrow \Fh \longrightarrow \eE' \longrightarrow \eE_2 \longrightarrow 0
\]
with $\eE' = \varphi^* \eE$ is exact. Since $\varphi$ is an isomorphism on the generic fibre and the generic fibre of $\Ih$ is $\Oh_Z$, the generic fibre of $\eE'$ resp. $\eE_2$ is isomorphic to $E$ resp. $E_2$. 

Now let $C_1 , \ldots , C_r$ be the irreducible components of the special fibre $\Zh'_{\kappa}$, and let $\tC_i \to C_i$ be their normalizations. By $\alpha_i : \tC_i \to C_i \to \Zh'$ we denote the corresponding morphisms. Since $\eE_2$ is locally free, the sequence
\[
0 \longrightarrow \alpha^*_i \Fh \longrightarrow \alpha^*_i \eE' \longrightarrow \alpha^*_i \eE_2 \longrightarrow 0
\]
is exact on $\tC_i$. Since the special fibre $\eE'_{\kappa}$ is trivial, the sheaf $\alpha^*_i \eE'$ is isomorphic to a power of the structure sheaf $\Oh_{\tC_i}$. In particular, it is a semistable sheaf of degree $0$ on the smooth, projective curve $\tC_i$ over $\kappa$.

Therefore, the quotient $\alpha^*_i \eE_2$ has degree $\ge 0$. By the degree formula in \cite{BLR}, 9.1, Proposition 5, it follows for the line bundle $(\det \eE_2)_{\kappa}$, that
\[
\chi (\Zh_{\kappa} , (\det \eE_2)_{\kappa}) - \chi (\Zh_{\kappa} , \Oh_{\Zh_{\kappa}}) \ge 0 \; .
\]
Since the Euler characteristics are constant in the fibres of $\Zh$, we deduce $\deg (E_2) = \deg ((\eE_2)_Q) = \deg \det ((\eE_2)_Q) \ge 0$. Hence $E$ is indeed semistable.
\end{proof}

\begin{rem}
  The indecomposable components $E_i$ of a semistable bundle $E$ of degree zero on $X_{\C_p}$ have degree zero since they are both sub- and quotient bundles of $E$ and hence have $\deg E_i \le 0$ and $\deg E_i \ge 0$. If $X = A$ is an elliptic curve over $\oQ_p$ the converse is true. A vector bundle $E$ on $A_{\C_p}$ is semistable of degree zero if and only if it is the direct sum of indecomposable bundles of degree zero. This follows from the splitting of the Harder--Narasimhan filtration on bundles over elliptic curves.
\end{rem}

By \cite{At}, Theorem 5, p. 432 every indecomposable vector bundle of degree zero on $A_{\C_p}$ is of the form $L \otimes F_r$ where $L$ is a line bundle of degree zero and $F_r$ is an iterated extension of trivial line bundles. Using theorems \ref{cd_t8} and \ref{cd_t9} it follows that $L \otimes F_r$ lies in $\eB^{\sharp}_{A_{\C_p}}$ and in the case where $A$ has good reduction even in $\eB_{A_{\C}}$. From all this one obtains:

\begin{cor}
  \label{cd_t15neu}
Let $A$ be an elliptic curve over $\oQ_p$.\\
{\bf a} The category $\eB^{\sharp}_{A_{\C_p}}$ consists of all semistable bundles of degree zero on $A_{\C_p}$. All of these are successive extensions of line bundles of degree zero.\\
{\bf b} If $A_K$ has good reduction we have in addition $\eB^{\sharp}_{A_{\C_p}} = \eB_{A_{\C_p}}$.
\end{cor}

The following result makes it substantially easier to verify that a vector bundle lies in $\eB_{\eX_{\eo},D}$.

\begin{theorem}
  \label{cd_t16_neu}
Let $\eX$ be a model over $\oZ_p$ of the smooth projective curve $X$ over $\oQ_p$. Let $k = \oF_p$ be the residue field of $\oZ_p$. A vector bundle $\eE$ on $\eX_{\eo}$ lies in $\eB_{\eX_{\eo} , D}$ if and only if there is an object $\pi : \Yh \to \eX$ of $\Sh_{\eX,D}$ such that $\pi^*_k \eE_k$ is a trivial bundle on $\Yh_k = \Yh \otimes_{\oZ_p} k$.
\end{theorem}

\begin{rem}
  In particular, every vector bundle $\eE$ on $\eX_{\eo}$ whose restriction $\eE_k$ to the special fibre $\eX_k$ of $\eX_{\eo}$ is trivial lies in $\eB_{\eX_{\eo}}$. As explained to us by Holger Brenner there exist examples of rank two bundles $\eE$ on smooth models of certain plane algebraic curves $X$ such that $\eE_k$ is trivial and $\eE_{\C_p}$ is stable of degree zero. They are constructed by restricting suitable syzygy bundles on $\Pa^2$. 
\end{rem}

\begin{proof}
  The necessity is clear. Consider a vector bundle $\eE$ on $\eX_{\eo}$ with $\pi^*_k \eE_k$ trivial. We may assume that $\pi$ is in $\Sh^{ss}_{\eX,D}$. The family $(\eX , D, \eE_1 , \pi : \Yh \to \eX)$ descends to a family $(\eX_0 , D_0 , \Fh , \pi_0 : \Yh_0 \to \eX_0)$ over $\eo_K$ for $K$ a finite extension of $\Q_p$. Here $\eX_{\eo}$ is a model of $X_0 = \eX_0 \otimes_{\eo_K} K$ and $\Fh$ is a vector bundle on $\eX_{\eo} \otimes \eo_K / p \eo_K$ whose restriction to the special fibre $\eX_0 \otimes \eo_K / \ep$ becomes trivial after pullback along $\pi_0 \otimes \eo_K / \ep$. Moreover $\pi_0$ is an object of $\Sh^{ss}_{\eX_0 , D_0}$. Let $e$ be the ramification index of $K$ over $\Q_p$ and set $\eo_{\nu / e} = \eo /\ep^{\nu} \eo = \oZ_p / \ep^{\nu} \oZ_p$. Note that this is compatible with our earlier notation $\eo_n = \oZ_p / p^n \oZ_p$. Let $\pi_{\nu / e} , E_{\nu / e}$ etc. be the base change with $\eo_{\nu / e}$. Since $\pi_{1/e}$ is also the base change of $\pi_0 \otimes \eo_K / \ep$ with $\eo_{1/e}$ it follows that $\pi^*_{1/e} \eE_{1/e}$ is trivial on $\Yh_{1/e}$. By induction it therefore suffices to prove the following assertion:\\
Given $\nu \ge 2$ and some $\pi : \Yh \to \eX$ in $\Sh^{ss}_{\eX , D}$ with $\pi^*_{(\nu-1)/e} \eE_{(\nu-1) / e}$ trivial, there exists an object $\mu : \Zh \to \eX$ in $\Sh^{ss}_{\eX,D}$ with $\mu^*_{\nu/e} \eE_{\nu/e}$ trivial on $\Zh_{\nu/e}$.

Consider the closed immersion $i : \Yh_{(\nu-1) / e} \hookrightarrow \Yh_{\nu/e}$ and set
\[
\Jh = \Imm (\omega^{\nu-1} : \Oh_{\Yh_{\nu/e}} \longrightarrow \Oh_{\Yh_{\nu/e}}) \; .
\]
Here $\omega$ is a prime element in $\eo_K$. Let $r$ be the rank of $\eE$, then we have a short exact sequence of (Zariski-)sheaves of groups on $\Yh_{\nu/e}$:
\[
0 \longrightarrow M_r (\Jh) \xrightarrow{f} \GL_r (\Oh_{\Yh_{\nu/e}}) \xrightarrow{\mathrm{adj}} i_* \GL_r (\Oh_{\Yh_{(\nu-1) / e}}) \longrightarrow 1 \; .
\]
Here $\adj$ is the adjunction map and $f (A) := 1 + A$. Observe that $f$ is a homomorphism, $f (A+A') = f (A) f (A')$ since $AA' = 0$ in $M_r (\Jh)$ because $\Jh^2 = 0$. Right exactness follows because $\GL_r$ is formally smooth over $\Z$. We obtain an exact sequence of pointed sets:
\[
H^1 (\Yh_{\nu/e} , M_r (\Jh)) \xrightarrow{f} H^1 (\Yh_{\nu/e} , \GL_r (\Oh)) \xrightarrow{i^*} H^1 (\Yh_{(\nu-1) / e} , \GL_r (\Oh)) \; .
\]
Exactness can be checked directly. Alternatively one may identify sheaf torsors for the {\it affine} group scheme $\GL_r$ with vector bundles and quote \cite{Gi} III Proposition 3.3.1 for the non-abelian cohomology sequence and \cite{Gi} V Proposition 3.1.3 for the isomorphism
\[
H^1 (\Yh_{\nu/e} , i_* \GL_r (\Oh)) = H^1 (\Yh_{(\nu-1) / e} , \GL_r (\Oh)) \; .
\]
Note here that for elementary reasons we have $R^1 i_* \GL_r (\Oh) = 0$.

Consider the class $\Omega$ of $\pi^*_{\nu/e} \eE_{\nu/e}$ in $H^1 (\Yh_{\nu/e} , \GL_r (\Oh))$. Via $i^*$ it is mapped to the class of $i^* \pi^*_{\nu/e} \eE_{\nu/e} = \pi^*_{(\nu-1) / e} \eE_{(\nu-1) / e}$ i.e. to the trivial class in $H^1 (\Yh_{(\nu-1) / e} , \GL_r (\Oh))$. Hence $\Omega$ is of the form $\Omega = f (A)$ for some class $A = (A_{kl})$ in
\[
H^1 (\Yh_{\nu/e} , M_r (\Jh)) = M_r (H^1 (\Yh_{\nu/e}, \Jh)) \; .
\]
Instead of recalling the argument from non-abelian cohomology we could also have quoted \cite{Gi} VII Th\'eor\`eme 1.3.1 for this conclusion.

The exact sequence on $\Yh_{\nu/e}$:
\[
0 \longrightarrow \Ker \omega^{\nu-1} \longrightarrow \Oh \xrightarrow{g} \Jh \longrightarrow 0
\]
where $g$ is multiplication by $\omega^{\nu-1}$ gives a surjection:
\[
H^1 (\Yh_{\nu/e} , \Oh) \overset{g}{\twoheadrightarrow} H^1 (\Yh_{\nu/e} , \Jh)
\]
because $\Yh_{\nu/e}$ is one-dimensional.

Hence we have $\Omega = fg(B)$ for some matrix $B = (B_{kl})$ with entries in $H^1 (\Yh_{\nu/e}, \Oh)$. If the genus of $Y$ is zero, the same argument as in the proof of theorem \ref{cd_t8} shows that $H^1 (\Yh_{\nu/e}, \Oh) = 0$ and we are done. If the genus of $Y$ is non-zero it was shown in the proof of theorem \ref{cd_t8} that there is a morphism $\rho : \eZ \to \Yh$ in $\Sh_{\eX,D}$ such that $\rho^* : H^1 (\Yh , \Oh) \to H^1 (\eZ , \Oh)$ satisfies $\rho^* (H^1 (\Yh , \Oh)) \subset p^{\nu} H^1 (\eZ , \Oh)$. By corollary \ref{cd_t3},\,3 we may assume that the object $\mu : \eZ \to \eX$ is even in $\Sh^{ss}_{\eX,D}$. Arguing as in the proof of theorem \ref{cd_t8} (reduction step for the {\bf claim}, with $p$ and $\eo_n$ replaced by $\omega$ and $\eo_{\nu/e}$) one sees that the induced map
\[
\rho^*_{\nu/e} : H^1 (\Yh_{\nu/e} , \Oh) \longrightarrow H^1 (\eZ_{\nu/e} , \Oh)
\]
is trivial. The commutative diagram
\[
\xymatrix{
H^1 (\Yh_{\nu/e} , M_r (\Oh)) \ar[r]^{fg} \ar[d]^{\rho^*_{\nu/e} = 0} & H^1 (\Yh_{\nu/e}, \GL_r (\Oh)) \ar[d]^{\rho^*_{\nu/e}} \\
H^1 (\eZ_{\nu/e} , M_r (\Oh)) \ar[r]^{fg} & H^1 (\eZ_{\nu/e} , \GL_r (\Oh))
}
\]
shows that $\rho^*_{\nu/e} \Omega$ is the trivial class. Hence
\[
\mu^*_{\nu/e} \eE_{\nu/e} = \rho^*_{\nu/e} (\pi^*_{\nu/e} \eE_{\nu/e})
\]
is a trivial bundle on $\eZ_{\nu/e}$, as was to be shown.
%
\end{proof}

\begin{rem}
  The proof shows that a vector bundle $\eE$ on $\eX_{\eo}$ lies in $\eB_{\eX_{\eo}}$ if the special fibre $\eE_k$ is trivial. In this case, for each $n \ge 1$ there is a trivializing cover $\pi$ in $\Sh_{\eX_{\eo}}$ whose generic fibre is a Galois covering of $X$ with solvable Galois group.
\end{rem}

\begin{definition}
  Let $R$ be a valuation ring with quotient field $Q$ and residue field $k$. Consider a model $\eX / R$ of a smooth projective curve $X / Q$ and let $\eE$ be a vector bundle on $\eX$. We say that $\eE$ has strongly semistable reduction of degree zero if the pullback of $\eE_k$ to the normalization $\tC$ of each irreducible component $C$ (with the reduced structure) of $\eX_k$ is strongly semistable of degree zero. Note here that each $\tC$ is a smooth projective curve over $k$.
\end{definition}

The following theorem is one of our main results.

\begin{theorem}
 \label{cd_t16_x}
Let $\eX / \oZ_p$ be a model of the smooth projective curve $X / \oQ_p$. Let $\eE$ be a vector bundle on $\eX_{\eo}$. Then $\eE$ belongs to $\eB_{\eX_{\eo} , D}$ for some divisor $D$ on $X$ if and only if $\eE$ has strongly semistable reduction of degree zero. In this case $\eE$ even belongs to $\eB_{\eX_{\eo} , D}$ and $\eB_{\eX_{\eo} , \tD}$ for two divisors $D$ and $\tD$ on $X$ with disjoint support.
\end{theorem}

The proof depends on the following result which for smooth projective curves is due to Lange and Stuhler \cite{LS} 1.9 Satz.

\begin{theorem}
  \label{cd_t17_x}
Let $E$ be a vector bundle on a purely one-dimensional proper scheme $X$ over $\F_q$. Then the following conditions are equivalent: \\
a) The pullback of $E$ to the normalization of each irreducible component of $X$ is strongly semistable of degree zero.\\
b) There is a finite surjective morphism $\varphi : Y \to X$ where $Y$ is a purely one-dimensional proper scheme over $\F_q$ such that $\varphi^* E$ is a trivial bundle.\\
c) Same as in b) but with $\varphi$ a composition $\varphi : Y \xrightarrow{F^s} Y \xrightarrow{\pi} X$ for some $s \ge 0$ where $\pi$ is finite \'etale and surjective and $F = \Fr_q = \Fr^r_p$ is the $q = p^r$-linear Frobenius on $Y$.
\end{theorem}

\begin{proof}
  If b) holds then every irreducible component $C$ of $X$ is finitely dominated by an irreducible component $D$ of $Y$. It follows that the pullback of $E$ to $\tC$ is trivialized by the finite surjective morphism $\tD \to \tC$. Since semistability can be verified after pullback to a finite covering and since the absolute Frobenius is functorial, assertion a) follows. 

It remains to show that a) implies c). There are only finitely many isomorphism classes of semistable vector bundles of degree zero on a smooth projective curve over a finite field. It follows that there are only finitely many isomorphism classes of vector bundles $E$ on $X$ whose pullbacks to the normalizations of the irreducible components of $X$ are semistable of degree zero. To see this, we first assume that $X$ is reduced. Let $X = \bigcup C_{\nu}$ be the decomposition of $X$ into its irreducible components and let $\pi : \tX = \coprod \tC_{\nu} \to X$ be the finite normalization morphism. Generalizing the arguments in the proofs of \cite{BLR}, Ch. 9, Propositions 9 and 10 or \cite{Liu} Lemma 7.5.12 one sees the following: The cokernel of the natural injection of sheaves of groups $\GL_r (\Oh_X) \to \pi_* \GL_r (\Oh_{\tX})$ is a skyscraper sheaf of sets $\prod_{x \in X^{\sing}} i_{x*} S_x$ where each set $S_x$ is finite. Using \cite{Gi} III Proposition 3.2.2 we obtain a non-abelian cohomology sequence
\[
\prod_{x \in X^{\sing}} S_x \longrightarrow H^1 (X, \GL_r (\Oh)) \longrightarrow H^1 (X , \pi_* \GL_r (\Oh_{\tX})) = \prod_{\nu} H^1 (\tC_{\nu} , \GL_r (\Oh)) \; .
\]
Here we have also used \cite{Gi} V Proposition 3.1.3 and the equation $R^1 \pi_* \GL_r (\Oh) = 1$ which follows because vector bundles are locally trivial. 
Using \cite{Gi} III Corollaire 3.2.4, it follows that there are only finitely many isomorphism classes of vector bundles on $X$ which induce given isomorphism classes of vector bundles on the curves $\tC_{\nu}$. 

If $X$ is not reduced, we have to show that the map
\[
H^1 (X , \GL_r (\Oh)) \longrightarrow H^1 (X^{\red} , \GL_r (\Oh))
\]
has finite fibres. By devissage it suffices to show that for every ideal $\Jh \subset \Oh_X$ with $\Jh^2 = 0$, the map
\[
\varphi : H^1 (X , \GL_r  (\Oh)) \longrightarrow H^1 (X' , \GL_r (\Oh))
\]
has finite fibres where $i : X' \hookrightarrow X$ is the closed subscheme of $X$ defined by $\Jh$. As in the proof of theorem \ref{cd_t16_neu}, we have a non-abelian cohomology sequence
\[
H^1 (X , M_r (\Jh)) \xrightarrow{f} H^1 (X , \GL_r (\Oh)) \xrightarrow{\varphi} H^1 (X' , \GL_r (\Oh)) \; .
\]
It follows that $\varphi$ has finite fibres.

Now, if we are given a vector bundle $E$ as in a) the pullbacks to $\tC_{\nu}$ of all the vector bundles $F^{*n}_X E$ on $X$ are semistable of degree zero. It follows that we have $F^{*s}_X E = F^{*t}_X E$ for some integers $t > s \ge 0$. For the bundle $E' = F^{s*}_X E$ we therefore have $F^{r*}_X E' = E'$ where $r = t-s \ge 1$. Now, the proof of \cite{LS} 1.4 Satz extends without change to an arbitrary $\F_q$-scheme (note that in \cite{LS} the proof that $\pi$ is finite is omitted, but this is not difficult). This shows that there exists a finite \'etale and surjective morphism $\pi : Y \to X$ such that $\pi^* E' = \pi^* F^{*s}_X E$ is a trivial bundle. With $X$, the scheme $Y$ is a purely one-dimensional proper $\F_q$-scheme as well. It follows that $(\pi \verk F^s_Y)^* E = (F^s_X \verk \pi)^* E = \pi^* F^{*s}_X E$ is a trivial bundle as was to be shown for c).
\end{proof}

\begin{proofof}
  {\bf theorem \ref{cd_t16_x}} For a vector bundle $\eE$ in $\eB_{\eX_{\eo}, D}$ choose a cover $\pi : \Yh \to \eX$ in $\eS^{\good}_{\eX, D}$ such that $\pi^*_k \eE_k$ is a trivial bundle. Let $\eX_k = \bigcup_{\nu} C_{\nu}$ be the decomposition of $\eX_k$ into irreducible components. Since $\eX$ is irreducible and $\pi (\Yh)$ is closed and contains the generic point of $\eX$, the map $\pi$ is surjective. Therefore any $C_{\nu}$ is finitely dominated by an irreducible component of $\Yh_k$. As above it follows that the pullbacks of $\eE_k$ to the $\tC_{\nu}$ are strongly semistable of degree zero. 

Now assume that the vector bundle $\eE$ on $\eX_{\eo}$ has strongly semistable reduction of degree zero. There is a finite extension $K$ of $\Q_p$ with ring of integers $\eo_K$ and residue field $\kappa \simeq \F_q$ such that the family $(X , \eX , C_{\nu} , \eE_k)$ descends to a family $(X_K , \eX_{\eo_K} , C_{\nu 0} , \eE_0)$ with corresponding properties. In particular $\eE_0$ is a vector bundle on the special fibre $\eX_0 = \eX_{\eo_K} \otimes \kappa$ whose pullbacks to the normalizations $\tC_{\nu 0}$ of the irreducible components $C_{\nu 0}$ of $\eX_0$ are strongly semistable of degree zero. Using theorem \ref{cd_t17_x} we obtain a finite \'etale morphism $\tpi_0 : \tYh_0 \to \eX_0$ such that for the composition $\tvarphi_0 : \tYh_0 \xrightarrow{F^s} \tYh_0 \xrightarrow{\tpi_0} \eX_0$ the pullback bundle $\tvarphi^*_0 \eE_0$ is trivial. Note that in this statement we may replace $s$ by any integer $s' \ge s$ and hence $F$ by any power of $F$. Next, using \cite{SGA1} IX th\'eor\`eme 1.10 we may lift $\tpi_0 : \tYh_0 \to \eX_0$ to a finite \'etale morphism $\tpi_{\eo_K} : \tYh_{\eo_K} \to \eX_{\eo_K}$ whose special fibre is $\tpi_0$. After replacing $K$ by a finite extension and performing a base extension to the new $\eo_K$, theorem \ref{cd_t1}, 5) allows us to dominate $\tpi_{\eo_K}$ by an object $\pi_{\eo_K} : \Yh_{\eo_K} \to \eX_{\eo_K}$ of $\Sh^{ss}_{\eX_{\eo_K}}$. By Lipman's desingularization theorem we may assume that $\Yh_{\eo_K}$ besides being semistable is also regular c.f. \cite{Liu} 10.3.25 and 10.3.26. Replacing $F$ by $\Fr^r_p$ where $q = p^r$ now denotes the order of the new residue field it follows that under the composition $\varphi_0 : \Yh_0 \xrightarrow{F^s} \Yh_0 \xrightarrow{\pi_0} \eX_0$ the pullback $\varphi^*_0 \eE_0$ is a trivial bundle. 


The irreducible regular surface $\Yh_{\eo_K}$ is proper and flat over $\eo_K$. Hence by a theorem of Lichtenbaum \cite{Lich} there exists a closed immersion $\Yh_{\eo_K} \hookrightarrow \Pa^N_{\eo_K}$ over $\eo_K$. Let $H_i$ be the coordinate hyperplane $x_i = 0$ in $\Pa^N_K$, and put $\Delta = \bigcup^N_{i=0} H_i$. Then $\Pa^N_K \ohne \Delta = \Ge^N_{m,K}$. We observe that for any finite set $S$ of closed points in $\Pa^N_K$ there is a linear isomorphism $f \in \PGL_N (\eo_K)$ of $\Pa^N_{\eo_K}$ such that its generic fibre $f_K$ maps $S$ to $\Pa^N_K \ohne \Delta$. Hence we can choose a closed immersion $\tau : \Yh_{\eo_K} \hookrightarrow \Pa^N_{\eo_K}$ in such a way that $Y_K$ is not contained in $\Delta$. Consider the finite morphism $F_{\eo_K} : \Pa^N_{\eo_K} \to \Pa^N_{\eo_K}$ given on $A$-valued points by mapping $[x_0 : \ldots : x_N]$ to $[x^q_0 : \ldots : x^q_N]$ for any $\eo_K$-algebra $A$. Over $\Ge^N_{m , K} = \Pa^N_K \ohne \Delta$ this morphism is \'etale. Define an $\eo_K$-scheme $\Yh'_{\eo_K}$ by the cartesian diagram:
\[
\xymatrix{
\Yh'_{\eo_K} \ar[r]^{ \rho_{\scriptstyle\eo_K}} \ar@{^{(}->}[d] & \Yh_{\eo_K} \ar@{^{(}->}[d]^{\tau} \\
\Pa^N_{\eo_K} \ar[r]^{ F^s_{\scriptstyle\eo_K}} & \Pa^N_{\eo_K}
}
\]
Then $\rho_{\eo_K}$ is finite and $\rho_K = \rho_{\eo_K} \otimes K : Y'_K \to Y_K$ is \'etale over $U_K = Y_K \cap \Ge^N_{m,K}$. Let $D'_K$ be a divisor on $Y_K$ whose support is $Y_K \ohne U_K$. Let $\rho_0 = \rho_{\eo_K} \otimes \kappa : \Yh'_0 \to \Yh_0$ be the special fibre of $\rho_{\eo_K}$. The reduction of $F^s_{\eo_K}$ is $F^s$ i.e. the $rs$-th power of the absolute Frobenius morphism on $\Pa^N_{\kappa}$. Define a morphism $i : \Yh_0 \to \Yh'_0$ over $\kappa$ by the commutative diagram
\[
\xymatrix{
\Yh_0 \ar@/^1pc/[rrd]^{F^s} \ar[dr]^i \ar@/_1pc/[ddr]_{\tau_0}  & & \\
 & \Yh'_0 \ar[r]^{\rho_0} \ar@{^{(}->}[d] & \Yh_0 \ar@{^{(}->}[d]^{\tau_0}\\
 & \Pa^N_{\kappa} \ar[r]^{F^s} & \Pa^N_{\kappa}
}
\]
Lemma \ref{cd_t18_neu} below implies that $i$ induces an isomorphism $\Yh_0 \silo \Yh'^{\red}_0$. Set $D_K = \pi_K (D'_K)$. Base extending the situation to $\eo_{\oK} = \oZ_p$ we obtain an object $\pi' : \Yh' \xrightarrow{\rho} \Yh \xrightarrow{\pi} \eX$ of $\Sh_{\eX,D}$. Moreover $\eE_k$ is trivialized by pullback via the composed map $\Yh_k \xrightarrow{i_k} \Yh'_k \xrightarrow{\pi'_k} \eX_k$ since we have $\pi'_k \verk i_k = \pi_k \verk (\rho_k \verk i_k) = \pi_k \verk (F^s \otimes_{\kappa} k)$ and $\eE_k = \eE_0 \otimes_{\kappa} k$. In addition $i_k$ induces an isomorphism of $\Yh_k$ onto $\Yh'^{\red}_k$. For this, note that $\Yh'^{\red}_k = (\Yh'_0 \otimes_{\kappa} k)^{\red} = \Yh'^{\red}_0 \otimes_{\kappa} k$ since $\Yh'^{\red}_0 \cong \Yh_0$ is geometrically reduced. By corollary \ref{cd_t3}, 3) there is an object $\mu : \Zh \to \eX$ of $\Sh^{ss}_{\eX,D}$ such that $\mu$ factors over $\pi'$:
\[
\mu : \Zh \xrightarrow{\psi} \Yh' \xrightarrow{\pi'} \eX \; .
\]
The special fibre of $\Zh$ is reduced. Hence the morphism $\psi_k : \Zh_k \to \Yh'_k$ factors over $i_k : \Yh_k \cong \Yh'^{\red}_k \to \Yh'_k$ and therefore $\mu_k$ factors over $\pi'_k \verk i_k$. It follows that $\mu^*_k \eE_k$ is trivial. Applying theorem \ref{cd_t16_neu} it follows that $\eE$ is an object of $\eB_{\eX_{\eo} , D}$. Let $\tau : \Yh_{\eo_K} \hookrightarrow \Pa^N_{\eo_K}$ be the projective embedding above. By the above observation on linear automorphisms, after changing $\tau$ by some $f \in \PGL_N (\eo_K)$ we can assume that $\tau_K$ maps the support of $\pi^*_K D_K$ to $\Pa^N_K \ohne \Delta = \Ge^N_{m,K}$. Let $\tD_K$ be a divisor on $Y_K$ with support equal to $Y_K \ohne Y_K \cap \tau^{-1}_K (\Ge^N_{m,K})$. Then we have seen that $\eE$ is in $\eB_{\eX_{\eo} , \tD}$ where $\tD$ is the base-change of $\pi_K (\tD_K)$. By construction, $\tD_K$ is disjoint from  $\pi^*_K D_K$ and hence $\tD$ is disjoint from $D$.
\end{proofof}

\begin{lemma}
  \label{cd_t18_neu}
Let $T$ be an $\F_p$-scheme and let $\tau : S \hookrightarrow T$ be a closed immersion of a reduced subscheme $S$ of $T$. For an integer $N \ge 1$ consider the canonical diagram where the square is cartesian:
\[
\xymatrix{
S \ar@/^1pc/[rrd]^{\Fr^N_p} \ar[dr]^i \ar@/_1pc/[ddr]_{\tau}  & & \\
 & S' \ar[r]^{\rho} \ar[d] & S \ar[d]^{\tau}\\
 & T \ar[r]^{\Fr^N_p} & T
}
\]
Then the induced map $i : S \to S'^{\red}$ is an isomorphism.
\end{lemma}

\begin{proof}
  We may assume that $T = \spec R$ is affine. Then we have $S = \spec R / \ea$ for an ideal $\ea$ with $\ea = \sqrt{\ea}$ and $S' = \spec R / \eb$ where $\eb$ is the ideal generated by all elements $r^{p^N}$ with $r \in \ea$. The homomorphism $i^{\sharp} : R / \eb \to R / \ea$ is given by $i^{\sharp} (r \bmod \eb) = r \bmod \ea$. It is immediate that $\sqrt{\eb} = \ea$. Hence $i : \spec R / \ea \to (\spec R / \eb)^{\red} = \spec (R / \sqrt{\eb})$ is an isomorphism.
\end{proof}

The following result due to M. Raynaud improves theorem \ref{cd_t16_x} in the case of good reduction. The proof is a modification of the argument for theorem \ref{cd_t16_x}.

\begin{theorem}
  \label{cd_t20_neu}
Let $\eX / \oZ_p$ be a smooth model of a smooth projective curve $X / \oQ_p$ of nonzero genus and let $\eE$ be a vector bundle on $\eX_{\eo}$. Then $\eE$ belongs to $\eB_{\eX_{\eo}}$ if and only if $\eE_k$ is strongly semistable of degree zero on the smooth projective curve $\eX_k$ over $k$.
\end{theorem}

\begin{proof}
  Assume that $\eE_k$ is strongly semistable of degree zero. As in the proof of theorem \ref{cd_t16_x} we descend $(X, \eX , \eE_k)$ to a family $(X_K , \eX_{\eo_K} , \eE_0)$ for a finite extension $K / \Q_p$ with residue field $\kappa \cong \F_q , q = p^r$. 

Since $\eE_0$ is strongly semistable of degree zero on the smooth projective curve $\eX_0 = \eX \otimes \kappa$ over $\kappa$, theorem \ref{cd_t17_x} or in fact the original result in \cite{LS} 1.9 Satz provides us with a smooth projective curve $\Yh_0$ over $\kappa$ and a composition $\varphi_0 : \Yh_0 \xrightarrow{F^s} \Yh_0 \xrightarrow{\pi_0} \eX_0$ with $s \ge 0$ and $\pi_0$ finite \'etale such that the bundle $\varphi^*_0 \eE_0$ is trivial. As before we may lift $\pi_0$ to a finite \'etale morphism $\pi_{\eo_K} : \Yh_{\eo_K} \to \eX_{\eo_K}$. Then $\Yh_{\eo_K}$ is a smooth and proper irreducible $\eo_K$-scheme. As in the proof of theorem \ref{cd_t16_x} we can replace $K$ by a finite extension and hence assume that we have a section $y \in \Yh_{\eo_K} (\eo_K) = Y_K (K)$. Set $\Bh = \Pic^0_{\Yh_{\eo_K} / \eo_K}$ and consider the Albanese map
\[
\tau : \Yh_{\eo_K} \hookrightarrow \Ah = \hat{\Bh}
\]
with $\tau (y) = 0$. Define $\Yh''_{\eo_K}$ by the cartesian diagram
\begin{equation}
  \label{eq:?}
  \xymatrix{
\Yh''_{\eo_K} \ar[r]^{\lambda} \ar@{^{(}->}[d] & \Yh_{\eo_K} \ar@{^{(}->}[d]^{\tau} \\
\Ah \ar[r]^{q^s} & \Ah
}
\end{equation}
After reduction, the $q^s$-multiplication map on $\Ah_0$ factors
\[
q^s : \Ah_0 \xrightarrow{V^s} \Ah_0 \xrightarrow{F^s} \Ah_0
\]
with $V$ the $r$-th power of Verschiebung and $F = \Fr^r_p$. Correspondingly the reduction of diagram (\ref{eq:?}) factors into two cartesian diagrams:
\[
\xymatrix{
\Yh''_0 \ar[r] \ar@{^{(}->}[d] & \Yh'_0 \ar[r] \ar@{^{(}->}[d] & \Yh_0 \ar@{^{(}->}[d]^{\tau} \\
\Ah_0 \ar[r]^{V^s} & \Ah_0 \ar[r]^{F^s} & \Ah_0
}
\]
By lemma \ref{cd_t18_neu} the diagram
\[
\xymatrix{
\Yh_0 \ar@/^1pc/[rrd]^{F^s} \ar[dr]^i \ar@/_1pc/[ddr]_{\tau}  & & \\
 & \Yh'_0 \ar[r] \ar@{^{(}->}[d] & \Yh_0 \ar@{^{(}->}[d]^{\tau}\\
 & \Ah_0 \ar[r]^{F^s} & \Ah_0
}
\]
induces an isomorphism $i : \Yh_0 \to \Yh'^{\red}_0$. Base extending to $\oZ_p$ resp. $k$ we can dominate $\Yh'' = \Yh''_{\eo_K} \otimes \oZ_p$ by an object $\mu : \eZ \to \eX$ of $\Sh^{ss}_{\eX}$. Since $\eZ_k$ is reduced, the reduction $\mu_k$ factors over $\Yh'^{\red}_k \hookrightarrow \Yh'_k \to \Yh_k$ and hence over $F^s \otimes k : \Yh_k \to \Yh_k$. Hence $\mu^*_k \eE_k$ is a trivial bundle and we conclude using theorem \ref{cd_t16_neu}.
\end{proof}

\section{\'Etale parallel transport for vector bundles in $\eB$}

For vector bundles in $\eB_{\eX_{\eo},D} , \eB_{X_{\C_p},D}$ and $\eB^{\sharp}_{X_{\C_p},D}$ we will now construct canonical isomorphisms of parallel transport along \'etale paths between geometric points of $X \ohne D$. We begin by recalling some facts about the fundamental groupoid. The general reference is \cite{SGA1}.

Let $Z$ be a variety over $\oQ_p$ and choose a geometric point $z$ in $Z (\C_p)$. Let $F_z$ be the functor from the category of finite \'etale coverings $Z'$ of $Z$ to the category of finite sets defined by $F_z = \Mor_Z (z, \_)$. It attaches to $Z'$ the set of $\C_p$-valued points of $Z'$ lying over $z$. The functor $F_z$ is known to be strictly pro-representable: There is a projective system $\tZ = (Z_i , z_i , \phi_{ij})_{i \in I}$ of pointed Galois coverings of $Z$ where $I$ is a directed set, and the $z_i \in Z_i (\C_p)$ are points over $z$. Moreover, for $i \ge j$ the map $\phi_{ij} : Z_i \to Z_j$ is an epimorphism over $Z$ such that $\phi_{ij} (z_i) = z_j$ and such that the natural map
\[
\varinjlim_i \Mor_Z (Z_i , Z') \longrightarrow F_z (Z')
\]
induced by evaluation on the $z_i$'s is a bijection for every $Z'$.

For our purposes, we define the \'etale fundamental groupoid $\Pi_1 (Z)$ of $Z$ as a topological category, as follows: The set of objects of $\Pi_1 (Z)$ is $Z (\C_p)$. For two $\C_p$-valued points $z$ and $z^*$ of $Z$ set
\begin{equation}
  \label{eq:cd4}
  \Mor_{\Pi_1 (Z)} (z, z^*) = \Iso (F_z , F_{z^*}) \; .
\end{equation}
Such an isomorphism of fibre functors will be called an \'etale path (up to homotopy) from $z$ to $z^*$. Using the pro-representability of $F_z$ and $F_{z^*}$, one sees that $\Mor_{\Pi_1 (Z)} (z , z^*)$ is a pro-finite set and as such a compact totally disconnected Hausdorff space. Moreover, composition of morphisms gives a {\it continuous} map
\[
\Mor_{\Pi_1 (Z)} (z , z^*) \times \Mor_{\Pi_1 (Z)} (z^* , z^{**}) \longrightarrow \Mor_{\Pi_1 (Z)} (z , z^{**}) \; .
\]
The \'etale fundamental group with base point $z$ is the profinite group
\[
\pi_1 (Z, z) = \Mor_{\Pi_1 (Z)} (z,z) = \Aut (F) \quad \mbox{where} \; F = F_z \; .
\]
There is an isomorphism of topological groups
\begin{equation}
  \label{eq:cd5}
  \pi_1 (Z,z) \silo \left( \varprojlim_i \Aut_Z (Z_i) \right)^{\op} \; .
\end{equation}
Here the natural transformation $\sigma_F : F \silo F$ given by the family of compatible bijections $\sigma_{F (Z_i)} : F (Z_i) \silo F (Z_i)$ for $i \in I$ is sent to the projective system $(\sigma_i)_{i \in I}$ where $\sigma_i \in \Aut_Z (Z_i)$ is uniquely determined by the relation:
\[
\sigma_i (z_i) = \sigma_{F (Z_i)} (z_i) \; .
\]
Let $\Rep_{\Pi_1 (Z)} (\eo)$ resp. $\Rep_{\Pi_1 (Z)} (\C_p)$ be the $\eo$-linear resp. $\C_p$-linear categories of continuous functors from $\Pi_1 (Z)$ into the category of free $\eo$-modules of finite rank resp. the category of finite dimensional $\C_p$-vector spaces. Here a functor between topological categories is called continuous if the induced maps between the topological spaces of morphisms are continuous.

We now make some remarks on the functoriality of $\Pi_1$.

Let $\alpha : Z_1 \to Z_2$ be a morphism of varieties over $\oQ_p$. There is an induced continuous functor $\alpha_* : \Pi_1 (Z_1) \to \Pi_1 (Z_2)$ defined as follows. On objects $\alpha_*$ is the map $\alpha : Z_1 (\C_p) \to Z_2 (\C_p)$. For points $z , z'$ of $Z_1 (\C_p)$ it remains to define continuous maps
\[
\alpha_* : \Iso (F_z , F_{z'}) \longrightarrow \Iso (F_{\alpha (z)} , F_{\alpha (z')}) \; .
\]
For a finite \'etale morphism $Y_2 \to Z_2$ consider the base change \\
$Y_1 = Y_2 \times_{Z_2} Z_1 \to Z_1$. There are natural bijections
\[
F_z (Y_1) \cong F_{\alpha (z)} (Y_2) \quad \mbox{and} \quad F_{z'} (Y_1) \cong F_{\alpha (z')} (Y_2) \; .
\]
For $\gamma \in \Iso (F_z , F_{z'})$ define $\alpha_* (\gamma) (Y_2)$ as the composition:
\[
\alpha_* (\gamma) (Y_2) : F_{\alpha (z)} (Y_2) \cong F_z (Y_1) \overset{\gamma (Y_1)}{\cong} F_{z'} (Y_1) \cong F_{\alpha (z')} (Y_2) \; .
\]
This defines an isomorphism $\alpha_* (\gamma)$ of fibre functors. By construction, the map $\gamma \mapsto \alpha_* (\gamma)$ is continuous. It is clear that $\alpha_*$ defined on objects and morphisms gives a functor.

For a second morphism $\beta : Z_2 \to Z_3$ of varieties over $\oQ_p$ we find that $(\beta \verk \alpha)_* = \beta_* \verk \alpha_* : \Pi_1 (Z_1) \to \Pi_1 (Z_3)$. Obviously $\id_* = \id$. 

Now we consider the effect of Galois conjugation on fundamental groupoids. For a scheme $Y$ over $\oQ_p$ and an automorphism $\sigma$ of $\oQ_p$ over $\Q_p$ set $^{\sigma} Y = Y \otimes_{\oQ_p , \sigma} \oQ_p$ and write $\sigma : Y \silo \hsigma Y$ for the inverse of the projection map.

We now define a continuous functor $\sigma_* : \Pi_1 (Z) \to \Pi_1 (\,\!^{\sigma} Z)$. On objects, $\sigma_*$ is defined by mapping $z \in Z (\C_p)$ to $^{\sigma} z = \sigma \verk z \verk \sigma^{-1} = \sigma \verk z \verk \spec \sigma$ in $^{\sigma} Z (\C_p)$.

The continuous map $\sigma_*$ between the spaces of morphisms
\[
\sigma_* : \Iso (F_z , F_{z'}) \silo \Iso (F_{\,\!^{\sigma}\!z} , F_{\,\!^{\sigma}\!z'})
\]
is obtained as follows: Every finite \'etale cover of $^{\sigma} Z$ is of the form $^{\sigma} Y$ for a finite \'etale cover $Y$ of $Z$. It is clear that $F_{\,\!^{\sigma}\! z} (\,\!^{\sigma}\!Y) \cong F_z (Y)$ naturally for every point $z$ of $Z (\C_p)$. Define $\sigma_* (\gamma) (\hsigma Y)$ as the composition:
\[
\sigma_* (\gamma) (\hsigma Y) : F_{\hsigma z} (\hsigma Y) \cong F_z (Y) \xrightarrow{\gamma} F_{z'} (Y) \cong F_{\hsigma z'} (\hsigma Y) \; .
\]
This defines an isomorphism of fibre functors $\sigma_* (\gamma)$. The map $\gamma \mapsto \sigma_* (\gamma)$ is continuous. The maps $\sigma_*$ on objects and morphisms define a functor $\sigma_*$. It is clear that we have $(\sigma \tau)_* = \sigma_* \verk \tau_*$ as functors from $\Pi_1 (Z)$ to $\Pi_1 (\,\!^{\sigma \tau} Z) = \Pi_1 (\hsigma (\htau Z))$.

If $Z$ is already defined over an extension $K \subset \oQ_p$ of $\Q_p$ i.e. $Z = Z_K \otimes_K \oQ_p$ for some variety $Z_K$ over $K$ then for every $\sigma \in G_K = \Gal (\oQ_p / K)$ the map $\id \times_{\spec K} \spec (\sigma^{-1})$ gives a $\oQ_p$-linear isomorphism $\hsigma Z \silo Z$. This will be used to identify $\hsigma Z$ with $Z$. It follows that for such $Z$ the group $G_K$ acts from the left by continuous automorphisms on the category $\Pi_1 (Z)$. 

For a topological group $\Sigma$ let $\Rep_{\Sigma} (\eo)$ be the category of continuous representations of $\Sigma$ on free $\eo$-modules of finite rank. We define $\Rep_{\Sigma} (\C_p)$ similarly.

All these categories are equipped with a tensor product, duals, internal homs and exterior powers. The $\eo$-categories are exact, the $\C_p$-categories are even abelian.

\begin{lemma}
  \label{cd_t11}
For a variety $Z$ as above and a fixed point $z_0 \in Z (\C_p)$ the natural forgetful functors
\[
\Rep_{\Pi_1 (Z)} (\eo) \longrightarrow \Rep_{\pi_1 (Z, z_0)} (\eo) \quad \mbox{and} \quad \Rep_{\Pi_1 (Z)} (\C_p) \longrightarrow \Rep_{\pi_1 (Z, z_0)} (\C_p)
\]
are fully faithful.
\end{lemma}

\begin{proof}
  Since $Z$ is connected all objects of $\Pi_1 (Z)$ are isomorphic to each other. Faithfulness follows. Given representations $V$ and $V'$ of $\Pi_1 (Z)$ let us write $V_z = V (z)$ and $V'_z = V' (z)$. Given a $\pi_1 (Z , z_0)$-equivariant homomorphism $f_{z_0} : V_{z_0} \to V'_{z_0}$ define $f_z : V_z \to V'_z$ for arbitrary $z \in Z (\C_p)$ as follows. Choose an \'etale path $\gamma \in \Mor_{\Pi_1 (Z)} (z,z_0)$ and set $f_z = V' (\gamma)^{-1} \verk f_{z_0} \verk V (\gamma)$. This is independent of $\gamma$ since $f_{z_0}$ is $\Aut_{\Pi_1 (Z)} (z_0 , z_0)$-equivariant. One checks that the family of homomorphisms $(f_z)_{z \in Z (\C_p)}$ defines a morphism of functors from $V$ to $V'$ which induces $f_{z_0}$. Hence the above forgetful functors are full.
\end{proof}

Consider as before a smooth projective curve $X$ over $\oQ_p$, a divisor $D$ in $X$ and a model $\eX$ of $X$ over $\oZ_p$. Set $U = X \ohne D$.

Given a bundle $\eE$ in $\eB_{\eX_{\eo},D}$, we will construct a continuous functor $\rho_{\eE}$ from $\Pi_1 (U)$ into the category of free $\eo$-modules of finite rank. By properness $X (\C_p) = \eX_{\eo} (\eo)$. Hence we may view any geometric point $x \in X (\C_p)$ as a section $x_{\eo} : \spec \eo \to \eX_{\eo}$ over $\spec \eo$. We write $\eE_{x_{\eo}} = x^*_{\eo} \eE$ viewed as a free $\eo$-module of $\rank \, r = \rank \, \eE$. The reduction $\eX_{\eo} (\eo) \to \eX_{\eo} (\eo_n) = \eX_n (\eo_n)$ maps $x_{\eo}$ to a morphism 
\[
x_n : \spec \eo_n \to \spec \eo \xrightarrow{x_{\eo}} \eX_{\eo}
\]
and we set $\eE_{x_n} = x^*_n \eE = \eE_{x_{\eo}} \otimes_{\eo} \eo_n$ viewed as a free $\eo_n$-module of rank $r$. We have
\[
\eE_{x_{\eo}} = \varprojlim_n \eE_{x_n}
\]
as topological $\eo$-modules, the topology on $\eE_{x_n}$ being the discrete one. 
We define $\rho_{\eE}$ on the set of objects $U (\C_p)$ of $\Pi_1 (U)$ by setting $\rho_{\eE} (x) = \eE_{x_{\eo}}$. It remains to define continuous maps:
\[
\rho_{\eE} : \Mor_{\Pi_1 (U)} (x , x') = \Iso (F_x , F_{x'}) \longrightarrow \Hom_{\eo} (\eE_{x_{\eo}} , \eE_{x'_{\eo}}) \; .
\]
These in turn will be obtained as the projective limit of maps
\[
\rho_{\eE,n} : \Iso (F_x , F_{x'}) \longrightarrow \Hom_{\eo_n} (\eE_{x_n} , \eE_{x'_n})
\]
for $n \ge 1$. By construction each map $\rho_{\eE,n}$ will factor over a finite quotient of the pro-finite set $\Iso (F_x , F_{x'})$. Hence each $\rho_{\eE,n}$ is continuous and therefore $\rho_{\eE}$ will be continuous as well.

Now, given $\gamma$ in $\Iso (F_x , F_{x'})$ and some $n \ge 1$, let us construct $\rho_{\eE,n} (\gamma)$. By definition of $\eB_{\eX_{\eo},D}$ and by Corollary \ref{cd_t3}, 3) there exists an object $\pi : \Yh \to \eX$ of $\Sh^{\good}_{\eX,D}$ such that $\pi^*_n \eE_n$ is a trivial bundle. Set $Y = \Yh \otimes \oQ_p$ and $V = Y \ohne \pi^* D$. Then $V \to U$ is a finite \'etale covering. Choose a point $y \in V (\C_p)$ above $x$ and let $y' = \gamma y \in V (\C_p)$ be the image of $y$ under the map
\[
\gamma_V : F_x (V) \to F_{x'} (V) \; .
\]
Then $y'$ lies over $x'$. Since the structural morphism $\lambda :\Yh \to \spec \oZ_p$ satisfies $\lambda_* \Oh_{\Yh} = \Oh_{\spec \oZ_p}$ universally, we find $\lambda_{n^*} \Oh_{\Yh_n} = \Oh_{\spec \eo_n}$ and therefore the pullback map under $y_n : \spec \eo_n \to \Yh_n$ is an isomorphism:
\[
y^*_n : \Gamma (\Yh_n , \pi^*_n \eE_n) \silo \Gamma (\spec \eo_n , y^*_n \pi^*_n \eE_n) = \Gamma (\spec \eo_n , x^*_n \eE_n) = \eE_{x_n} \; .
\]
We can now define $\rho_{\eE,n} (\gamma)$ to be the composition:
\[
\rho_{\eE,n} (\gamma) = \gamma (y)^*_n \verk (y^*_n)^{-1} = y'^*_n \verk (y^*_n)^{-1} : \eE_{x_n} \longrightarrow \eE_{x'_n} \; .
\]
Note that by construction $\rho_{\eE,n}$ factors over the finite set $\Iso (F_x (V) , F_{x'} (V))$.

\begin{theorem}
  \label{cd_t12}
The preceeding constructions are independent of all choices and define a continuous functor $\rho_{\eE}$ from $\Pi_1 (X \ohne D)$ into the category of free $\eo$-modules of finite rank.
\end{theorem}

\begin{proof}
  We first check that $\rho_{\eE,n} (\gamma)$ does not depend on the choice of the point $y$ above $x$. So let $z$ be another point in $V (\C_p)$ over $x$. By theorem \ref{cd_t4} there are a finite group $G$ and a $G$-equivariant morphism $\tpi : \tYh \to \eX$ defining an object of $\eT^{\good}_{\eX,D}$, together with a morphism $\varphi : \tYh \to \Yh$ with $\tpi = \pi \verk \varphi$. In particular $\tilde{V} = \tY \ohne \tpi^* D$ is a Galois covering of $U$ with group $G$. Here $\tY$ is the generic fibre of $\tYh$. Choose points $\ty$ and $\tz$ in $\tilde{V} (\C_p)$ above $y$ resp. $z$. Then the points $\gamma \ty$ and $\gamma \tz$ lie above $\gamma y$ resp. $\gamma z$. Since $\ty$ and $\tz$ both lie above $x$, there is a unique $\sigma$ in $G$ with $\sigma \ty = \tz$ and hence with $\sigma \ty_{\eo} = \tz_{\eo}$ and $\sigma \ty_n = \tz_n$ as well. By construction the following diagram is commutative:
\[
\xymatrix{
\eE_{x_n} \ar@{=}[d] & \Gamma (\Yh_n , \pi^*_n \eE_n) \ar[l]_-{\overset{y^*_n}{\sim}} \ar[d]^{\wr \, \varphi^*_n} \ar[r]^-{\overset{(\gamma y)^*_n}{\sim}} & \eE_{x'_n} \ar@{=}[d] \\
\eE_{x_n} & \Gamma (\tYh_n , \tpi^*_n \eE_n)\ar[l]_-{\overset{\ty^*_n}{\sim}} \ar[r]^-{\overset{(\gamma \ty)^*_n}{\sim}} & \eE_{x'_n}
}
\]
Hence we have the formula:
\[
(\gamma y)^*_n \verk (y^*_n)^{-1} = (\gamma \ty)^*_n \verk (\ty^*_n)^{-1}
\]
and similarly 
\[
(\gamma z)^*_n \verk (z^*_n)^{-1} = (\gamma \tz)^*_n \verk (\tz^*_n)^{-1} \; .
\]
Now, $\tz = \sigma \verk \ty$ implies that $\tz^*_n$ equals the composition
\[
\Gamma (\tYh_n , \tpi^*_n \eE_n) \xrightarrow{\sigma^*} \Gamma (\tYh_n , \tpi^*_n \eE_n) \xrightarrow{\ty^*_n} \eE_{x_n} \; .
\]
By naturality of $\gamma$ we have $\gamma \tz = \sigma \verk \gamma \ty$ and as before $(\gamma \tz)^*_n = (\gamma \ty)^*_n \verk \sigma^*$. Thus we find
\[
(\gamma z)^*_n \verk (z^*_n)^{-1} = (\gamma \tz)^*_n \verk (\tz^*_n)^{-1} = (\gamma \ty)^*_n \verk \sigma^* \verk (\sigma^*)^{-1} \verk (\ty^*_n)^{-1} = (\gamma y)^*_n \verk (y^*_n)^{-1} \; .
\]
Now we prove that $\rho_{\eE,n} (\gamma)$ does not depend on the trivializing cover $\pi : \Yh \to \eX$. So, let $\tpi : \tYh \to \eX$ be another object of $\Sh^{\good}_{\eX,D}$ such that $\tpi^*_n \eE_n$ is a trivial bundle.

By Corollary \ref{cd_t3}, 3) we may assume that there is a morphism $\varphi : \tYh \to \Yh$ with $\tpi = \pi \verk \varphi$. With notations as above choose a point $\ty \in \tilde{V} (\C_p)$ above $x$ and set $y = \varphi_{\oQ_p} (\ty)$ where $\varphi_{\oQ_p} : \tY \to Y$ is the induced map on the generic fibres. It follows that $\varphi_{\oQ_p} (\gamma \ty) = \gamma y$ and by properness of $\Yh$ and $\tYh$ over $\spec \oZ_p$ that $y_{\eo} = \varphi (\ty_{\eo})$ and $(\gamma y)_{\eo} = \varphi ((\gamma \ty)_{\eo})$. One obtains the same diagram as above.

Hence we have
\[
(\gamma y)^*_n \verk (y^*_n)^{-1} = (\gamma \ty)^*_n \verk (\ty^*_n)^{-1}
\]
and this implies that $\rho_{\eE,n} (\gamma)$ does not depend on the trivializing good cover. Hence $\rho_{\eE,n} (\gamma)$ is well defined.

It is clear that we have $\rho_{\eE,n} (\id) = \id$ for the trivial path $\id \in \Iso (F_x , F_x)$. For paths $\gamma \in \Iso (F_x , F_{x'})$ and $\gamma' \in \Iso (F_{x'} , F_{x''})$, choosing a point $y \in V (\C_p)$ over $x$, the point $\gamma y$ lies over $x'$ and hence we have
\[
\rho_{\eE,n} (\gamma) = (\gamma y)^*_n \verk  (y^*_n)^{-1} \quad \mbox{and} \quad \rho_{\eE,n} (\gamma') = (\gamma' (\gamma y))^*_n \verk (\gamma y)^{*-1}_n \; .
\]
This implies the equation:
\begin{equation}
  \label{eq:cd6}
  \rho_{\eE,n} (\gamma') \verk \rho_{\eE,n} (\gamma) = ((\gamma' \verk \gamma) (y))^*_n \verk (y^*_n)^{-1} = \rho_{\eE,n} (\gamma' \verk \gamma) \; .
\end{equation}
We now check that the maps
\[
\rho_{\eE,n} : \Iso (F_x , F_{x'}) \longrightarrow \Hom_{\eo_n} (\eE_{x_n} , \eE_{x'_n})
\]
form a projective system with respect to the natural projections
\[
\lambda_{n+1} : \Hom_{\eo_{n+1}} (\eE_{x_{n+1}} , \eE_{x'_{n+1}}) \longrightarrow \Hom_{\eo_{n+1}} (\eE_{x_{n+1}} , \eE_{x'_{n+1}}) \otimes_{\eo_{n+1}} \eo_n = \Hom_{\eo_n} (\eE_{x_n} , \eE_{x'_n}) 
\]
i.e. that $\lambda_{n+1} \verk \rho_{\eE,n+1} = \rho_{\eE,n}$. 

For a given $n \ge 1$ choose $\pi : \Yh \to \eX$ in $\Sh^{\good}_{\eX,D}$ such that $\pi^*_{n+1} \eE_{n+1}$ is a trivial bundle. Then $\pi^*_n \eE_n$ is trivial as well. For $y$ in $V (\C_p)$ over $x$ and $\gamma \in \Iso (F_x , F_{x'})$ consider the commutative diagram where $a$ and $b$  are the natural maps
\[
\xymatrix{
\spec \eo_n \ar[r]^-{y_n} \ar[d]_b & \Yh_n \ar[d]_a & \spec \eo_n \ar[l]_-{(\gamma y)_n} \ar[d]^b \\
\spec \eo_{n+1} \ar[r]^-{y_{n+1}} & \Yh_{n+1}  & \spec \eo_{n+1} \ar[l]_-{(\gamma y)_{n+1}}
}
\]
It induces a commutative diagram
\[
\xymatrix{
\eE_{x_{n+1}} \ar[d]_{b^*} & \Gamma (\Yh_{n+1} , \pi^*_{n+1} \eE_{n+1}) \ar[l]_-{\overset{y^*_{n+1}}{\sim}} \ar[d]^{a^*} \ar[r]^-{\overset{(\gamma y)^*_{n+1}}{\sim}} & \eE_{x'_{n+1}} \ar[d]^{b^*} \\
\eE_{x_n} & \Gamma (\Yh_n , \pi^*_n \eE_n) \ar[l]_-{\overset{y^*_n}{\sim}} \ar[r]^-{\overset{(\gamma y)^*_n}{\sim}} & \eE_{x'_n}
}
\]
The maps $b^*$ are just the natural reduction maps from the $\eo_{n+1}$-module $\eE_{x_{n+1}}$ resp. $\eE_{x'_{n+1}}$ to the $\eo_n$-module $\eE_{x_n} = \eE_{x_{n+1}} \otimes_{\eo_{n+1}} \eo_n$ resp. $\eE_{x'_n } = \eE_{x'_{n+1}} \otimes_{\eo_{n+1}} \eo_n$. Hence the map $\rho_{\eE,n} (\gamma) = (\gamma y)^*_n \verk (y^*_n)^{-1}$ is the reduction $\mod p^n$ of the map $\rho_{\eE , n+1} (\gamma) = (\gamma y)^*_{n+1} \verk (y^*_{n+1})^{-1}$. 

Let
\[
\rho_{\eE} : \Iso (F_x , F_{x'}) \longrightarrow \Hom_{\eo} (\eE_{x_{\eo}} , \eE_{x'_{\eo}})
\]
be the projective limit of the maps $\rho_{\eE,n}$. Using (\ref{eq:cd6}) it follows that together with the previously defined map $\rho_{\eE}$ on objects, we obtain a continuous functor from $\Pi_1 (X \ohne D)$ into the category of free $\eo$-modules of finite rank.
\end{proof}

For a fixed $\C_p$-valued point $x$ of $X \ohne D$ the continuous functor $\rho_{\eE}$ induces in particular a continuous representation
\[
\rho_{\eE} : \pi_1 (X \ohne D , x) = \Aut_{\Pi_1 (X \ohne D)} (x) \longrightarrow \Aut_{\eo} (\eE_{x_{\eo}}) \; .
\]
In a preliminary version (\cite{De-We}) we defined a representation $\rho_{\eE}$ in the following way: Choose a $G$-equivariant morphism $\pi : \Yh \to \eX$ in $\eT^{\good}_{\eX}$ such that $\pi^*_n \eE_n$ is a trivial bundle. The choice of a point $y \in Y (\C_p)$ above $x$ determines a homomorphism
\[
\pi_1 (X , x) \xrightarrow{\varphi_y} \Aut^{\op}_X Y = G^{\op} \longrightarrow \Aut_{\eo_n} \Gamma (\Yh_n , \pi^*_n \eE_n) \; ,
\]
i.e. a left action of $\pi_1 (X , x)$ on $\Gamma (\Yh_n , \pi^*_n \eE_n)$. Transporting this action to $\eE_{x_n}$ via the isomorphism
\[
\Gamma (\Yh_n , \pi^*_n \eE_n) \xrightarrow{\overset{y^*_n}{\sim}} \eE_{x_n}
\]
gives a representation $\trho_{E,n}$ of $\pi_1 (X , x)$ on $\eE_{x_n}$. The projective limit of the $\trho_{\eE,n}$ defines a representation $\trho_E : \pi_1 (X,x) \to \Aut_{\eo} (\eE_{x_0})$. 

\begin{prop}
  \label{cd_t13}
The representations $\rho_E$ and $\trho_E : \pi_1 (X , x) \to \Aut_{\eo} (\eE_{x_{\eo}})$ agree with each other.
\end{prop}

\begin{proof}
  The present construction obtains $\rho_{\eE} (\gamma)$ as the limit of $\rho_{\eE,n} (\gamma)$, where $\rho_{\eE,n} (\gamma)$ is the composition 
\[
\eE_{x_n} \xrightarrow{\overset{(y^*_n)^{-1}}{\sim}} \Gamma (\Yh_n , \pi^*_n \eE_n) \xrightarrow{(\gamma y)^*_n} \eE_{x_n} \; .
\]
We now show that $\trho_{\eE,n} (\gamma) = \rho_{\eE,n} (\gamma)$. For this, note that under the map $\varphi_y$ the natural transformation $\gamma$ is sent to the unique automorphism $\sigma \in G$ of $Y$ which sends $y$ to $\gamma y$. Hence we have
\[
\trho_{\eE,n} (\gamma) = y^*_n \verk \sigma^* \verk (y^*_n)^{-1} = (\sigma y)^*_n \verk (y^*_n)^{-1} = (\gamma y)^*_n \verk (y^*_n)^{-1} = \rho_{\eE,n} (\gamma) \; .
\]
\end{proof}

We now turn the map $\eE \mapsto \rho_{\eE}$ into a functor $\rho$ from $\eB_{\eX_{\eo},D}$ into $\Rep_{\Pi_1 (X \ohne D)} (\eo)$. Let $f : \eE \to \eE'$ be a morphism in $\eB_{\eX_{\eo},D}$. We claim that the family of $\eo$-module homomorphisms
\[
f_{x_{\eo}} = x^*_{\eo} f : \eE_{x_{\eo}} \longrightarrow \eE'_{x_{\eo}} \quad \mbox{for all} \; x \in U (\C_p) = \Ob \Pi_1 (U)
\]
defines a natural transformation, denoted by $\rho_f$ from $\rho_{\eE}$ to $\rho_{\eE'}$. So, let $\gamma \in \Mor_{\Pi_1 (U)} (x , x')$ be an \'etale path. For a given $n \ge 1$ there is an object $\pi : \Yh \to \eX$ in $\Sh^{\good}_{\eX,D}$ such that both $\pi^*_n \eE_n$ and $\pi^*_n \eE'_n$ are trivial bundles. This follows from Corollary \ref{cd_t3}, 3). Let $f_n$ be the reduction of $f \mod p^n$ and set $f_{x_n} = x^*_n (f)$. Choose a point $y$ above $x$ and set $y' = \gamma y$. Then the commutative diagram
\[
\xymatrix{
\eE_{x_n} \ar[rr]^{f_{x_n}} & & \eE'_{x_n} \\
\Gamma (\Yh_n , \pi^*_n \eE_n) \ar[u]^{y^*_n}_{\wr} \ar[rr]^-{\Gamma (\Yh_n , \pi^*_n f_n)} \ar[d]^{\wr}_{y'^{*}_n} & & \Gamma (\Yh_n , \pi^*_n \eE'_n) \ar[u]_{\wr\, y^*_n} \ar[d]^{\wr \, y'^*_n} \\
\eE_{x'_n} \ar[rr]^-{f_{x'_n}} & & \eE'_{x'_n}
}
\]
shows that we have $f_{x'_n} \verk \rho_{E,n} (\gamma) = \rho_{E',n} (\gamma) \verk f_{x_n}$. In the limit we obtain that $f_{x'_{\eo}} \verk \rho_{\eE} (\gamma) = \rho_{\eE'} (\gamma) \verk f_{x_{\eo}}$. Hence $\rho_f = (f_{x_{\eo}})$ is a morphism from $\rho_{\eE}$ to $\rho_{\eE'}$. It is clear that in this way we obtain a functor $\rho = \rho^{\eX}$. The proof of the following proposition is easy:

\begin{prop}
  \label{cd_t14}
The functor $\rho = \rho^{\eX} : \eB_{\eX_{\eo},D} \to \Rep_{\Pi_1 (X \ohne D)} (\eo)$ is $\eo$-linear and commutes with tensor products, duals, internal homs and exterior powers of vector bundles. Exact sequences of bundles are mapped to exact sequences of representations of $\Pi_1 (X \ohne D)$. 
\end{prop}

We now describe the effect of Galois conjugation on $\rho$. Consider an automorphism $\sigma$ of $\oQ_p$ over $\Q_p$. It induces a $\sigma$-linear functor $\sigma_*$ from $\uVec_{\eX_{\eo}}$ to $\uVec_{\hsigma \eX_{\eo}}$. Here $\hsigma \eX = \eX \otimes_{\oZ_p , \sigma} \oZ_p$ and hence $(\hsigma \eX)_{\eo} = \eX_{\eo} \otimes_{\eo, \sigma} \eo = \hsigma \eX_{\eo}$. The functor $\sigma_*$ sends the vector bundle $\eE$ over $\eX_{\eo}$ to the vector bundle $\hsigma \eE = \eE \otimes_{\eo , \sigma} \eo$ over $\hsigma \eX_{\eo}$. A morphism $f : \eE_1 \to \eE_2$ is sent to $\hsigma f : \hsigma \eE_1 \to \hsigma \eE_2$. On the other hand we have a $\sigma$-linear functor:
\[
\bC_{\sigma} : \Rep_{\Pi_1 (U)} (\eo) \longrightarrow \Rep_{\Pi_1 (\hsigma U)} (\eo) \; .
\]
It is obtained as follows. Let $\Mod_{\eo}$ be the category of free $\eo$-modules of finite rank. We define a continuous $\sigma$-linear functor
\[
\sigma_* : \Mod_{\eo} \silo \Mod_{\eo}
\]
by mapping an $\eo$-module $\Gamma$ to $\sigma_* (\Gamma) = \hsigma \Gamma$, which is $\Gamma$ as a set but with the twisted $\eo$-module structure $\lambda \cdot \gamma = \sigma^{-1} (\lambda) \gamma$ for $\lambda \in \eo$ and $\gamma \in \Gamma$ i.e. $^{\sigma} \Gamma = \Gamma \otimes_{\eo , \sigma} \eo$. We write the identity map $\Gamma \xrightarrow{\id} \hsigma \Gamma$ as $\sigma : \Gamma \to \hsigma \Gamma$ since it is $\sigma$-linear. An $\eo$-module homomorphism $f : \Gamma_1 \to \Gamma_2$ is sent to $\hsigma f = f$ in the first description of $\hsigma \Gamma$ and to $f \otimes_{\eo , \sigma} \eo$ in the second.

On objects of $\Rep_{\Pi_1 (U)} (\eo)$, i.e. on continuous functors $\Gamma : \Pi_1 (U) \to \Mod_{\eo}$, the functor $\bC_{\sigma}$ is defined by setting $\bC_{\sigma} (\Gamma) = \sigma_* \verk \Gamma \verk (\sigma_*)^{-1}$ where \\
$\sigma_* : \Pi_1 (U) \silo \Pi_1 (\hsigma U)$ is the isomorphism of categories recalled above. For a morphism $f : \Gamma \to \Gamma'$ in $\Rep_{\Pi_1 (U)} (\eo)$, i.e. a family of $\eo$-module homomorphisms $f_x : \Gamma_x \to \Gamma'_x$ for $x \in U (\C_p)$ with $\Gamma' (\gamma) \verk f_x = f_{x'} \verk \Gamma (\gamma)$ for all $\gamma \in \Mor_{\Pi_1 (U)} (x,x')$ we define $\bC_{\sigma} (f) : \bC_{\sigma} (\Gamma) \to \bC_{\sigma} (\Gamma')$ as follows. Every point of $\hsigma U (\C_p)$ is of the form $\hsigma x = \sigma_* (x)$ for some point $x$ of $U (\C_p)$. Hence we have to define an $\eo$-linear map
\[
\bC_{\sigma} (f)_{\hsigma x} : \bC_{\sigma} (\Gamma)_{\hsigma x} = (\sigma_* \verk \Gamma) (x) \longrightarrow (\sigma_* \verk \Gamma') (x) = \bC_{\sigma} (\Gamma')_{\hsigma x}
\]
for every $x \in U (\C_p)$ i.e. a map
\[
\bC_{\sigma} (f)_{\hsigma x} : \hsigma \Gamma_x \longrightarrow \hsigma \Gamma'_x \; .
\]
In the above notation we set $\bC_{\sigma} (f)_{\hsigma x} = \sigma \verk f_x \verk \sigma^{-1}$. The family $(\bC_{\sigma} (f)_{\hsigma x})$ defines the desired natural transformation $\bC_{\sigma} (f)$ and $\bC_{\sigma}$ becomes a functor which is easily checked to be $\eo$-linear. 

Moreover we have $\bC_{\tau\sigma} = \bC_{\tau} \verk \bC_{\sigma}$ and $\bC_{\id} = \id$ in an obvious sense.

With trivial changes we also get analogous functors $\sigma_* : \eB_{X_{\C_p},D} \to \eB_{\hsigma X_{\C_p}, \hsigma D}$ and $\sigma_* : \uVec_{\C_p} \to \uVec_{\C_p}$ and $\bC_{\sigma} : \Rep_{\Pi_1 (U)} (\C_p) \to \Rep_{\Pi_1 (\hsigma U)} (\C_p)$. 

The proof of the following proposition is routine.

\begin{prop}
  \label{cd_t15}
In the above situation the diagram of categories and functors is commutative (up to canonical isomorphisms of functors):
\[
\xymatrix{
\eB_{\eX_{\eo},D} \ar[r]^-{\rho^{\eX}} \ar[d]_{\sigma_*} & \Rep_{\Pi_1 (U)} (\eo) \ar[d]^{\bC_{\sigma}} \\
\eB_{\hsigma \eX_{\eo}, \hsigma D} \ar[r]^-{\rho^{\hsigma \eX}} & \Rep_{\Pi_1 (\hsigma U)} (\eo)
}
\]
In particular, we have for $\eE$ in $\eB_{\eX_{\eo},D}$ that
\[
\rho_{\hsigma \eE} = \sigma_* \verk \rho_{\eE} \verk (\sigma_*)^{-1}
\]
as functors from $\Pi_1 (\hsigma U)$ to $\Mod_{\eo}$.
\end{prop}

\begin{rem}
  It also follows that if $\eX$ and $D$ are defined over $\eo_K$, so that $(\hsigma \eX , \hsigma D)$ can be identified with $(\eX, D)$ for all $\sigma \in G_K$, the functor
\[
\rho : \eB_{\eX_{\eo} ,D} \longrightarrow \Rep_{\Pi_1 (U)} (\eo)
\]
commutes with the left $G_K$-actions on these categories defined by letting $\sigma$ act via $\sigma_*$ resp. via $\bC_{\sigma}$.
\end{rem}

The next type of functoriality will be used all the time later. Let $\alpha : \eX \to \eX'$ be a morphism over $\oZ_p$ of models and let $D'$ be a divisor on $X'$. Set $U' = X' \ohne D'$ and $U = X \ohne a^* D'$. The generic fibre of $\alpha$ induces a functor
\[
A (\alpha) : \Rep_{\Pi_1 (U')} (\eo) \longrightarrow \Rep_{\Pi_1 (U)} (\eo)
\]
as follows: For an object $\Gamma$ of $\Rep_{\Pi_1 (U')} (\eo)$ we define $A (\alpha) (\Gamma)$ to be the composed functor:
\[
A (\alpha) (\Gamma) : \Pi_1 (U) \xrightarrow{\alpha_*} \Pi_1 (U') \xrightarrow{V} \Mod_{\eo} \; .
\]
For a morphism $f : \Gamma_1 \to \Gamma_2$ in $\Rep_{\Pi_1 (U')} (\eo)$ given by a family of $\eo$-linear maps $f_{x'} : \Gamma_{1x'} \to \Gamma_{2x'}$ for $x' \in U' (\C_p)$ we define $A (\alpha) (f)$ to be the family of maps
\[
A (\alpha) (\Gamma_1)_x = \Gamma_{1 , \alpha (x)} \xrightarrow{f_{\alpha (x)}} \Gamma_{2 , \alpha (x)} = A (\alpha) (\Gamma_2)_x \; .
\]
It is clear that $A (\alpha)$ so defined gives a functor and that for a second map $\alpha' : \eX' \to \eX''$ we have $A (\alpha' \verk \alpha) = A (\alpha) \verk A (\alpha')$.

\begin{prop}
  \label{cd_t16}
For a morphism $\alpha : \eX \to \eX'$ as above the pullback along $\alpha^*$ induces a functor $\alpha^* : \eB_{\eX'_{\eo},D'} \to \eB_{\eX_{\eo}, \alpha^* D'}$ and the following diagram of categories and functors commutes (up to canonical isomorphisms):
\begin{equation}
  \label{eq:cd7}
\vcenter{\xymatrix{
\eB_{\eX'_{\eo},D'} \ar[r]^-{\rho} \ar[d]_{\alpha^*} & \Rep_{\Pi_1 (U')} (\eo) \ar[d]^{A (\alpha)} \\
\eB_{\eX_{\eo},\alpha^* D'} \ar[r]^-{\rho} & \Rep_{\Pi_1 (U)} (\eo)
}}  
\end{equation}
In particular, for every $\eE$ in $\eB_{\eX'_{\eo},D'}$ we have
\begin{equation}
  \label{eq:cd8}
  \rho_{\alpha^* \eE} = \rho_{\eE} \verk \alpha_*
\end{equation}
as functors from $\Pi_1 (U)$ to $\Mod_{\eo}$.
\end{prop}

\begin{proof}
  Let $\eE$ be a vector bundle in $\eB_{\eX'_{\eo} , D'}$. By proposition \ref{cd_t6}, $\alpha^* \eE$ lies in $\eB_{\eX_{\eo},\alpha^* D'}$. We have $(\rho \verk \alpha^*) (\eE) = \rho_{\alpha^* \eE}$ and $(A (\alpha) \verk \rho) (\eE) = \rho_{\eE} \verk \alpha_*$. Commutativity of (\ref{eq:cd7}) on objects is therefore equivalent to (\ref{eq:cd8}). On objects, relation (\ref{eq:cd8}) simply says that $(\alpha^* \eE)_x = \eE_{\alpha (x)}$, a canonical isomorphism. For $\gamma \in \Mor_{\Pi_1 (U)} (x, z)$ it suffices to show that for every $n \ge 1$, we have
  \begin{equation}
    \label{eq:cd9}
    \rho_{\alpha^* \eE,n} (\gamma) = \rho_{\eE,n} (\alpha_* (\gamma)) \; .
  \end{equation}
Let $\pi' : \Yh' \to \eX'$ be an object of $\Sh^{\good}_{\eX', D'}$ such that $\pi'^*_n \eE_n$ is trivial. Choose some $\pi : \Yh \to \eX$ in $\Sh^{\good}_{\eX,\alpha^*D'}$ covering the object $\tpi : \tYh = \Yh' \times_{\eX'} \eX \to \eX$ of $\Sh_{\eX , \alpha^* D'}$, so that we get a commutative diagram
\[
\xymatrix{
\Yh \ar[r]^{\psi} \ar[d]_{\pi} & \Yh' \ar[d]^{\pi'} \\
\eX \ar[r]^{\alpha} & \eX'
}
\]
Let $y$ be a point in $V (\C_p)$ above $x$ and set $y' = \psi (y)$, a point in $V' (\C_p)$ above $\alpha (x)$. Now, $\rho_{\alpha^* \eE,n}$ is the composition
\[
(\alpha^* \eE)_{x_n} \xrightarrow{(y^*_n)^{-1}} \Gamma (\Yh_n , \pi^*_n \alpha^*_n \eE_n) \xrightarrow{(\gamma y)^*_n} (\alpha^* \eE)_{z_n}
\]
and $\rho_{\eE,n} (\alpha_* (\gamma))$ is the composition
\[
\eE_{\alpha (x)_n} \xrightarrow{(y'^*_n)^{-1}} \Gamma (\Yh'_n , \pi'^*_n \eE_n) \xrightarrow{(\alpha_* (\gamma) y')^*_n} \eE_{\alpha (z)_n} \; .
\]
Hence, for (\ref{eq:cd8}) it suffices to show that the following diagram commutes:
\[
\xymatrix{
(\alpha^* \eE)_{x_n} \ar@{=}[d] & \Gamma (\Yh_n , \pi^*_n \alpha^*_n \eE_n) \ar[l]_-{y^*_n} \ar[r]^-{(\gamma y)^*_n} & (\alpha^* \eE)_{z_n}  \ar@{=}[d] \\
\eE_{\alpha (x)_n} & \Gamma (\Yh'_n , \pi'^*_n \eE_n) \ar[l]_-{y'^*_n} \ar[r]^-{(\alpha_* (\gamma) y')^*_n} \ar[u]_{\psi^*_n} & \eE_{\alpha (z)_n}
}
\]
For the left square this follows from the relation $y' = \psi \verk y$ as morphisms from $\spec \C_p$ to $V' \subset \Yh'$. Similarly the right square commutes because we have $\alpha_* (\gamma) y' = \psi \verk (\gamma y)$. Namely, factoring $\psi$ as a composition $\psi : \Yh \xrightarrow{b} \tYh \xrightarrow{a} \Yh'$ and setting $\ty = b (y)$, we have $\alpha_* (\gamma) (y') = \alpha_* (\gamma) (a (\ty)) = a (\gamma \ty) = a (\gamma (by)) = (a \verk b) (\gamma y) = \psi (\gamma y)$. It is an immediate consequence of the definitions, that diagram (\ref{eq:cd7}) commutes for morphisms i.e. that $A (\alpha) \rho_f = \rho_{\alpha^* f}$ for all $f : \eE_1 \to \eE_2$ in $\eB_{\eX'_{\eo},D'}$.
 \end{proof}

We can now define a parallel transport along \'etale paths for the bundles in $\eB_{X_{\C_p},D}$. 

\begin{prop}
  \label{cd_t17}
Let $X$ be a smooth, projective curve over $\oQ_p$ with models $\eX_1$ and $\eX_2$ over $\oZ_p$. Then there is a third model $\eX_3$ of $X$ together with morphisms
\[
\eX_1 \xleftarrow{p_1} \eX_3 \xrightarrow{p_2} \eX_2
\]
restricting to the identity on the generic fibres (after their identification with $X$). For any divisor $D$ on $X$ we have a commutative diagram of fully faithful functors
\[
\xymatrix{
\eB_{\eX_{1\eo},D} \otimes \Q \ar[rd]^{p^*_1} \ar@/^/[rrd]^{j^*_{\eX_{1\eo}}} & & \\
 & \eB_{\eX_{3\eo},D} \otimes \Q \ar[r]^{j^*_{\eX_{3\eo}}} & \eB_{X_{\C_p},D} \; .\\
\eB_{\eX_{2\eo},D} \otimes \Q \ar[ru]_{p^*_2} \ar@/_/[rru]_{j^*_{\eX_{2\eo}}} & &
}
\]
\end{prop}

\begin{proof}
  Descend $X$ to a smooth projective curve $X_K$ over a finite extension $K$ of $\Q_p$, and $\eX_1 , \eX_2$ to models $\eX_{1 , \eo_K} , \eX_{2 , \eo_K}$ of $X_K$ over $\eo_K$. Let $\eX^*_{\eo_K}$ be the closure of the image of the morphism
\[
X_K \xrightarrow{\Delta} X_K \times_{\spec K} X_K \longrightarrow \eX_{1 , \eo_K} \times_{\spec \eo_K} \eX_{2, \eo_K}
\]
endowed with the reduced subscheme structure. Let $\eX_{3, \eo_K}$ be the normalization of $\eX^*_{\eo_K}$. Then there are natural morphisms $\eX_{3 , \eo_K} \to \eX^*_{\eo_K} \to \eX_{1 , \eo_K}$ and $\eX_{3 , \eo_K} \to \eX^*_{\eo_K} \to \eX_{2 , \eo_K}$ restricting to the identity on the generic fibres. Now the first claim follows by base change. It remains to show that for any model $\eX$ of $X$ the functor $j^*_{\eX_{\eo}} : \eB_{\eX_{\eo} , D} \otimes \Q \to \eB_{X_{\C_p} , D}$ induced by the canonical morphism $j_{\eX_{\eo}} : X_{\C_p} \to \eX_{\eo}$ is fully faithful. For bundles $\eE_1$ and $\eE_2$ on $\eX_{\eo}$ set $F = \uHom_{\eX_{\eo}} (\eE_1 , \eE_2)$. Then flat base change applied to the global sections of $F$ implies that
\[
\Hom_{\eX_{\eo}} (\eE_1 , \eE_2) \otimes \Q \xrightarrow[j^*_{\eX_{\eo}}]{} \Hom_{X_{\C_p}} (j^*_{\eX_{\eo}} \eE_1 , j^*_{\eX_{\eo}} \eE_2)
\]
is an isomorphism.
\end{proof}

For every morphism $f : \eX \to \eX'$ over $\oZ_p$ of models of $X$ restricting to the identity on generic fibres, proposition \ref{cd_t16} gives a commutative diagram:
\begin{equation}
  \label{eq:cd10}
  \vcenter{\xymatrix{
\eB_{\eX'_{\eo},D} \ar[rr]^{f^*} \ar[dr]_{\rho} & & \eB_{\eX_{\eo},D} \ar[dl]^{\rho} \\
 & \Rep_{\Pi_1 (U)} (\eo) & 
}}
\end{equation}
Next, note that there is a canonical functor:
\begin{equation}
  \label{eq:cd11}
  \Rep_{\Pi_1 (U)} (\eo) \otimes \Q \longrightarrow \Rep_{\Pi_1 (U)} (\C_p) \; .
\end{equation}
Thus we get a commutative diagram:
\begin{equation}
  \label{eq:cd12}
  \vcenter{\xymatrix{
\eB_{\eX'_{\eo},D} \otimes \Q \ar[rr]^{f^*} \ar[dr]_{\rho} & & \eB_{\eX_{\eo},D}  \otimes \Q \ar[dl]^{\rho} \\
 & \Rep_{\Pi_1 (U)} (\C_p) & 
}}
\end{equation}
Together with proposition \ref{cd_t17}, we obtain a functor
\[
\rho = \rho^{X} : \eB_{X_{\C_p},D} \longrightarrow \Rep_{\Pi_1 (U)} ( \C_p) \; .
\]
Explicitly, it is given as follows: For an object $E$ of $\eB_{X_{\C_p},D}$ we obtain the continuous functor $\rho (E) = \rho_E : \Pi_1 (U) \to \uVec_{\C_p}$ by setting on the one hand $\rho_E (x) = E_x = x^* E$ for $x \in U (\C_p) = \Ob \Pi_1 (U)$. On the other hand, for $x , x' \in U (\C_p)$ the continuous map 
\[
\rho_E = \rho_{E, x, x'} : \Mor_{\Pi_1 (U)} (x, x') \longrightarrow \Hom_{\C_p} (E_x , E_{x'})
\]
is given by
\[
\rho_E (\gamma) = \psi^{-1}_{x'} \verk (\rho_{\eE} (\gamma) \otimes_{\eo} \C_p) \verk \psi_x \; .
\]
Here we have chosen a model $\eX$ of $X$ over $\oZ_p$ and a bundle $\eE$ in $\eB_{\eX_{\eo},D}$ together with an isomorphism $\psi : E \to j^*_{\eX_{\eo}} \eE$ in $\uVec_{X_{\C_p}}$. Moreover $\psi_x$ is the fibre map:
\[
\psi_x = x^* (\psi) : E_x \silo (j^*_{\eX_{\eo}} \eE)_x = \eE_{x_{\eo}} \otimes_{\eo} \C_p  = \eE_{x_{\eo}} \otimes_{\Z} \Q \; .
\]
For a morphism $f : E \to E'$ in $\eB_{X_{\C_p},D}$ the morphism $\rho (f) = \rho_f : \rho_E \to \rho_{E'}$ is given by the family of linear maps $f_x = x^* (f) : E_x \to E'_x$ for all $x \in U (\C_p)$.

The main properties of parallel transport for bundles of class $\eB$ on $p$-adic curves are collected in the next result:

\begin{theorem}
  \label{cd_t18}
Let $X , X'$ be smooth projective curves over $\oQ_p$ and let \linebreak
$f : X \to X'$ be a morphism between them. Let $D$ and $D'$ be divisors on $X$ and $X'$.\\
{\bf a} The functor
\[
\rho : \eB_{X_{\C_p},D} \longrightarrow \Rep_{\Pi_1 (U)} (\C_p)
\]
is $\C_p$-linear, exact and commutes with tensor products, duals, internal homs and exterior powers.\\
{\bf b} Pullback of vector bundles induces an additive and exact functor \linebreak
$f^* : \eB_{X'_{\C_p},D'} \to \eB_{X_{\C_p},f^* D'}$ which commutes with tensor products, duals, internal homs and exterior powers. The following diagram is commutative:
\begin{equation}
  \label{eq:cd13}
  \vcenter{\xymatrix{
\eB_{X'_{\C_p},D'} \ar[r]^-{\rho} \ar[d]_{f^*} & \Rep_{\Pi_1 (X' \ohne D')} (\C_p) \ar[d]^{A (f)} \\
\eB_{X_{\C_p},f^*D'} \ar[r]^-{\rho} & \Rep_{\Pi_1 (X \ohne f^* D')} (\C_p)
}}
\end{equation}
In particular, for $E$ in $\eB_{X_{\C_p},f^* D'}$ we have
\begin{equation}
  \label{eq:cd14}
  \rho_{f^* E} = \rho_E \verk f_*
\end{equation}
as functors from $\Pi_1 (X \ohne f^* D')$ to $\uVec_{\C_p}$.\\
{\bf c} For every automorphism $\sigma$ of $\oQ_p$ over $\Q_p$ the following diagram commutes
\begin{equation}
  \label{eq:cd15}
\vcenter{  \xymatrix{
\eB_{X_{\C_p},D} \ar[r]^-{\rho} \ar[d]_{\sigma_*} & \Rep_{\Pi_1 (U)} (\C_p) \ar[d]^{\bC_{\sigma}} \\
\eB_{\hsigma X_{\C_p} , \hsigma D} \ar[r]^-{\rho} & \Rep_{\Pi_1 (\hsigma U)} (\C_p)
}}
\end{equation}
In particular, we have for $E$ in $\eB_{X_{\C_p},D}$ that 
\[
\rho_{\hsigma E} = \sigma_* \verk \rho_E \verk (\sigma_*)^{-1}
\]
as functors from $\Pi_1 (\hsigma U)$ to $\uVec_{\C_p}$. If $X = X_K \otimes_K \oQ_p$ and $D = D_K \otimes_K \Q_p$ for some field $\Q_p \subset K \subset \oQ_p$, so that $(\hsigma X , \hsigma D)$ is canonically identified with $(X,D)$ over $\oQ_p$ for all $\sigma \in G_K$, the functor
\[
\rho : \eB_{X_{\C_p},D} \longrightarrow \Rep_{\Pi_1 (U)} (\C_p)
\]
commutes with the left $G_K$-actions on these categories, defined by letting $\sigma$ act via $\sigma_*$ resp. via $\bC_{\sigma}$. 
\end{theorem}

\begin{rem}
  As usual, for diagrams of functors to commute means to commute up to canonical isomorphisms. 
\end{rem}

\begin{proof}
  Assertions {\bf a} and {\bf c} follow from propositions \ref{cd_t14} and \ref{cd_t15}, respectively. Assertion {\bf b} follows from proposition \ref{cd_t16} and lemma \ref{cd_t7neu}.
\end{proof}

Assume that $X$ has a smooth model $\eX$ over $\oZ_p$. Then by Theorem \ref{cd_t9}, every line bundle of degree zero on $X_{\C_p}$ lies in $\eB_{X_{\C_p}}$. Hence our functor $\rho$ induces a homomorphism
\[
\Pic^0 (X_{\C_p}) \longrightarrow \Hom_{\cont} (\pi_1 (X,x) , \C^*_p) \; .
\]
In \cite{De-We2} we show that on a certain open subgroup of $\Pic^0 (X_{\C_p})$ this homomorphism coincides with the one constructed by Tate in \cite{Ta} \S\,4 using the $p$-divisible group of the abelian scheme $\Pic^0_{\eX / \oZ_p}$ and its Cartier dual. 

The following theorem gives another relation to Tate's work \cite{Ta}. A proof is contained in \cite{De-We2}.

\begin{theorem}
  \label{cd_t25_neu}
Let $X$ be a smooth, projective curve over $\oQ_p$ with a smooth model $\eX$ over $\oZ_p$. We write $\Ext^1_{\eB_{X_{\C_p}}} (\Oh , \Oh)$ for the Yoneda groups of isomorphy classes of extensions $0 \to \Oh \to \Oh (E) \to \Oh \to 0$, where $E$ lies in $\eB_{X_{\C_p}}$. Since $\rho$ is exact, it induces a homomorphism 
\[
\rho_* : \Ext^1_{\eB_{X_{\C_p}}} (\Oh , \Oh) \longrightarrow \Ext_{\Rep_{\pi_1 (X,x)} (\C_p)} (\C_p , \C_p) \; .
\]
Then the following diagram commutes:
\[
\xymatrix{
\Ext^1_{\eB_{X_{\C_p}}} (\Oh , \Oh) \ar[r]^-{\rho_*} \ar@{=}[d] & \Ext^1_{\Rep_{\pi_1 (X,x)} (\C_p)} (\C_p , \C_p) \ar@{=}[d] \\
H^1 (X , \Oh) \otimes_{\oQ_p} \C_p \ar[r]^{\alpha} & H^1_{\et} (X , \Q_p) \otimes_{\oQ_p} \C_p
}
\]
where $\alpha$ is the Hodge--Tate map coming from the Hodge--Tate decomposition of \linebreak
$H^1_{\et} (X , \Q_p) \otimes_{\oQ_p} \C_p$. 
\end{theorem}

\begin{prop}
  \label{cd_t20}
For a fixed point $x_0 \in U (\C_p)$ the functor ``fibre in $x_0$''
\[
\eB_{X_{\C_p},D} \longrightarrow \uVec_{\C_p} \; , \; E \longmapsto E_{x_0} \; , \; f \longmapsto f_{x_0}
\]
is faithful. In particular, the evaluation map
\[
\Gamma (X_{\C_p} , E) \longrightarrow E_{x_0} \; , \; s \longmapsto s (x_0)
\]
is injective for all bundles $E$ in $\eB_{X_{\C_p},D}$.
\end{prop}

\begin{proof}
  The functor $\rho : \eB_{X_{\C_p},D} \to \Rep_{\Pi_1 (U)} (\C_p)$ is faithful because a morphism of vector bundles $f : E \to E'$ is determined by the collection of linear maps $f_x : E_x \to E'_x$ for all $\C_p$-valued points $x$ of $U_{\C_p}$ cf. \cite{EGA1} 7.2.2.1. Using lemma \ref{cd_t11} it follows that the functor ``fibre in $x_0$'' is faithful as well. In particular, the map
\[
\Gamma (X_{\C_p} , E) = \Hom_{X_{\C_p}} (\Oh , E) \longrightarrow \Hom_{\C_p} (\C_p , E_{x_0}) = E_{x_0}
\]
is injective, where $\Oh$ denotes the trivial line bundle on $X_{\C_p}$.
\end{proof}

In order to extend the preceeding results to the category $\eB^{\sharp}_{X_{\C_p},D}$ we need the following result.

\begin{prop}
  \label{cd_t21}
Consider a Galois covering $\alpha : Y \to X$ between varieties over $\oQ_p$. A (continuous) functor $W : \Pi_1 (Y) \to \Ch$ into a (topological) category $\Ch$ factors as $W = V \verk \alpha_*$ for some (continuous) functor $V : \Pi_1 (X) \to \Ch$ if and only if we have $W \verk \sigma_* = W$ for all $\sigma \in G = \Gal (Y / X)$. If $\alpha$ is only finite and \'etale but not necessarily Galois, the relation $W = V \verk \alpha_*$ already determines $V$ uniquely. 
\end{prop}

\begin{proof}
  The condition $W \verk \sigma_* = W$ is necessary for the existence of $V$ since $\alpha_* \verk \sigma_* = (\alpha \verk \sigma)_* = \alpha_*$. Now assume that we have $W \verk \sigma_* = W$ for all $\sigma$. It implies that $W_{\sigma (y)} = W_y$ for all $y \in Y (\C_p)$. Hence we may define
\[
V : X (\C_p) = \Ob \Pi_1 (X) \longrightarrow \Ob \Ch
\]
by setting $V_x = W_y$ for an arbitrary $y \in Y (\C_p)$ with $\alpha (y) = x$. We define
\[
V : \Mor_{\Pi_1 (X)} (x_1 , x_2) \longrightarrow \Mor_{\Ch} (V_{x_1} , V_{x_2})
\]
as follows. Let $y_1 \in Y (\C_p)$ be a point with $\alpha (y_1) = x_1$. For any finite \'etale map $\alpha$, the natural map
\begin{equation}
  \label{eq:cd16}
  \coprod_{\alpha (y_2) = x_2} \Mor_{\Pi_1 (Y)} (y_1 , y_2) \overset{\alpha_*}{\silo} \Mor_{\Pi_1 (X)} (x_1 , x_2)
\end{equation}
is a homeomorphism. Hence any \'etale path $\gamma$ from $x_1$ to $x_2$ has a unique lifting to an \'etale path $\gamma'$ from $y_1$ to some point $y_2$ above $x_2$. The desired relation $W = V \verk \alpha_*$ forces us to set $V (\gamma) = W (\gamma')$, a morphism from $V_{x_1} = W_{y_1}$ to $V_{x_2} = W_{y_2}$. We have to check that this is well defined i.e. independent of the choice of $y_1$. Let $y'_1$ be another point above $x_1$ and let $\sigma \in G$ be the automorphism with $\sigma y_1 = y'_1$. Then $\sigma_* (\gamma')$ is the unique path above $\gamma$ from $y'_1$ to some point $y'_2$ above $x_2$. Thus we have to show that $W (\gamma') = W (\sigma_* (\gamma'))$. But this follows from the relation $W \verk \sigma_* = W$ on morphisms. It is clear that $V$ is a functor with $W = V \verk \alpha_*$. We have also seen that this property determines $V$ uniquely. The continuity assertions are clear. 
\end{proof}

\begin{rem}
  In particular the proposition applies to representations of $\Pi_1 (Y)$ on $\C_p$-vector spaces. There is no analogous result if one only considers representations of the fundamental group $\pi_1 (Y,y)$. For example, consider a smooth surface $X$ with finite fundamental group and universal covering $\alpha : Y \to X$. Then a representation of the trivial group $\pi_1 (Y,y)$ carries no information whereas a representation $W$ of $\Pi_1 (Y)$ defines a transitive set of isomorphisms between the vector spaces $W_y$ for all $y$ in $Y (\C_p)$. 
\end{rem}

We can now define a functor $\rho : \eB^{\sharp}_{X_{\C_p},D} \to \Rep_{\Pi_1 (U)} (\C_p)$ extending the functor $\rho$ previously defined on $\eB_{X_{\C_p},D}$. Thus let $E$ be a vector bundle in $\eB^{\sharp}_{X_{\C_p},D}$. Choose a ramified Galois covering $\alpha : Y \to X$ which is \'etale over $U = X \ohne D$ such that $\alpha^* E$ lies in $\eB_{Y_{\C_p},\alpha^* D}$. By theorem \ref{cd_t18} {\bf b} we have
\[
\rho_{\alpha^* E} \verk \sigma_* = \rho_{\sigma^* (\alpha^* E)} = \rho_{\alpha^* E}
\]
for every $\sigma$ in the Galois group of $V = Y \ohne \alpha^* D$ over $U$. Using proposition \ref{cd_t21} it follows that there is a unique functor $\rho (E) = \rho_E : \Pi_1 (U) \to \uVec_{\C_p}$ such that we have
\begin{equation}
  \label{eq:cd17}
  \rho_{\alpha^* E} = \rho_E \verk \alpha_* \; .
\end{equation}
This functor is continuous. (In order to apply proposition \ref{cd_t21}, we view canonical isomorphisms such as $\sigma^* (\alpha^* E) = (\alpha \sigma)^* E$ as identifications.)\\
In particular, we have $\rho_E (x) = E_x$ for all $x \in U (\C_p)$. For an \'etale path $\gamma$ from $x_1$ to $x_2$ in $U$ we have
\begin{equation}
  \label{eq:cd18}
  \rho_E (\gamma) = \rho_{\alpha^* E} (\gamma') : E_{x_1} = (\alpha^* E)_{y_1} \to (\alpha^* E)_{y_2} = E_{x_2} \; .
\end{equation}
Here $y_1 \in V (\C_p)$ lies above $x_1$ and $\gamma'$ is the unique path in $V$ with $\alpha_* \gamma' = \gamma$ from $y_1$ to a point $y_2$ above $x_2$. For a morphism $f : E \to E'$ of vector bundles in $\eB^{\sharp}_{X_{\C_p},D}$ the morphism $\rho (f) = \rho_f : \rho_E \to \rho_{E'}$ is defined to be the family of linear maps $f_x : E_x \to E'_x$ for all $x \in U (\C_p)$.

\begin{prop}
\label{cd_t22}
The preceeding constructions give a well defined functor $\rho : \eB^{\sharp}_{X_{\C_p},D} \to \Rep_{\Pi_1 (U)} (\C_p)$ which extends the previously defined functor $\rho$ on $\eB_{X_{\C_p},D}$. 
\end{prop}

\begin{proof}
  We first have to show that the definition of $\rho_E$ is independent of $\alpha$. If we are given ramified Galois coverings of smooth projective curves $\alpha_1 : Y_1 \to X$ and $\alpha_2 : Y_2 \to X$ which are \'etale over $U$ there is a third one $\alpha_3 : Y_3 \to X$ covering $\alpha_1$ and $\alpha_2$ i.e. $\alpha_3 = \alpha_i \verk \pi_i$ for morphisms $\pi_i : Y_3 \to Y_i$ where $i = 1 ,2$. Now assume that $\alpha^*_i E \in \eB_{Y_{i \C_p},\alpha^*_i D}$. By the above we have
\[
\rho_{\alpha^*_i E} = \rho_i \verk \alpha_{i*}
\]
for functors $\rho_i : \Pi_1 (U) \to \uVec_{\C_p}$ where $i = 1,2$. We have to show that $\rho_1 = \rho_2$. By theorem \ref{cd_t18}, {\bf b} we find for $i = 1,2$ that
\begin{eqnarray*}
  \rho_{\alpha^*_3 E} & = & \rho_{\pi^*_i (\alpha^*_i E)} = \rho_{\alpha^*_i E} \verk \pi_{i*} = \rho_i \verk \alpha_{i*} \verk \pi_{i*} \\
 & = & \rho_i \verk \alpha_{3*} \; .
\end{eqnarray*}
The uniqueness assertion of proposition \ref{cd_t21} now implies that $\rho_1 = \rho_2$. 

Next we have to check that for a morphism $f : E \to E'$ in $\eB^{\sharp}_{X_{\C_p},D}$ the family of maps $f_x : E_x \to E'_x$ defines a morphism in $\Rep_{\Pi_1 (U)} (\C_p)$. We may assume that both $\alpha^* E$ and $\alpha^* E'$ lie in $\eB_{Y_{\C_p},\alpha^*D}$. Then $\rho_{\alpha^* f}$, i.e. the family of maps $(\alpha^* f)_y : (\alpha^* E)_y \to (\alpha^* E')_y$, defines a morphism in $\Rep_{\Pi_1 (V)} (\C_p)$. 
Using (\ref{eq:cd18}) we see that
\[
\xymatrix{
E_{x_1} \ar[r]^{f_{x_1}} \ar[d]_{\rho_E (\gamma)} & E'_{x_1} \ar[d]^{\rho_{E'} (\gamma)} \\
E_{x_2} \ar[r]^{f_{x_2}} & E'_{x_2}
}
\]
commutes for every $\gamma$, as desired. It is clear that $\rho$ is a functor and that it extends $\rho : \eB_{X_{\C_p},D} \to \Rep_{\Pi_1 (U)} (\C_p)$. 
\end{proof}


\begin{theorem}
  \label{cd_t23}
Assertions {\bf a}, {\bf b} and {\bf c} of theorem \ref{cd_t18} hold for $\eB^{\sharp}_{X_{\C_p},D}$ instead of $\eB_{X_{\C_p},D}$ as well. For any point $x_0 \in U (\C_p)$ the fibre functor
\[
\eB^{\sharp}_{X_{\C_p},D} \longrightarrow \uVec_{\C_p} \; , \; E \longmapsto E_{x_0} \; , \; f \longmapsto f_{x_0}
\]
is faithful.
\end{theorem}

\begin{proof}
  {\bf a} Exactness of $\rho$ is clear from its definition. For two vector bundles $E_1$ and $E_2$ in $\eB^{\sharp}_{X_{\C_p},D}$ choose a ramified Galois covering $\alpha : Y \to X$, \'etale over $X \ohne D$ such that $\alpha^* E_1$ and $\alpha^* E_2$ are in $\eB_{X_{\C_p},D}$. Then we have
\begin{eqnarray*}
  \rho_{\alpha^* (E_1 \otimes E_2)} & = & \rho_{\alpha^* E_1 \otimes \alpha^* E_2} = \rho_{\alpha^* E_1} \otimes \rho_{\alpha^* E_2} = (\rho_{E_1} \verk \alpha_*) \otimes (\rho_{E_2} \verk \alpha_*) \\
& = & (\rho_{E_1} \otimes \rho_{E_2}) \verk \alpha_*
\end{eqnarray*}
by theorem \ref{cd_t18} {\bf a}. It follows that we have $\rho_{E_1 \otimes E_2} = \rho_{E_1} \otimes \rho_{E_2}$. It is clear that we also have $\rho_{f_1 \otimes f_2} = \rho_{f_1} \otimes \rho_{f_2}$ for morphisms $f_1 , f_2$ of vector bundles. Hence $\rho$ commutes with $\otimes$-products and similarly with direct sums, duals, internal homs and exterior powers.

Let $f : X \to X'$ be a morphism of smooth projective curves over $\oQ_p$. On objects, we have to show that $\rho_{f^* E} = \rho_E \verk f_*$ for all $E$ in $\eB^{\sharp}_{X'_{\C_p}}$. Thus let $\alpha' : Y' \to X'$ be a ramified Galois covering \'etale over $X' \ohne D'$ with $\alpha'^* E$ in $\eB_{Y'_{\C_p},\alpha'^* D'}$. With notations as in the proof of proposition \ref{cd_t6} we see that $\alpha^* f^* E = g^* \alpha'^* E$ lies in $\eB_{Y_{\C_p},\alpha^* f^* D'}$. Moreover
\[
\rho_{\alpha^* f^* E} = \rho_{g^* \alpha'^* E} = \rho_{\alpha'^* E} \verk g_* = \rho_E \verk \alpha'_* \verk g_* = \rho_E \verk f_* \verk \alpha_*
\]
by theorem \ref{cd_t18} {\bf a} and the definition of $\rho_E$. Now on the other hand, $\alpha : Y \to X$ is a ramified Galois covering, \'etale over $X \ohne f^* D'$. Hence $f^* E$ lies in $\eB^{\sharp}_{X_{\C_p} , f^* D'}$ and $\rho_{f^* E}$ is by definition the unique functor with $\rho_{\alpha^* f^* E} = \rho_{f^* E} \verk \alpha_*$. It follows that we have $\rho_{f^* E} = \rho_E \verk f_*$. It is immediate from the definitions that diagram (\ref{eq:cd13}) for $\eB^{\sharp}$ commutes on the level of morphisms. The proof that $\rho$ behaves functorially with respect to automorphisms is deduced similarly from theorem \ref{cd_t18}, {\bf c}. The last assertion is proved in the same way as proposition \ref{cd_t20}.
\end{proof}

\begin{rem}
  It is known that the fibre functor in a point is faithful on the category of stable bundles of degree zero on a compact Riemann surface cf. \cite{Sesh} Ch. 1, IV. By an induction on the length of the Jordan--H\"older filtration one gets faithfulness also on the category of semistable bundles of degree zero. The analogous assertion therefore holds on smooth projective curves over fields that can be embedded into $\C$, e.g. over $\C_p$. Together with theorem \ref{cd_t13neu} one thus gets another proof of theorem \ref{cd_t23} {\bf b}. 
\end{rem}

We will now explain how to glue the representations $\rho_E$ attached to a vector bundle $E$ on $X_{\C_p}$ which belongs to $\eB^{\sharp}_{X_{\C_p} , D}$ for several divisors $D$. For this we need the following Seifert--van Kampen theorem for \'etale groupoids:

\begin{prop}
  \label{cd_t30_neu}
Given open subschemes $U_1$ and $U_2$ of a curve $X$, let $i_1 : U_1 \cap U_2 \hookrightarrow U_1 , i_2 : U_1 \cap U_2 \hookrightarrow U_2$ and $j_1 : U_1 \hookrightarrow U_2 \cup U_2 , j_2 : U_2 \hookrightarrow U_1 \cup U_2$ be the corresponding immersions and consider the commutative diagram of fundamental groupoids
\[
\xymatrix{
\Pi_1 (U_1 \cap U_2) \ar[r]^{i_{1*}} \ar[d]_{i_{2*}} & \Pi_1 (U_1) \ar[d]^{j_{1*}} \\
\Pi_1 (U_2) \ar[r]^{j_{2*}} & \Pi_1 (U_1 \cup U_2)
}
\]
Then for any Hausdorff topological category $\Ch$ and continuous functors $\rho_1 : \Pi_1 (U_1) \to \Ch$ and $\rho_2 : \Pi_1 (U_2) \to \Ch$ such that $\rho_1 \verk i_{1*} = \rho_2 \verk i_{2*}$ there is a unique continuous functor $\rho : \Pi_1 (U_1 \cup U_2) \to \Ch$ such that $\rho \verk j_{1*} = \rho_1$ and $\rho \verk j_{2*} = \rho_2$.
\end{prop}

\begin{proof}
  We may assume that $U_1$ and $U_2$ are nonempty. Let $\gamma : x_1 \to x_2$ be an \'etale path in $U_1 \cup U_2$ with $x_1 \in U_1$ and $x_2 \in U_2$. Choose a point $x'$ in $U_1 \cap U_2$. Then since $U_1 \cup U_2$ is connected we may write $\gamma$ as the composition of a path $\gamma_1 : x_1 \to x'$ with a path $\gamma_2 : x' \to x_2$ in $U_1 \cup U_2$. The homomorphisms $\pi_1 (U_{\nu} , x_{\nu}) \to \pi_1 (U_1 \cup U_2 , x_{\nu})$ for $\nu = 1,2$ are known to be surjective since $X$ is a curve. We deduce that there are paths $\tgamma_1 : x_1 \to x'$ in $U_1$ and $\tgamma_2 : x' \to x_2$ in $U_2$ such that $j_{\nu *} (\tgamma_{\nu}) = \gamma_{\nu}$. Hence if $\rho$ exists, we have $\rho (\gamma) = \rho (\gamma_2 \cdot \gamma_1) = \rho (\gamma_2) \rho (\gamma_1) = \rho_2 (\tgamma_2) \rho_1 (\tgamma_1)$ and similarly for paths from $x_2$ to $x_1$. For paths $\gamma$ in $X$ whose endpoints are both in $U_{\nu}$ we have $\rho (\gamma) = \rho_{\nu} (\tgamma)$ where $j_{\nu *} (\tgamma) = \gamma$. Hence the functor $\rho$ is uniquely determined. As for existence, it is clear how to define $\rho$ on objects and it remains to check that $\rho$ given on morphisms by the above formulas is well defined. This follows from the Seifert--van Kampen theorem for the \'etale fundamental group, cf. \cite{SGA1} IX, corollare 5.6. There is a subtlety here. The pushout property holds only in the category of profinite groups. But the maps $i_{\nu *} : \pi_1 (U_1 \cap U_2 , x') \to \pi_1 (U_{\nu} , x')$ are surjective and hence the maps $\rho_{\nu}$ on $\pi_1 (U_{\nu} , x')$ have the same images for $\nu = 1,2$. This common image is the quotient of a profinite group by a closed subgroup and hence profinite. Here we used the Hausdorff assumption on (the spaces of morphisms of) $\Ch$.
\end{proof}

\begin{prop}
  \label{cd_t32_neu}
Let $D_1$ and $D_2$ be divisors on $X$ and set $U_1 = X \ohne D_1$ and $U_2 = X \ohne D_2$. For a vector bundle $E$ on $X_{\C_p}$ let $\rho^1_E$ and $\rho^2_E$ be the continuous representations of $\Pi_1 (U_1)$ resp. $\Pi_1 (U_2)$ on $\C_p$-vector spaces constructed before. Then there is a unique continuous representation $\rho_E$ of $\Pi_1 (U_1 \cup U_2)$ which induces $\rho^{\nu}_E$ on $\Pi_1 (U_{\nu})$ for $\nu = 1,2$. For the induced functor where $U = U_1 \cup U_2$
\[
\rho : \eB^{\sharp}_{X_{\C_p} , D_1} \cap \eB^{\sharp}_{X_{\C_p} , D_2} \longrightarrow \Rep_{\Pi_1 (U)} (\C_p)
\]
the analogue of theorem \ref{cd_t23} holds.
\end{prop}

{\bf Variant} For $\eE$ in $\eB_{\eX_{\eo} , D_{\nu}}$ for $\nu = 1,2$ we obtain a well defined representation of $\Pi_1 (U)$ on free $\eo$-modules of finite rank.

\begin{proof}
  On objects $\rho_E$ is defined by $\rho_E (x) = E_x$ as before. The assertions are a formal consequence of proposition \ref{cd_t30_neu} and theorem \ref{cd_t23}.
\end{proof}

Let $\eB^s_{X_{\C_p}}$ be the category of vector bundles on $X_{\C_p}$ with strongly semistable reduction, as defined in the introduction.

\begin{theorem}
  \label{cd_t36_neu}
We have $\eB^s_{X_{\C_p}} = \bigcup_D \eB_{X_{\C_p}, D}$. Every vector bundle $E$ in $\eB^s_{X_{\C_p}}$ lies both in $\eB_{X_{\C_p} , D}$ and $\eB_{X_{\C_p} , \tD}$ for suitable divisors $D$ and $\tD$ with disjoint support. There is a unique representation $\rho_E$ of $\Pi_1 (X)$ on finite dimensional $\C_p$-vector spaces such that $\rho_E (x) = E_x$ for all $x \in X (\C_p)$ and such that $\rho_E$ is compatible with the representations $\rho_E$ of $\Pi_1 (X \ohne D)$ constructed earlier for those $D$ with $E$ in $\eB_{X_{\C_p} ,D}$. As before one obtains an exact additive functor $\rho : \eB^s_{X_{\C_p}} \to \Rep_{\Pi_1 (X)} (\C_p)$ which commutes with tensor products, duals, internal homs and exterior powers. Moreover, it behaves functorially with respect to morphisms of curves over $\oQ_p$ and automorphisms of $\oQ_p$ over $\Q_p$. For any point $x_0 \in X (\C_p)$, the fibre functor
\[
\eB^s_{X_{\C_p}} \longrightarrow \uVec_{\C_p} \; , \; E \longmapsto E_{x_0} \; , \; f \longmapsto f_{x_0}
\]
is faithful.
\end{theorem}

\begin{proof}
  This follows from theorem \ref{cd_t16_x} and \ref{cd_t23} together with proposition \ref{cd_t32_neu}. 
\end{proof}

\begin{rem}
Arguing as in the proofs of propositions \ref{cd_t6} and \ref{cd_t22} this result implies the theorem in the introduction.
\end{rem}

Let $\Rep^{\infty}_{\pi_1 (X , x)} (\C_p)$ be the full subcategory of $\Rep_{\pi_1 (X , x)} (\C_p)$ of those representations $\lambda : \pi_1 (X , x) \to \GL (V)$ which are continuous if $V$ is given the discrete topology. Equivalently $\lambda$ has to factor over a finite quotient of $\pi_1 (X , x)$.

\begin{prop}
  \label{cd_t24}
The category $\Rep^{\infty}_{\pi_1 (X , x)} (\C_p)$ is contained in the essential image of $\rho : \eB^{\sharp}_{X_{\C_p}} \to \Rep_{\pi_1 (X , x)} (\C_p)$.
\end{prop}

\begin{proof}
  Let $\lambda$ be a representation as above. Let $G$ be the image of $\lambda$ in $\GL (V)$. Let $\alpha : Y \to X$ be a Galois extension of $X$ with group $G$ such that $\alpha_* : \pi_1 (Y , y) \to \pi_1 (X , x)$ with $y \in Y (\C_p)$ above $x$ induces an isomorphism $\pi_1 (Y , y) \cong \ker \lambda$. Define a vector bundle $E$ on $X$ by setting $E = Y \times^G \bV$ where $\bV$ is the affine space over $\oQ_p$ attached to $V$. Then $E$ lies in $\eB^{\sharp}_{X_{\C_p}}$ because $\alpha^* E$ is a trivial bundle on $Y$ and hence lies in $\eB_{Y_{\C_p}}$. On $S$-valued points of $Y$ a trivialization
\[
\psi : \bV_Y = Y \times \bV \silo \alpha^* E
\]
is described by mapping $(y , v)$ to the pair $([y,v] , y)$ in $\alpha^* E$. Here $[y, v] \in E_{\alpha (y)}$ is the class of $(y,v) \mod G$. We can now calculate $\rho_E$. For $\gamma \in \pi_1 (X , x)$ there is a unique \'etale path $\gamma'$ in $Y$ from $y$ to $\sigma y$ for a uniquely determined $\sigma \in G$. The commutative diagram
\[
\xymatrix{
E_x \ar@{=}[r] \ar[d]_{\rho_E (\gamma)} & (\alpha^* E)_y \ar[d]^{\rho_{\alpha^* E} (\gamma')} & (\bV_Y)_y \ar[l]_-{\overset{\psi_y}{\sim}} \ar[d]^{\rho_{\bV_Y} (\gamma')} \ar@{=}[r] & V \ar@{=}[d] \\
E_x \ar@{=}[r] & (\alpha^* E)_{\sigma y} & (\bV_Y)_{\sigma y} \ar[l]_-{\overset{\psi_{y^{\sigma}}}{\sim}} \ar@{=}[r] & V
}
\]
shows that if we identify $E_x$ with $V$ via $\psi_y$ the automorphism $\rho_E (\gamma) : V \to V$ is given by $\rho_E (\gamma) = \psi^{-1}_y \verk \psi_{y^{\sigma}}$. Thus we have
\[
\rho_E (\gamma) (v) = \psi^{-1}_y [y^{\sigma} , v] = \psi^{-1}_y [y , \sigma v] = \sigma v = \lambda (\gamma) v \; .
\]
Hence we have $\rho_E = \lambda$ as representations of $\pi_1 (X , x)$ on $V \overset{\psi_y}{\cong} E_x$.
\end{proof}

\begin{prop}
  \label{cd_t25}
Let $\alpha : Y \to X$ be a finite \'etale covering of smooth projective curves over $\oQ_p$ and let $E$ be a vector bundle in $\eB^{\sharp}_{Y_{\C_p}}$. Then
\[
\rho_{\alpha_* E} : \pi_1 (X , x) \to \GL ((\pi_* E)_x)
\]
is the representation obtained from $\rho_E : \pi_1 (Y , y) \to \GL (E_y)$ by induction via the inclusion $\alpha_* : \pi_1 (Y , y) \hookrightarrow \pi_1 (X , x)$. Here $y$ is any point in $Y (\C_p)$ above $x \in X (\C_p)$.
\end{prop}

\begin{proof}
  Under the natural injection $E \hookrightarrow \alpha^* \alpha_* E$ we may view $E_y$ as a subspace of $(\alpha^* \alpha_* E)_y = (\alpha_* E)_x$. From theorem \ref{cd_t23} {\bf a} we get that
\[
\rho_{\alpha_* E} \verk \alpha_* = \rho_{\alpha^* \alpha_* E}
\]
as representations of $\pi_1 (Y , y)$ on $\GL ((\alpha_* E)_x)$. Thus $\rho_{\alpha^* \alpha_* E}$ is the restriction of $\rho_{\alpha_* E}$ to the subgroup (via $\alpha_*$) $\pi_1 (Y , y)$ of $\pi_1 (X , x)$. Since $\rho_E \hookrightarrow \rho_{\alpha^* \alpha_* E}$ by the exactness of $\rho$, it follows that the restriction of $\rho_{\alpha_* E}$ to $\pi_1 (Y , y)$ leaves the subspace $E_y$ of $(\alpha_* E)_x$ invariant and gives the representation $\rho_E$ there. It remains to show that $(\alpha_* E)_x$ is the direct sum of the translates $\rho_{\alpha_* E} ([\gamma]) E_y$ for $[\gamma]$ running over the cosets of $\alpha_* \pi_1 (Y , y)$ in $\pi_1 (X , x)$. There is a bijection
\[
\pi_1 (X , x) / \alpha_* \pi_1 (Y , y) \silo \{ y' \in Y (\C_p) \tei \alpha (y') = x \}
\]
given by mapping $[\gamma]$ to the ``endpoint'' of the unique lifting of $\gamma$ to an \'etale path in $Y$ starting at $y$, cf. (\ref{eq:cd16}). Together with the natural isomorphism
\[
(\alpha_* E)_x = \bigoplus_{\alpha (y') = x} E_{y'}
\]
the assertion follows. Namely, we have:
\[
\rho_{\alpha_* E} (\gamma) E_y = \rho_{\alpha^* \alpha_* E} (\gamma') E_y = \rho_E (\gamma') E_y = E_{y'} \; .
\]
\end{proof}

We conclude this section with some general observations on the structure of representations.

Consider a continuous representation $\rho : G \to \GL_r (\eo)$ of a profinite group $G$. Then $\rho_1 : G \to \GL_r (\eo_1)$, the reduction $\mod p$ of $\rho$ has finite image since $G$ is compact and $\eo_1$ is discrete. Hence the image of $\rho_1$ is contained in $\GL_r (\eo_K / p \eo_K)$ for some finite extension $K$ of $\Q_p$. Let $\ep$ be the prime ideal of $\eo_K$ and consider the reductions $\rho_{(n)} : G \to \GL_r (\eo_{(n)})$ of $\rho \mod \ep^n$ where we have set $\eo_{(n)} = \eo / \ep^n \eo$. By construction $\rho_{(1)}$ factors:
\[
\rho_{(1)} : G \to \GL_r (\eo_K / \ep) \subset \GL_r (\eo_{(1)}) \; .
\]
Extending scalars to $k = \eo / \emm$ the modular representation $G \to \GL_r (\eo_K / \ep)$ becomes $\rho_k$, the reduction of $\rho \mod \emm$. For every $n \ge 1$ the image of $\rho_{(n)}$ is finite. Hence $G_n = \Ker \rho_{(n)}$ is an open normal subgroup of $G$. Let \\
$\orho_{(n)} : G / G_n \to \GL_r (\eo_{(n)})$ be the induced representation. We have a commutative diagram
\[
\xymatrix{
1 \ar[r] & G_n / G_{n+1} \ar@{^{(}->}[d]^{\lambda_n} \ar[r] & G / G_{n+1} \ar@{^{(}->}[d]^{\orho_{(n+1)}} \ar[r]  & G / G_n \ar@{^{(}->}[d]^{\orho_{(n)}} \ar[r] & 1 \\
1 \ar[r] & M_r (\ep^n \eo_{(n+1)}) \ar[r]^f & \GL_r (\eo_{(n+1)}) \ar[r] & \GL_r (\eo_{(n)}) \ar[r] & 1 
}
\]
Here $f$ is the homomorphism $f (A) = 1 + A$ and $\lambda_n$ is induced by $\orho_{(n+1)}$. Since $G_n / G_{n+1}$ is finite and $M_r (\ep^n \eo_{(n+1)})$ abelian and annihilated by $p$ it follows that $G_n / G_{n+1}$ is isomorphic to $(\Z / p)^{\delta_n}$ for some integer $\delta_n \ge 0$.

Thus $\rho$ is built up from the modular representation $\orho_{(1)}$ of the finite group $G / G_1$ in $\GL_1 (\eo_K / \ep) \subset \GL_r (k)$ via successive extensions by representations of elementary abelian $p$-groups. It is instructive to compare this fact with the proof of theorem \ref{cd_t16_neu}: The way a bundle $\eE$ in $\eB_{\eX , D}$ is built up from $\eE_k$ is similar to the way a continuous representation $\rho$ is built up from its residual representation $\rho_k$. 

Let us call a representation on a free $\eo_{(n)}$-module $M$  ``irreducible'' if every invariant free and cofree submodule of $M$ is either trivial or equal to $M$.

\begin{prop}
  If $\rho_{(n)}$ is ``irreducible'' for some $n \ge 1$, e.g. if $\rho_k$ is irreducible then the representation $\rho_{\C_p} : G \to \GL_r (\C_p)$ is irreducible. 
\end{prop}

\begin{proof}
  Let $V \subset \C^r_p$ be a $\rho_{\C_p}$-invariant subspace of dimension $s \neq 0,r$. Then $\Gamma = V \cap \eo^r$ is a $\rho$-invariant $\eo$-submodule of $\eo^r$ for which $\eo^r / \Gamma$ is $\eo$-torsionfree. Since $\eo^r / \Gamma$ is finitely generated it follows from \cite{B} Lemma 3.9 that $\eo^r / \Gamma$ is a free $\eo$-module of rank $t$ say. Hence we get an exact sequence:
\[
0 \longrightarrow \Gamma \longrightarrow \eo^r \longrightarrow \eo^t \longrightarrow 0 \; .
\]
By an induction on $t$ it follows that $\Gamma$ is a free $\eo$-module. Because of $V = \Gamma \otimes \C_p$ the rank of $\Gamma$ is equal to $s$. Hence $\Gamma_{(n)} = \Gamma \otimes \eo_{(n)}$ is a free and cofree $\eo_{(n)}$-module of rank $s$ and therefore $\rho_{(n)}$ is ``reducible''. Note that the rank of a free $\eo_{(n)}$-module $\Gamma_{(n)}$ is well defined because it equals the dimension of $\Gamma_{(n)} \otimes k$ over $k$.
\end{proof}

{\bf Example} Let $\rho : G \to \GL_2 (\eo)$ be a representation for which the image of $\rho_2 : G \to \GL_2 (\eo / p^2 \eo)$ contains the two commuting matrices $\left( 
  \begin{smallmatrix}
    1 & p \\ 0 & 1
  \end{smallmatrix} \right)$ and $\left( 
  \begin{smallmatrix}
    1 & 0 \\ p & 1
  \end{smallmatrix} \right)$. Then $\rho_2$ is ``irreducible'' and hence $\rho_{\C_p}$ is irreducible as well. For example $\rho_1$ could be trivial and we could have $G_1 / G_2 = G / G_2 \cong (\Z / p)^2$ with $\rho_2$ given by $\rho_2 (i,j) = \left( 
  \begin{smallmatrix}
    1 & ip \\ jp & 1
  \end{smallmatrix} \right)$.


\begin{minipage}[t]{6cm}
Mathematisches Institut\\
Einsteinstr. 62\\
48149 M\"unster, Germany\\
deninger@math.uni-muenster.de
\end{minipage} \hspace*{\fill}
\begin{minipage}[t]{7cm}
Institut f\"ur Algebra und Zahlentheorie\\
Pfaffenwaldring 57\\
70569 Stuttgart, Germany\\
Annette.Werner@mathematik.uni-stuttgart.de
\end{minipage}
\end{document}